\documentclass[11pt]{article}
 \usepackage{benstyle}
 
\usepackage{tikz}
\usetikzlibrary{shapes,arrows}
\usetikzlibrary{decorations.pathmorphing}

\usepackage{bm}
\usepackage{marginnote}

\usepackage[font=small]{caption}

\title{
Approximating smooth, multivariate functions on irregular domains
}
\author{Ben Adcock\footnote{Department of Mathematics,	Simon Fraser University, 8888 University Drive, Burnaby, BC V5A 1S6, Canada (\texttt{ben\_adcock@sfu.ca}, \texttt{http://www.benadcock.ca})}
\and 
Daan Huybrechs\footnote{Department of Computer Science, KU Leuven, Celestijnenlaan 200A, BE-3001 Leuven, Belgium (\texttt{daan.huybrechs@cs.kuleuven.be}, \texttt{http://people.cs.kuleuven.be/\textasciitilde daan.huybrechs/})}
}

\begin{document}

\maketitle

\begin{abstract}
In this paper, we introduce a method known as \textit{polynomial frame approximation} for approximating smooth, multivariate functions defined on irregular domains in $d$ dimensions, where $d$ can be arbitrary.  This method is simple, and relies only on orthogonal polynomials on a bounding tensor-product domain.  In particular, the domain of the function need not be known in advance.  When restricted to a subdomain, an orthonormal basis is no longer a basis, but a frame.  Numerical computations with frames present potential difficulties, due to the near-linear dependence of the truncated approximation system.  Nevertheless, well-conditioned approximations can be obtained via regularization, for instance, truncated singular value decompositions.  We comprehensively analyze such approximations in this paper, providing error estimates for functions with both classical and mixed Sobolev regularity, with the latter being particularly suitable for higher-dimensional problems.  We also analyze the sample complexity of the approximation for sample points chosen randomly according to a probability measure, providing estimates in terms of the corresponding \textit{Nikolskii inequality} for the domain.  In particular, we show that the sample complexity for points drawn from the uniform measure is quadratic (up to a log factor) in the dimension of the polynomial space, independently of $d$, for a large class of nontrivial domains.  This extends a well-known result for polynomial approximation in hypercubes.
\end{abstract}

\vspace{0.5pc} \noindent

\section{Introduction}\label{s:introduction}

Many problems in computational science call for the approximation of smooth, multivariate functions.  This problem is often challenging, due to the \textit{curse of dimensionality}.  Yet significant strides have been made over the last several decades towards its mitigation, typically by assuming some anisotropic behaviour of the function being approximated.  Approaches such as sparse grids \cite{bunggrieb} have enjoyed substantial success in the numerical solution of high-dimensional PDEs, and more recently techniques based on computing multivariate polynomial approximations -- often referred to as generalized polynomial chaos expansions \cite{XiuKarniadakisPC} -- have begun to be widely used for problems in Uncertainty Quantification (UQ) (see  \cite{BASBCWMatheon,ChkifaEtAl,ChkifaEtAlBreaking,CohenDeVoreApproxPDEs,DavenportEtAlLeastSquares,DoostanOwhadiSparse,NarayanZhouCCP,KarniadakisUQCS} and references therein).

The majority of algorithms for high-dimensional approximation assume the underlying function $f$ is defined over a tensor-product domain.  The key benefit of doing so is simplicity.  Indeed, the orthogonal polynomials on a tensor-product domain with respect to a tensor-product measure are precisely tensor products of the corresponding one-dimensional orthogonal polynomials.  Yet there are many practical instances where the domain of interest is not of tensor-product type.  One example is surrogate model construction in UQ.  In practice, it is often the case that the random variables are correlated \cite{WitteveenIaccarinoSimplex}, which leads to an irregular domain.  Alternatively or in addition, the given forward model may not be well-defined over the whole of the assumed tensor-product domain, or may produce values in certain regions that are known to be unphysical  (e.g.\ negative pressures).  This in effect leads to failed evaluations, resulting once more in an irregular domain \cite{SargsyanEtAlDimUQ}.  Similarly, in model order reduction, techniques such as active subspaces \cite{ConstantineBook} lead to approximation problems over irregular domains.  For example, when a function defined on a high-dimensional hypercube is projected to a function of a reduced set of parameters, the resulting domain (the projection of the hypercube) is generally polyhedral, a so-called zonotope \cite{StinsonZonotop}.  Finally, many applications in UQ also involve forward models which are piecewise smooth (see \cite{SargsyanEtAlUQDisc,SargsyanEtAlDimUQ} and references therein).  Unless such discontinuities happen to be aligned along coordinate axes, this results in an approximation problem involving two or more smooth functions defined over irregular domains.

With this issue in mind, the purpose of this paper is to present a systematic study of a simple but effective technique for approximating high-dimensional functions defined on irregular domains. It is based on using tensor-product orthogonal polynomials on a bounding box, and is referred to as \textit{polynomial frame approximation}.  The approach corresponds to approximation in a frame, rather than a basis, since there are potentially many ways the unknown function on the irregular domain can be represented in a basis on the bounding box. Our main results demonstrate that this procedure achieves (to a significant degree) the four primary criteria for a numerical approximation scheme: namely, \textit{simplicity}, \textit{accuracy}, \textit{stability} and \textit{efficiency}.  We elaborate on the meaning of these terms in the next section, however we note in passing that simplicity means that the same procedure can be applied to {a broad class of irregular domains}.  In particular, no costly parametrization of the domain or its boundary (a potentially infeasible task in high dimensions) is required to construct the approximation. {Instead, we will make the less restrictive assumption that samples can be randomly drawn from a certain measure on $\Omega$, related to the orthogonality measure on the bounding box.  Typically, this is taken as the uniform measure on $\Omega$.  Note that our focus in this paper is neither on the best choice of measure, nor on the question of how to sample from a given measure.  These are challenging, and potentially highly domain-dependent, issues, whereas in this work we strive for generality.  We return briefly to this question in \S \ref{s:conclusion}.}

The main contribution of this paper is the rigorous analysis of polynomial frame approximations.  Central to this is the notion of \text{frames} of Hilbert spaces, as opposed to more conventional orthogonal bases.  We stress at this point that our technique does not attempt to orthogonalize a basis.  Instead, it relies on the particular properties of frames to achieve accurate and stable approximations.  A key facet of frame approximations (not just of polynomial type) is that they lead to highly ill-conditioned linear systems of equations.  However, by using regularization we are able to obtain a mapping from the sample points to the polynomial space that is both well-conditioned and accurate.  We also determine approximation rates and sample complexity estimates that scale well with the underlying dimension, thus {(on the proviso that samples can be drawn efficiently from the desired measure -- see above)} mitigating the curse of dimensionality to a significant extent.

Before we proceed further, it is worth noting that polynomial frame approximation, and variations thereof, are in essence already used in many of the aforementioned applications.  Indeed, any approach to surrogate model construction in UQ which computes a generalized polynomial chaos expansion from function evaluations which are limited (due to the particular problem at hand) to a non-tensorial subdomain is equivalent to polynomial frame approximation.  See \S \ref{ss:related} for further details.  However, a thorough analysis of the accuracy, stability and efficiency of such approximations -- in particular, exploiting the connections to frame theory as we do in this paper -- is, to the best of our knowledge, lacking.  
Besides providing the first clear theoretical explanation for why these algorithms work in practical setting of irregular domains, we also expect the results of this paper to shed light on ways in which to improve them.  For example, the problem of designing better sampling sets for irregular domains -- a topic of significant practical interest.

\section{Overview of the paper}
We commence with a short overview of the paper.

\subsection{Polynomial frame approximations}\label{ss:polyframe_intro}
This paper concerns the approximation of a smooth multivariate function $f : \Omega \rightarrow \bbC$ defined over a non-tensor product domain $\Omega \subset \bbR^d$.  The approximation is based on four key steps:
\begin{enumerate}
\item[(i)] Choose a tensor-product domain $D$ such that $\Omega \subseteq D$.
\item[(ii)] Choose a tensor-product probability measure $\nu$ on $D$, a tensor-product orthonormal basis $\{ \psi_{\bm{n}} \}_{\bm{n}}$ of $L^2(D,\nu)$ and a finite index set $\Lambda$ with $| \Lambda | = N$.
\item[(iii)] Take $M$ samples of $f$ of the form $f(\bm{y}_1),\ldots,f(\bm{y}_M)$ where $\Upsilon = \{\bm{y}_1,\ldots,\bm{y}_M \} \subset \Omega$.
\item[(iv)] Compute an approximation to $f$ of the form $f_{\Upsilon,\Lambda} = \sum_{\bm{n} \in \Lambda} c_{\bm{n}} \phi_{\bm{n}}$, where $\phi_{\bm{n}} = \psi_{\bm{n}} |_{\Omega}$.
\end{enumerate}
This immediately raises a number of questions, which are now discussed:

\pbk
\textit{1.\ How to compute the approximation.}  There are several options for doing this, including interpolation if $M = | \Lambda | = N$, sparse regularization (i.e.\ compressed sensing) if $M < N$ and least-squares fitting if $M > N$.  We consider the latter.  Interpolation requires good choices of nodes $\bm{y}_1,\ldots,\bm{y}_N$ so as to maintain small Lebesgue constants, and it is unclear how to design such nodes for general irregular domains.  Compressed sensing is an interesting option, however beyond the scope of this paper (see \S \ref{s:conclusion} for some further discussion).  Least-squares fitting, on the other hand, is a popular tool for high-dimensional approximation on tensor-product domains \cite{ChkifaEtAl,MiglioratiCohenOptimal,MiglioratiThesis,MiglioratiEtAlFoCM,NarayanJakemanZhouChristoffelLS,NarayanZhouCCP,ZhouNarayanXiu}, and has the twin benefits of being simple to both implement and analyze.  Note that the least-squares approximation $f_{\Upsilon,\Lambda}$ is given by
\be{
\label{ls_approximation}
f_{\Upsilon,\Lambda} = \argmin{p \in P_{\Lambda}} \frac{1}{M} \sum_{\bm{y} \in \Upsilon} \left | f(\bm{y}) - p(\bm{y}) \right |^2,
}
where $P_{\Lambda} = \spn \{ \phi_{\bm{n}} : \bm{n} \in \Lambda \}$ is the finite-dimensional approximation space.  Equivalently, the coefficients $\bm{c} = ( c_{\bm{n}} )_{\bm{n} \in \Lambda} $ of $f_{\Upsilon,\Lambda}$ are the solution of the algebraic least-squares problem 
\be{
\label{algLSintro}
\bm{c} = \argmin{\bm{x} \in \bbC^N} \nm{\bm{A} \bm{x} - \bm{b}}_{2},
}
where $\bm{A}= \left \{ \frac{1}{\sqrt{M}} \phi_{\bm{n}}(\bm{y}) \right \}_{\bm{y} \in \Upsilon,\bm{n} \in \Lambda} \in \bbC^{M \times N}$ and $\bm{b} =  \left \{ \frac{1}{\sqrt{M}} f(\bm{y}) \right \}_{\bm{y} \in \Upsilon} \in \bbC^{M}$.

\pbk
\textit{2.\ How to choose the orthonormal basis $\{ \psi_{\bm{n}} \}_{\bm{n}}$ and index set $\Lambda$.}  Smooth functions are typically well-approximated by polynomials, so we shall generally take $\{ \psi_{\bm{n}} \}_{\bm{n}}$ to be an orthonormal tensor-product polynomial basis.  Our main numerical examples consider tensor-product Legendre polynomials.  We also highlight the possibility of nonpolynomial approximations, for example using a cosine basis when $\Omega$ is compactly contained in $D = (-1,1)^d$.  Given the basis $\{ \psi_{\bm{n}} \}_{\bm{n}}$, we consider several standard choices for $\Lambda$, including total degree and hyperbolic cross index sets, or more generally, so-called \textit{lower} sets.  These sets have been used quite extensively for multivariate polynomial approximations in tensor-product domains (see \cite{AdcockCSFunInterp,ChkifaEtAl,ChkifaEtAlBreaking,ChkifaDownwardsCS,CohenDeVoreSchwabFoCM,MiglioratiJAT,MiglioratiEtAlFoCM} and references therein).

\pbk
\textit{3.\ How to choose the sample points $\Upsilon$.}  Our primary concern in this regard lies with the \textit{sampling efficiency} (or \textit{sample complexity}) of the approximation: namely, how large $M$ must be in relation to $N = |\Lambda|$ to ensure a good approximation.  The problem of designing optimal sampling points for high-dimensional polynomial approximation remains open even in tensor-product domains (although we note in passing some recent quasi-optimal constructions \cite{MiglioratiCohenOptimal}).  We shall therefore not attempt to solve it for irregular domains.  Instead, we consider straightforward random samplings.  Specifically, we draw $\bm{y}_1,\ldots,\bm{y}_M$ independently according to a suitable probability measure on $\Omega$ (for example, the uniform measure whenever $\Omega$ is compact).  {We throughout assume that it is computationally feasible to draw samples from this measure.}  Although simple, this approach permits concrete sample complexity estimates for a large class of domains $\Omega$ which are quadratic (up to a log factor) in $N = | \Lambda |$ for any dimension $d$.  Up to a domain-dependent constant which we determine, this log-quadratic sample complexity is the same as the corresponding result for compact tensor-product domains when the sample points are drawn from the uniform measure \cite{ChkifaEtAl}.\footnote{This scaling is essentially sharp.  As discussed in \cite{AdcockNecSamp} (based on a result of \cite{TrefPlatteIllCond}), in one dimension if the sample points are deterministic and exactly equispaced, then the least-squares approximation is ill-conditioned unless the number of sample points scales quadratically in the polynomial degree $N$.}

\subsection{Conditioning and stability}\label{ss:intro_stability}
The approach outlined above is certainly simple, and it is tempting to think that it can achieve high accuracy.  After all, the method computes a polynomial approximation in a domain, albeit an irregular one.  Unfortunately, there is an issue.  The matrix $\bm{A}$ of the system \R{algLSintro} is extremely ill-conditioned, even when $M \gg N$ (we estimate this ill-conditioning later in the paper for relevant examples).  This is due to the fact that the set $\{ \phi_{\bm{n}} \}_{\bm{n}}$ is not a basis for the space of square-integrable functions over $\Omega$, but rather a frame.  See \S \ref{ss:polyframe} for the definition of a frame.  Frames are typically redundant, meaning that any function $f$ has infinitely-many expansions of the form $f = \sum_{\bm{n}} c_{\bm{n}} \phi_{\bm{n}}$ with coefficients $\{ c_{\bm{n}} \}_{\bm{n}}$ in $\ell^2$.  When translated to the finite setting, this redundancy means that the truncated Gram matrix
\be{
\label{Grammatrix}
\bm{G}_{\Lambda} = \left \{ \ip{\phi_{\bm{m}}}{\phi_{\bm{n}}}_{L^2(\Omega,\mu)} \right \}_{\bm{m},\bm{n} \in \Lambda},
}
where $\mu$ is the measure defined in \R{mu_def}, is typically extremely poorly conditioned for large $N$ \cite{BADHframespart}.  Note that $\bbE(\bm{A}^* \bm{A} ) = \bm{G}_{\Lambda}$ if the sample points $\bm{y}_i$ are drawn independently according to $\mu$.  Hence the least-squares matrix $\bm{A}$ is expected to inherit similar ill-conditioning.

In the face of such ill-conditioning, one would usually expect it to be impossible to achieve high accuracy in floating point arithmetic.  Fortunately, this expectation turns out to be incorrect.  The frame property endows the problem with sufficient structure so that accurate, well-conditioned approximations can be computed via a simple regularization procedure.  In this paper we show that regularized least-squares solutions, computed via hard thresholding of the singular values of $\bm{A}$, yield well-conditioned approximations which converge rapidly down to the thresholding parameter $\epsilon$.  This parameter is typically set according to some desired target accuracy.

\rem{
We stress that the frame property is crucial in endowing the approximation with these properties, hence why we refer to this approach as \textit{polynomial frame approximations}.  Choosing $\{ \phi_{\bm{n}} \}_{\bm{n}}$ to be the monomial basis also leads to an exceedingly ill-conditioned problem, but one where high accuracy may not be possible.  The underlying reason for this is that the frame property guarantees existence of expansions $f = \sum_{\bm{n}} c_{\bm{n}} \phi_{\bm{n}}$ for which the coefficients $\{ c_{\bm{n}}  \}_{\bm{n}}$ decay (accuracy) and have bounded $\ell^2$-norm (stability).  See \cite{BADHFramesPart2,BADHframespart} for further discussion.
}

\subsection{Main results}
We now summarize our main results.  

\pbk
\textbf{Accuracy and conditioning.}\ Our first result concerns the accuracy and condition number of the regularized least-squares approximation.  As mentioned above, this approximation is constructed using a truncated SVD of the least-squares matrix $\bm{A}$ with a threshold parameter $\epsilon > 0$.  We write $f_{\Upsilon,\Lambda,\epsilon}$ for this approximation and $\bm{c}^{\epsilon}$ for its coefficients in the system $\{ \phi_{\bm{n}} \}_{\bm{n} \in \Lambda}$.

\thm{[Accuracy and conditioning]
\label{t:intro_acc_stab_cond}
There exists a constant $C_{\Upsilon,\Lambda,\epsilon} > 0$ such that
\be{
\label{coffeemug}
\| f - f_{\Upsilon,\Lambda,\epsilon} \|_{L^2(\Omega,\mu)} \leq \left ( 1 + C_{\Upsilon,\Lambda,\epsilon} \right ) E_{\Lambda,\epsilon}(f),
}
where
\bes{
E_{\Lambda,\epsilon}(f) = \inf \left \{ \nm{f - p }_{L^\infty(\Omega)} + \epsilon \nm{p}_{L^2(D,\nu)} : \ p  \in P_{\Lambda} \right \},
}
and $\mu$ is the measure given by \R{mu_def}.  Moreover, the coefficients $\bm{c}^{\epsilon}$ of $f_{\Upsilon,\Lambda,\epsilon}$ satisfy
\be{
\label{coeffboundintro}
\nm{\bm{c}^{\epsilon}}_2 = \| f_{\Upsilon,\Lambda,\epsilon} \|_{L^2(D,\nu)} \leq \frac{E_{\Lambda,\epsilon}(f)}{\epsilon},
}
and the absolute $(\ell^2,L^2)$-condition number of the reconstruction operator $\cL_{\Upsilon,\Lambda,\epsilon} : \bbC^M \rightarrow P_{\Lambda}, \bm{b} \mapsto f_{\Upsilon,\Lambda,\epsilon}$, where $\bm{b}$ is as in \R{algLSintro}, is at most $C_{\Upsilon,\Lambda,\epsilon}$.
}

See \S \ref{s:acc_stab}.  Several remarks are in order.  First, the bound \R{coffeemug} separates the accuracy of the regularized least-squares approximation into an approximation error term $E_{\Lambda,\epsilon}(f) $ depending only on $\epsilon$ and the space $P_{\Lambda}$ and independent of the samples $\Upsilon$, and a constant $C_{\Upsilon,\Lambda,\epsilon}$ depending on $\epsilon$, $\Upsilon$ and $P_{\Lambda}$.  In other words, $E_{\Lambda,\epsilon}(f) $ determines the rate of approximation, whereas $C_{\Upsilon,\Lambda,\epsilon}$ (more specifically, the requirement that $C_{\Upsilon,\Lambda,\epsilon} \lesssim 1$) determines the sample complexity.

Second, notice that $E_{\Lambda,\epsilon}(f) $ depends on how well $f$ can be approximated in $\Omega$ by polynomials $p \in P_{\Lambda}$ that do not grow too large on $D$.  The latter requirement -- which stems from the regularization carried out -- is an expression of stability, since a polynomial growing large on $D$ would necessarily have large coefficients.  Our main estimates for $E_{\Lambda,\epsilon}(f)$, given below, are derived by constructing polynomials which approximate $f$ at specified rates in $\Omega$ (depending on the smoothness of $f$), and which remain bounded on $D$.

Third, note that \R{coeffboundintro} ensures the stored values -- namely, the coefficients $\bm{c}^{\epsilon}$ -- cannot be too large in magnitude. This would otherwise result in ill-conditioning of the evaluation map $\bm{c}^{\epsilon} \mapsto f_{\Upsilon,\Lambda,\epsilon}(\bm{x})$.  While $\nm{\bm{c}^{\epsilon}}_{2}$ may be of magnitude roughly $1/\epsilon$ initially, once the approximation error $E_{\Lambda,\epsilon}(f)$ reaches close to the target accuracy $\epsilon$ we have $\nm{\bm{c}^{\epsilon}}_{2} \lesssim 1$.

\pbk
\textbf{Rate of approximation.}\  In \S \ref{s:approx_err} we analyze $E_{\Lambda,\epsilon}(f) $ for the main example considered in this paper, Legendre polynomials on $D = (-1,1)^d$.  We consider two standard choices of index sets $\Lambda$: the total degree index set $\Lambda = \Lambda^{\mathrm{TD}}_{n}$ defined in \R{LambdaTD} and the hyperbolic cross index set $\Lambda = \Lambda^{\mathrm{HC}}_{n}$ defined in \R{LambdaHC}.  The former is suitable for low-dimensional problems, but quickly becomes too large as $d$ increases.  The cardinality of the latter on the other hand scales much more mildly with $d$.  

Our main results are split into two cases:

\pbk
\textit{(i) $f$ smooth in $\Omega$ only.}\
In the first case, $f$ is smooth in $\Omega$ but may be nonsmooth, or even undefined in $D \backslash \Omega$.  If $\Omega$ is a Lipschitz domain and $f \in H^{m}(\Omega,\mu)$, where $H^{m}(\Omega,\mu)$ is the classical Sobolev space of order $m$ (see \R{SobSpace}), then we show that
\be{
\label{case1intro}
E_{\Lambda,\epsilon}(f) \leq \left \{ \begin{array}{ll}   c_{m,d,\Omega} \left ( n^{d-m} + \epsilon \right ) \| f \|_{H^{m}(\Omega,\mu)} & \Lambda = \Lambda^{\mathrm{TD}}_n \\   c_{m,d,\Omega}\left ( n^{\frac{d-m}{d}} + \epsilon \right ) \| f \|_{H^{m}(\Omega,\mu)} & \Lambda = \Lambda^{\mathrm{HC}}_n \end{array} \right . ,
}
where $c_{m,d,\Omega}>0$ is a constant depending on $m$, $d$ and $\Omega$ but independent of $f$.
See Theorem \ref{t:approxerr1} (we note in passing that the factor $d-m$ can be improved slightly to $\theta(d)-m$ where $\theta(d)$ is a particular constant satisfying $\theta(d) \leq d$).  This result asserts convergence at an algebraic rate depending on the smoothness of $f$ in $\Omega$ only.  However, it also exhibits the familiar curse of dimensionality.  In the case 
of the total degree index set $\Lambda^{\mathrm{TD}}$ the cardinality $N = |\Lambda^{\mathrm{TD}}| \asymp n^d$ as $n \rightarrow \infty$, and therefore
\bes{
n^{d-m} \asymp N^{\frac{d-m}{d}},\qquad n \rightarrow \infty,
}
whereas for the hyperbolic cross space (wherein $N = |\Lambda^{\mathrm{HC}}_n| \asymp n (\log (n) )^{d-1}$) one has
\bes{
n^{\frac{d-m}{d}} \asymp N^{\frac{d-m}{d}} ( \log(N))^{\frac{(m-d)(d-1)}{d}},\qquad n \rightarrow \infty.
}

\pbk
\textit{(ii) $f$ smooth in $D$.}\
In high-dimensional approximation a standard way to overcome the $d$-dependence in results such as \R{case1intro} is to assume certain anisotropic smoothness.  As we discuss in \S \ref{s:conclusion} it is currently unknown how to do this within the setting of case (i).  However, when $f$ has appropriate regularity over the whole of $D$ -- or equivalently, $f$ is the restriction to $\Omega$ of some appropriately regular function defined on $D$ --  then we have the following result.  If $f \in H^{m}_{\mix}(D,\nu) $, where $H^{m}_{\mix}(D,\nu)$ is the Sobolev space of dominating mixed smoothness on $D$ (see \R{MixSobSpace}), then
\be{
\label{case2intro}
E_{\Lambda,\epsilon}(f) \leq \left \{ \begin{array}{ll}  c_{m,d} \| f \|_{H^{m}_{\mix}(D,\nu)} n^{1-m} + \epsilon \| f \|_{L^2(D,\nu)}   & \Lambda = \Lambda^{\mathrm{TD}}_n \\ c_{m,d} \| f \|_{H^{m}_{\mix}(D,\nu)} n^{1-m} (\log (n))^{\frac{d-1}{2}}  + \epsilon \| f \|_{L^2(D,\nu)}   & \Lambda = \Lambda^{\mathrm{HC}}_n \end{array} \right . ,
}
where $c_{m,d} > 0$ is a constant depending on $m$ and $d$ but independent of $\Omega$ and $f$.
See Theorem \ref{t:approxerr2}.  Observe that
\bes{ \
n^{1-m} \asymp N^{\frac{1-m}{d}},\qquad N = | \Lambda^{\mathrm{TD}}_{n} |,
}
whereas
\bes{
n^{1-m} (\log (n))^{\frac{d-1}{2}} \asymp N^{1-m} (\log(N))^{(d-1)(m-1/2)},\qquad N = | \Lambda^{\mathrm{HC}}_{n} |.
}
Hence, up to the logarithmic factor, the hyperbolic cross index set $\Lambda^{\mathrm{HC}}_{n}$ achieves an algebraic rate of convergence that is independent of the dimension $d$, and therefore suitable for higher-dimensional problems.  Our numerical results in \S \ref{s:numerics} show computations using the hyperbolic cross index set for dimensions up to $d = 15$.

\pbk
Let us make several remarks.  First, we note that case (ii) requires absolutely no conditions on the domain $\Omega$, besides being measureable.  In particular, the domain can be extremely rough, as long as $f$ is smooth over the whole extended domain $D$.  In \S \ref{s:numerics} we show some numerical results of this type.  Second, \R{case2intro} behaves like $n^{1-m}$, not $n^{-m}$ as might be expected.  The additional power of $n$ stems from the presence of the $L^\infty(\Omega)$ norm in $E_{\Lambda,\epsilon}(f)$.  This factor can be improved whenever $\Omega$ is compactly contained in $D$, in which case one obtains a factor of the form $n^{1/2-m}$ (see Theorems \ref{t:approxerr1} and \ref{t:approxerr2}).  Third, when the sample points are drawn randomly and independently (as they are in this paper) it is possible to prove an estimate in expectation for the squared $L^2$-error of a slightly modified least-squares estimator (see \S \ref{s:L2bounds}) involving the $L^2$-norm approximation error 
\bes{
\tilde{E}_{\Lambda,\epsilon}(f) = \inf \left \{ \nm{f - p }_{L^2(\Omega,\mu)} + \epsilon \| p \|_{L^2(D,\nu)} : p \in P_{\Lambda} \right \}.
}
See Theorem \ref{t:truncated_est}.  Analogous to \R{case1intro} and \R{case2intro}, this quantity admits the following estimates.  First, if $\Omega$ is Lipschitz and $f \in H^{m}(\Omega,\mu)$ then
\be{
\label{case1introL2}
\tilde{E}_{\Lambda,\epsilon}(f) \leq \left \{ \begin{array}{ll} c_{m,d,\Omega} \left ( n^{-m} + \epsilon \right ) \| f \|_{H^{m}(\Omega,\mu)} & \Lambda = \Lambda^{\mathrm{TD}}_n \\ c_{m,d,\Omega} \left ( n^{-\frac{m}{d}} + \epsilon \right ) \| f \|_{H^{m}(\Omega,\mu)} & \Lambda = \Lambda^{\mathrm{HC}}_n \end{array} \right . .
}
Conversely, if $f \in H^{m}_{\mix}(D,\nu)$ then
\be{
\label{case2introL2}
\tilde{E}_{\Lambda,\epsilon}(f) \leq c_{m,d} \| f \|_{H^{m}_{\mix}(D,\nu)} n^{-m} + \epsilon \| f \|_{L^2(D,\nu)} ,\quad \mbox{$\Lambda = \Lambda^{\mathrm{TD}}_n$ or $\Lambda = \Lambda^{\mathrm{HC}}_n$}.
}
See Theorems \ref{t:approxerr1L2} and \ref{t:approxerr2L2}.  As with $E_{\Lambda,\epsilon}(f)$ above, this latter result for the hyperbolic cross index set shows how polynomial frame approximation can mitigate the curse of dimensionality.

\pbk
\textbf{Sample complexity.}\ 
Our final result concerns efficiency, i.e.\ sample complexity, of the approximation.  In view of Theorem \ref{t:intro_acc_stab_cond}, this corresponds to determining how large $M$ must be in order for the condition $C_{\Upsilon,\Lambda,\epsilon} \lesssim 1$ to hold.  Our main contribution is for the following class of domains $\Omega$:

\defn{
[$\lambda$-rectangle property]
A compact domain $\Omega$ has the $\lambda$-rectangle property for some $0<\lambda<1 $ if it can be written as a (possibly overlapping and uncountable) union
\bes{
\Omega = \bigcup_{R \in \cR} R,
}
of hyperrectangles $R$ satisfying
\bes{
 \inf_{R \in \cR} \mathrm{Vol}(R) = \lambda \mathrm{Vol}(\Omega).
}
}

Note that many domains of practical interest have this property.  However there are notable exceptions, including simplices and balls.  See \S \ref{ss:lambdarectangle} for further discussion.  

As we show in \S \ref{ss:Nikolskii}, when the samples $\bm{y}_m$ are chosen randomly and independently according to the uniform measure on $\Omega$ the sample complexity of the approximation {can in general be related to the constant of the $(L^2(\Omega,\mu),L^\infty(\Omega))$-Nikolskii inequality for the space $P_{\Lambda}$}.  We use the $\lambda$-rectangle property to get concrete estimates for this constant, culminating in the following result:

\thm{[Sample complexity]
Suppose that $\Omega \subseteq (-1,1)^d$ has the $\lambda$-rectangle property and let $P_{\Lambda}$ be constructed from the tensor Legendre polynomial basis on $(-1,1)^d$, where $\Lambda \subset \bbN^d_0$ is any lower set (see Definition \ref{d:lowerset}) of cardinality $| \Lambda | = N$.  Let $0 < \delta,\gamma < 1$ and $\bm{y}_1,\ldots,\bm{y}_M$ be independent and randomly drawn according to the uniform probability measure on $\Omega$.  Then
\bes{
C_{\Upsilon,\Lambda,\epsilon} \leq \frac{1}{\sqrt{1-\delta}},\qquad \forall \epsilon > 0,
}
with probability at least $1-\gamma$, provided
\bes{
M \geq N^2 \lambda^{-1} \left ( (1-\delta) \log(1-\delta)+\delta) \right )^{-1} \log(N/\gamma).
}
}

See Corollary \ref{c:sample_complexity}.  This result establishes log-quadratic scaling of the number of samples with the dimension of the polynomial space, extending a well-known result for tensor-product domains to a large class of irregular domains.  Note that this result holds for all lower sets, and in particular, the total degree and hyperbolic cross index sets discussed above.

\subsection{Related work}\label{ss:related}

The idea of approximating a function on an irregular domain by using an orthogonal basis on a bounding tensor-product domain is well established within the context of embedded or fictitious domain methods in numerical PDEs \cite{pasquettiFourEmbed} (see also \cite{boyd2005fourier}).  So-called \textit{Fourier extensions} or \textit{Fourier continuations} were studied in detail in \cite{BoydFourCont,brunoFEP}.  Applications to surface parametrization and numerical PDEs in complex geometries were considered in \cite{brunoFEP}  and \cite{albin2011,bruno2010high,lyon2010high} respectively.  Our work can be considered an extension of \cite{FEStability} from the univariate to the multivariate setting, although we use algebraic as opposed to trigonometric polynomials since these are more common in applications such as UQ.
Our work also extends recent research on computing polynomial approximations of functions defined on high-dimensional tensor-product domains.  This approach has received substantial interest recently, due to its applications in, notably, UQ.  See \cite{BASBCWMatheon,ChkifaEtAl,ChkifaEtAlBreaking,CohenDeVoreApproxPDEs,DavenportEtAlLeastSquares,DoostanOwhadiSparse,NarayanZhouCCP,KarniadakisUQCS} and references therein.  A consequence of this paper is that an irregular domain (either known or unknown) often presents no barrier to polynomial approximation of high-dimensional functions.  As noted, polynomial approximations are frequently used in practical UQ studies even when the domain is non-tensorial  (see \cite{SargsyanEtAlUQDisc,SargsyanEtAlDimUQ} and references therein).  Our work therefore provides a theoretical basis for these approaches.
Finally, we note that polynomial frame approximation is just once example of so-called \textit{numerical frame approximation}.  For a broader perspective on the uses of frames in numerical analysis and approximation, see \cite{BADHFramesPart2,BADHframespart}.

\section{Polynomial frame approximations}

\subsection{Notation}
We first require some further notation.  Throughout this paper $D \subseteq \bbR^d$ will be a domain with a probability measure $\nu$.  Typically, $D$ will be of tensor-product type, i.e.\
\be{
\label{Dtensor}
D = [a_1,b_1] \otimes \cdots \otimes [a_d,b_d] \subseteq \bbR^d,
}
where $-\infty \leq a_k < b_k \leq \infty$ and $\nu = \nu^{(1)} \otimes \cdots \otimes \nu^{(d)}$ will be a tensor-product of one-dimensional probability measures.  We write $L^2(D,\nu)$ for the space of square-integrable functions on $D$.

The $d$-dimensional variable is denoted by $\bm{y} = (y_1,\ldots,y_d) \in \bbR^d$.  Given $D$, we let $\Omega \subseteq D$ be a domain and define the probability measure $\mu$ by
\be{
\label{mu_def}
\D \mu(\bm{y}) = \frac{\bbI_{\Omega}(\bm{y})}{v_{\Omega}}   \D \nu(\bm{y}),\qquad v_{\Omega} = \int_{\Omega} \D \nu,
}
where $\bbI_{\Omega}$ is the indicator function of $\Omega$.  {We write $L^2(\Omega,\mu)$ for the space of square-integrable functions on $\Omega$ with respect to $\mu$.}

Throughout, $\bm{n} = (n_1,\ldots,n_d) \in \bbN^d_0$ denotes a multi-index.  Let $I \subseteq \bbN^d_0$ be a countable set of multi-indices and $\{ \psi_{\bm{n}} \}_{\bm{n} \in I}$ be an orthonormal basis of $L^2(D,\nu)$.  If $D$ is of the form \R{Dtensor}, then this basis will usually be of tensor product-type, i.e.
\bes{
\psi_{\bm{n}}(\bm{y}) = \prod^{d}_{k=1} \psi^{(k)}_{n_k}(y_k),
}
where $\{ \psi^{(k)}_{n_k} \}$ is an orthonormal basis of $L^2((a_k,b_k),\nu^{(k)})$.  Given $\{ \psi_{\bm{n}} \}_{\bm{n} \in I}$ we let 
\be{
\label{PolyFrame}
\phi_{\bm{n}} = \psi_{\bm{n}} \big |_{\Omega},\qquad \bm{n} \in I,
}
be the corresponding functions defined on $\Omega$.

\subsection{Multi-index sets}
Our interest lies in computing finite approximations in the system \R{PolyFrame}.  To this end, let $\Lambda \subset I$
be a finite multi-index set and define
\bes{
P_{\Lambda} = \spn \left \{ \phi_{\bm{n}} : \bm{n} \in \Lambda \right \} \subset L^2(\Omega,\mu),
}
as the finite-dimensional space within which we seek an approximation to $f$.  We consider the following three standard choices of multi-index sets.  The \textit{tensor product} set
\be{
\label{LambdaTP}
\Lambda = \Lambda^{\mathrm{TP}}_n = \left \{ \bm{n}\in \bbN^d_0 : | \bm{n} |_{\infty} \leq n \right \},
}
where $n \in \bbN_0$ and $| \bm{n} |_{\infty} = \max_{k=1,\ldots,d} | n_k |$, the \textit{total degree} set
\be{
\label{LambdaTD}
\Lambda = \Lambda^{\mathrm{TD}}_n = \left \{ \bm{n}\in \bbN^d_0 : | \bm{n} |_1 \leq n \right \},
}
where $| \bm{n} |_{1} = |n_1|+\ldots+|n_d|$, and the (isotropic) \textit{hyperbolic cross} set
\be{
\label{LambdaHC}
\Lambda = \Lambda^{\mathrm{HC}}_n = \left \{ \bm{n} \in \bbN^d_0 : |\bm{n}|_{\mathrm{hc}} \leq n+1 \right \},\qquad |\bm{n}|_{\mathrm{hc}} = \prod^{d}_{k=1} (n_k+1).
}
Note that the cardinality $N = |\Lambda^{\mathrm{TP}}_n | = (n+1)^d$ usually grows too quickly with $n$ in high dimensions to be practical.  The total degree set, with cardinality
\bes{
N = |\Lambda^{\mathrm{TD}}_n | = \left ( \begin{array}{c} n+d \\ d \end{array} \right ),
}
mitigates this issue to some extent, but still typically grows too rapidly for moderate to high-dimensional problems.  Hyperbolic cross index sets are a practical alternative in this case.  An exact formula for the cardinality of the hyperbolic cross $\Lambda^{\mathrm{HC}}_n$ in terms of $n$ and $d$ is not known, but there are a variety of upper bounds, including:
\bes{
\left | \Lambda^{\mathrm{HC}}_n  \right | \leq \left \lfloor (n+1) (1 + \log(n+1) )^{d-1} \right \rfloor.
}
See, for example, \cite[Prop.\ A.1]{MiglioratiThesis}.  

The above three multi-index sets are all examples of so-called \textit{lower} sets (also known as \textit{downward closed} or \textit{monotone} sets -- see, for example, \cite{CohenDeVoreApproxPDEs,DynFloaterDownward}):

\defn{
\label{d:lowerset}
A multi-index set $\Lambda \subseteq \bbN^d_0$ is lower if whenever $\bm{n} = (n_1,\ldots,n_d) \in \Lambda$ and $\bm{n}' = (n'_1,\ldots,n'_d)$ satisfies $n'_k \leq n_k $ for all $k$ then $\bm{n}' \in \Lambda$.
}

In our main results regarding efficiency, we establish sample complexity estimates which are valid for arbitrary lower sets.  While we shall not do it in this paper, such generality allows for the possibility of considering other multi-index sets, e.g.\ anisotropic hyperbolic cross index sets, which may be defined by \textit{a priori} or \textit{a posteriori} estimates, or computed adaptively.

\subsection{Polynomial frames}\label{ss:polyframe}
We first recall the definition of a frame (see, for example, \cite{christensen2003introduction}):

\defn{
A countable set $\{ \phi_{\bm{n}}\}_{\bm{n} \in I}$ of a Hilbert space $H$ is a frame if there exist constants $0 < A \leq B < \infty$ such that
\be{
\label{FrameBounds}
A \nm{f}^2 \leq \sum_{\bm{n} \in I} | \ip{f}{\phi_{\bm{n}}} |^2 \leq B \nm{f}^2,\quad \forall f \in H,
}
where $\ip{\cdot}{\cdot}$ and $\nm{\cdot}$ are the inner product and norm respectively on $H$.
}
Let $\{ \phi_{\bm{n}} \}_{\bm{n} \in I}$ be the system constructed in \R{PolyFrame}.  It is straightforward to see that this is a frame for $H = L^2(\Omega,\mu)$.  Indeed, let $f \in L^2(\Omega,\mu)$ and $\tilde{f}$ be its extension by zero to $D$.  Then by Parseval's relation for the orthonormal basis $\{ \psi_{\bm{n}} \}_{\bm{n} \in I}$,
\bes{
\sum_{\bm{n} \in I} | \ip{f}{\phi_{\bm{n}}}_{L^2(\Omega,\mu)} |^2 =  \sum_{\bm{n} \in I} \frac{| \ip{\tilde{f}}{\psi_{\bm{n}}}_{L^2(D,\nu)} |^2}{(v_{\Omega})^2} = \frac{\| \tilde{f} \|^2_{L^2(D,\nu)}}{(v_{\Omega})^2} = \frac{\| f \|^2_{L^2(\Omega,\mu)}}{v_{\Omega}},
}
where $v_{\Omega}$ is given by \R{mu_def}. 
Hence \R{FrameBounds} holds with $A = B = 1/v_{\Omega}$, making this system a frame.  Frames such as this for which $A = B$ are known as \textit{tight} frames.

A general property of frames is their redundancy: any $f \in H$ can have infinitely-many expansions $f = \sum_{\bm{n} \in I} c_{\bm{n}} \phi_{\bm{n}}$ with coefficients $\{ c_{\bm{n}} \}_{\bm{n} \in I} \in \ell^2(I)$.  It is easy to see how redundancy occurs in the polynomial frame.  Indeed, let $\tilde{f}$ be any extension of $f$ to $L^2(D,\nu)$ and define
\bes{
c_{\bm{n}} = \ip{\tilde{f}}{\psi_{\bm{n}}}_{L^2(D,\nu)},
}
as the coefficients of $\tilde{f}$ in the orthonormal basis $\{ \psi_{\bm{n}} \}_{\bm{n} \in I }$.  Then $\{ c_{\bm{n}} \}_{\bm{n} \in I} \in \ell^2(I)$ and
\bes{
\sum_{\bm{n} \in I} c_{\bm{n}} \phi_{\bm{n}} = \sum_{\bm{n} \in I} c_{\bm{n}} \psi_{\bm{n}}  \Big |_{\Omega} = \tilde{f} \big |_{\Omega} = f.
}
Since there are infinitely many extensions of $f$ to $L^2(D,\nu)$, each with distinct coefficients $\{ c_{\bm{n}} \}_{\bm{n} \in I}$, it follows that there are infinitely many representations of $f$ in the frame $\{ \phi_{\bm{n}} \}_{\bm{n} \in I}$.

\subsection{Least-squares polynomial frame approximations}
Let $\Upsilon = \{ \bm{y}_1,\ldots,\bm{y}_M \} \subset \Omega$ be a set of $M$ distinct points (for the moment we choose not to specify their distribution) and $\Lambda$ be a finite set of multi-indices of size $|\Lambda|  = N$, where $N \leq M$.  {Consider the} approximation to $f$ in the space $P_{\Lambda}$ by discrete least-squares fitting:
\bes{
f_{\Upsilon,\Lambda} = \argmin{p \in P_{\Lambda}} \frac{1}{M} \sum_{\bm{y} \in \Upsilon} \left | f(\bm{y}) - p(\bm{y}) \right |^2.
}
If $f_{\Upsilon,\Lambda} $ is expressed as
\bes{
f_{\Upsilon,\Lambda} = \sum_{\bm{n} \in \Lambda} c_{\bm{n}} \phi_{\bm{n}},
}
then this is equivalent to the algebraic least-squares problem
\be{
\label{algLS}
\bm{c} = ( c_{\bm{n}} )_{\bm{n} \in \Lambda} = \argmin{\bm{x} \in \bbC^N} \nm{\bm{A} \bm{x} - \bm{b}}_{2},
}
where
\bes{
\bm{A} = \bm{A}_{\Upsilon,\Lambda} = \left ( \frac{1}{\sqrt{M}} \phi_{\bm{n}}(\bm{y}) \right )_{\bm{y} \in \Upsilon,\bm{n} \in \Lambda} \in \bbC^{M \times N},\qquad \bm{b} = \bm{b}_{\Upsilon} =  \left ( \frac{1}{\sqrt{M}} f(\bm{y}) \right )_{\bm{y} \in \Upsilon} \in \bbC^{M}.
}
{Note that $\bm{A}$ may fail to be full rank -- e.g.\ if the points $\Upsilon$ are chosen poorly or the functions $\phi_{\bm{n}}$, $\bm{n} \in \Lambda$, are linearly dependent -- in which case \R{algLS} does not have a unique solution.  However, even if it is full rank, as} mentioned in \S \ref{ss:intro_stability} and shown explicitly in \S \ref{ss:ill_conditioning} below, $\bm{A}$ is typically severely ill-conditioned for large $N$.  Hence it is necessary to regularize \R{algLS}.  We shall do this via truncated singular value decompositions (i.e.\ spectral filtering).\footnote{Related strategies such as Tikhonov regularizatioxn could be used instead, with some changes to the ensuing presentation.}

To this end, suppose that $\bm{A}$ has singular values $\{ \sigma_{\bm{n}} \}_{\bm{n} \in \Lambda}$ and singular value decomposition $\bm{A} = \bm{U} \bm{\Sigma} \bm{V}^*$,
where $\bm{U} \in \bbC^{M \times M}$, $\bm{\Sigma} \in \bbR^{M \times N}$ and $\bm{V} \in \bbC^{N \times N}$.  Define
\bes{
\bm{A}_{\epsilon} = \bm{A}_{\Upsilon,\Lambda,\epsilon} = \bm{U} \bm{\Sigma}_{\epsilon} \bm{V}^*,
}
where the diagonal matrix $\bm{\Sigma}_{\epsilon} $ has $\bm{n}^{\rth}$ entry $\sigma_{\bm{n}}$ if $\sigma_{\bm{n}} > \epsilon$ and zero otherwise.  Then the truncated SVD least-squares approximation is defined as
\be{
\label{truncLSapprox}
f_{\Upsilon,\Lambda,\epsilon} = \sum_{\bm{n} \in \Lambda} (\bm{c}_{\epsilon})_{\bm{n}} \phi_{\bm{n}},
}
where its coefficients $\bm{c}_{\epsilon}$ are given by
\bes{
 \bm{c}_{\epsilon} = \left ( \bm{A}_{\Upsilon,\Lambda,\epsilon} \right )^{\dag} \bm{b}_{\Upsilon} = \bm{V} (\bm{\Sigma}_{\epsilon} )^{\dag} \bm{U}^* \bm{b}_{\Upsilon}.
}
Here $\dag$ denotes the pseudoinverse.  We consider this approximation from now on.  Note that the regularization parameter $\epsilon$ is usually set in relation to some desired target accuracy (see \S \ref{s:numerics}).

\subsection{Main example}
We end this section by introducing our main example.  This is the case where $\Omega$ is bounded and, without loss of generality, contained in $D = (-1,1)^d$, and where $\{ \psi_{\bm{n}} \}_{\bm{n} \in \bbN^d_0}$ is the tensor Legendre polynomial basis on $D$ corresponding to the uniform probability measure $\D \nu(\bm{y}) = 2^{-d} \D \bm{y}$.  When normalized with respect to the uniform probability measure on $D$, this basis is defined by
\bes{
\psi_{\bm{n}}(\bm{y}) = \prod^{d}_{k=1} \sqrt{2 n_k+1} P_{n_k}(y_k),
}
where $P_{n}$ is the $n^{\rth}$ classical Legendre polynomial (see Appendix \ref{s:Legbackground}).  For the truncated index set, we let $\Lambda = \Lambda^{\mathrm{TD}}_n$ or $\Lambda = \Lambda^{\mathrm{HC}}_n$ be either the total degree \R{LambdaTD} or hyperbolic cross \R{LambdaHC} index set with index $n$.  We also assume that the sampling points $\bm{y}_{1},\ldots,\bm{y}_{M}$ are drawn independently according to the measure $\mu$, which in this case is the uniform probability measure on $\Omega$:
\be{
\label{UnifProbMeas}
\D \mu(\bm{y}) = \frac{1}{\mathrm{Vol}(\Omega)} \D \bm{y}.
}
While this approach leads to concrete, $d$-independent sample complexity estimates for many domains, we do not claim that it is an optimal sampling procedure.  See \S \ref{s:numerics}--\ref{s:conclusion} for further discussion.

{
\rem{
As mentioned, we assume that it is computationally feasible to draw samples from $\mu$.  For the numerical examples shown later, this is achieved by rejection sampling.  Depending on the domain, however, and especially in high dimensions, this may be a substantial challenge.  
}
}

\section{Accuracy and conditioning}\label{s:acc_stab}
We now investigate the accuracy and conditioning of the approximation \R{truncLSapprox}.  In \S \ref{ss:ill_conditioning} we show that least-squares matrix $\bm{A}$ is ill-conditioned for large $N$, thus explaining why regularization is needed.  Next, in \S \ref{ss:accstabconsts} we introduce the key constant $C_{\Upsilon,\Lambda,\epsilon}$, and in \S \ref{ss:mainresaccstab} we give the main result of this section.  Note that the approach in \S \ref{ss:accstabconsts}--\ref{ss:mainresaccstab} follows that of \cite{BADHFramesPart2} (which applies to general frames).

\subsection{Ill-conditioning of the matrix $A$}\label{ss:ill_conditioning}

{Unless the frame happens to be a Riesz basis (which is not the case in our setting) frame approximations} always lead to ill-conditioned least-squares matrices {for sufficiently large truncation space $\Lambda$} \cite[Lem.\ 5]{BADHframespart}.  In the case of the polynomial frame, this is related to the Remez inequality for the polynomial space $P_{\Lambda}$ over $\Omega$ and $D$.  To see this, observe that the minimal and maximal singular values of $\bm{A}$ are
\eas{
\sigma_{\min}(\bm{A}) = \inf_{\substack{p \in P_{\Lambda} \\ p \neq 0}} \left \{ \frac{\sqrt{\frac1M \sum_{\bm{y} \in \Upsilon} | p(\bm{y}) |^2 }}{\| p \|_{L^2(D,\nu)}} \right \},
\quad
\sigma_{\max}(\bm{A})= \sup_{\substack{p \in P_{\Lambda} \\ p \neq 0}} \left \{ \frac{\sqrt{\frac1M \sum_{\bm{y} \in \Upsilon} | p(\bm{y}) |^2 }}{\| p \|_{L^2(D,\nu)}} \right \}.
}
For simplicity, assume that the constant function is contained in $P_{\Lambda}$.  This will hold in all examples considered later.  Letting $p(\bm{y}) = 1$ we get $\sigma_{\max}(\bm{A}) \geq 1$.  Conversely, note that $\frac1M \sum_{\bm{y} \in \Upsilon} | p(\bm{y}) |^2 \leq \| p \|^2_{L^\infty(\Omega)}$ and let $\rN(P_{\Lambda},D,\nu) > 0$ be the optimal constant such that
\bes{
\| p \|_{L^\infty(D)} \leq \rN(P_{\Lambda},D,\nu) \| p \|_{L^2(D,\nu)},\quad \forall p \in P_{\Lambda}.
}
We refer to this as an $(L^2(D,\nu),L^{\infty}(D))$-\textit{Nikolskii inequality} for the space $P_{\Lambda}$.  Inequalities such as these will be discussed further in \S \ref{s:sampcomp}, since they are pivotal in estimating the sample complexity of the approximation.  This gives
\bes{
\frac{1}{\sigma_{\min}(\bm{A})} \geq \left ( \rN(P_{\Lambda},D,\nu) \right )^{-1}  \sup \left \{ \frac{\| p \|_{L^\infty(D)}}{\| p \|_{L^{\infty}(\Omega)}} : p \in P_{\Lambda},\ p \neq 0 \right \},
}
and therefore
\be{
\label{cond_lower}
\mathrm{cond}(\bm{A}) \geq \frac{\rR(P_{\Lambda},\Omega,D)}{\rN(P_{\Lambda},D,\nu)},
}
where $\rR(P_{\Lambda},\Omega,D)$ is the constant in Remez's inequality for the domains $\Omega$ and $D$:
\bes{
\| p \|_{L^\infty(D)} \leq \rR(P_{\Lambda},\Omega,D) \| p \|_{L^\infty(\Omega)},\quad p \in P_{\Lambda}.
}
Note that the bound \R{cond_lower} is {completely deterministic, and independent of the samples $\Upsilon$.}

Typically, the right-hand side of \R{cond_lower} will grow rapidly with $N$.  To see why, note first that the Nikolskii constant is usually at most algebraic in $N = | \Lambda |$.  In particular, if $\nu$ is the uniform measure on $D$ and $\Lambda$ is a lower set, then $\rN(P_{\Lambda},D,\nu) \leq N^2$ \cite[Thm.\ 6]{MiglioratiJAT} (see also the proof of Theorem \ref{t:Nikolskii_lambda_rect}).  Similar bounds are found in \cite{MiglioratiJAT} for other ultraspherical and Jacobi measures.  Conversely, the constant $\rR(P_{\Lambda},\Omega,D)$ is typically exponentially large in $N$.  Its precise behaviour depends on the domain $\Omega$ and the index set $\Lambda$, and for the sake of brevity, we will not consider this issue in depth.  However, we note in passing that in the one-dimensional case for example, if $\Lambda = \{0,\ldots,N-1\}$ and $D = (-1,1)$ then
\bes{
R(P_{\Lambda},\Omega,(-1,1)) \leq T_{N-1}(4/|\Omega|-1),
}
where $T_{N-1}$ is the $(N-1)^{\rth}$ Chebyshev polynomial and $| \Omega |$ denotes the Lebesgue measure of $\Omega$.  Moreover, equality holds if $\Omega = [-1,-1+|\Omega|]$ in which case one has the exponential growth
\bes{
R(P_{\Lambda},\Omega,(-1,1)) \geq \frac12 \left ( \frac{4}{|\Omega|}-1 \right )^{N-1}.
}
We refer to \cite{GanzburgRemez} for further information, including results in higher dimensions, as well as to \cite{TemlyakovRemezHC} for results on multivariate Remez inequalities for hyperbolic cross index sets.

\subsection{Key constants}\label{ss:accstabconsts}
For convenience we now define the following operator
\bes{
\cT_{\Lambda} : \bbC^{N} \rightarrow P_{\Lambda},\ \bm{c} = \{ c_{\bm{n}} \}_{\bm{n} \in \Lambda} \mapsto \sum_{\bm{n} \in \Lambda} c_{\bm{n}} \phi_{\bm{n}}.
}
This is commonly referred to as the \textit{synthesis} operator in frame theory.  We now let
\be{
\label{Cconst}
C_{\Upsilon,\Lambda,\epsilon} = \max \left \{ C'_{\Upsilon,\Lambda,\epsilon} , C''_{\Upsilon,\Lambda,\epsilon} \right \},
}
where
\be{
\begin{split}
C'_{\Upsilon,\Lambda,\epsilon} &= \max_{\substack{\bm{b} \in \bbC^M \\ \nm{\bm{b}}_2 = 1 }} \nm{\cT_{\Lambda} (\bm{A}_{\Upsilon,\Lambda,\epsilon})^{\dag} \bm{b}}_{L^2(\Omega,\mu)},
\\
C''_{\Upsilon,\Lambda,\epsilon} &= \epsilon^{-1} \max_{\substack{\bm{d} \in \bbC^N \\ \nm{\bm{d}}_2 = 1 }}  \nm{\cT_{\Lambda} \bm{d} - \cT_{\Lambda} (\bm{A}_{\Upsilon,\Lambda,\epsilon})^{\dag} \bm{A}_{\Upsilon,\Lambda} \bm{d} }_{L^2(\Omega,\mu)}. 
\end{split}
\label{C12const}
}
It is useful to interpret these constants.  First, define the \textit{reconstruction} operator
\be{
\label{Loperator}
\cL_{\Upsilon,\Lambda,\epsilon}: \bbC^{M} \rightarrow P_{\Lambda} ; \bm{b} \mapsto \cT_{\Lambda} (\bm{A}_{\Upsilon,\Lambda,\epsilon})^{\dag} \bm{b}.
}
This operator takes a vector of samples $\bm{b} \in \bbC^{M}$ to its truncated SVD approximation in $P_{\Lambda}$.  In particular, if
\be{
\label{SUpsilon}
\cS_{\Upsilon} : L^\infty(\Omega) \rightarrow \bbC^{M}; f \mapsto \left \{ \frac{1}{\sqrt{M}} f(\bm{y}) \right \}_{\bm{y} \in \Upsilon},
}
is the operator taking a function $f$ to its samples then
\be{
\label{fapproxoperators}
f_{\Upsilon,\Lambda,\epsilon} = \cL_{\Upsilon,\Lambda,\epsilon} \cS_{\Upsilon} f.
}
The constant $C'_{\Upsilon,\Lambda,\epsilon}$ is precisely the operator norm -- or equivalently, since it is a linear operator, the absolute condition number -- of $\cL_{\Upsilon,\Lambda,\epsilon}$ with resect to the $\ell^2$- and $L^2(\Omega,\mu)$-norms:
\bes{
C'_{\Upsilon,\Lambda,\epsilon} =  \max_{\substack{\bm{b} \in \bbC^M \\ \nm{\bm{b}}_2 = 1 }} \nm{\cL_{\Upsilon,\Lambda,\epsilon}\bm{b}}_{L^2(\Omega,\mu)}.
}
In other words, boundedness of $C_{\Upsilon,\Lambda,\epsilon}$ implies robustness of the approximation to perturbations in the data (e.g.\ noise).  On the other hand, $C''_{\Upsilon,\Lambda,\epsilon} $ also has the equivalent definition
\bes{
C''_{\Upsilon,\Lambda,\epsilon} = \epsilon^{-1} \sup \left \{ \nm{p - p_{\Upsilon,\Lambda,\epsilon} }_{L^2(\Omega,\mu)} : p \in P_{\Lambda},\ \| p \|_{L^2(D,\nu)} = 1 \right \}.
}
In particular, $C''_{\Upsilon,\Lambda,0} = 0$ since the unregularized mapping $f \mapsto f_{\Upsilon,\Lambda,0}$ is a projection onto $P_{\Lambda}$.  When $\epsilon > 0$ this constant measures how close the map $f \mapsto f_{\Upsilon,\Lambda,\epsilon}$ is to being a projection onto $P_{\Lambda}$.

\subsection{Main result on accuracy and conditioning}\label{ss:mainresaccstab}

\thm{
\label{t:acc_stab_cond}
Let $f \in L^{\infty}(\Omega,\mu)$ and suppose that $f_{\Upsilon,\Lambda,\epsilon}$ is the truncated SVD least-squares approximation.  Then
\eas{
\nm{f - f_{\Upsilon,\Lambda,\epsilon}}_{L^2(\Omega,\mu)} \leq  \left ( 1 + C'_{\Upsilon,\Lambda,\epsilon} \right ) \| f - p \|_{L^\infty(\Omega)} + \epsilon C''_{\Upsilon,\Lambda,\epsilon} \| p \|_{L^2(D,\nu)} 
 \leq  \left ( 1 + C_{\Upsilon,\Lambda,\epsilon} \right ) E_{\Lambda,\epsilon}(f),
}
where $C'_{\Upsilon,\Lambda,\epsilon}$, $C''_{\Upsilon,\Lambda,\epsilon}$ and $C_{\Upsilon,\Lambda,\epsilon}$ are as in \R{C12const} and \R{Cconst} respectively, and 
\be{
\label{ELambdaDef}
E_{\Lambda,\epsilon}(f) = \inf \left \{ \nm{f - p }_{L^\infty(\Omega)} + \epsilon \| p \|_{L^2(D,\nu)} : \ p \in P_{\Lambda} \right \}.
}
Moreover, the coefficients $\bm{c}^{\epsilon}$ of $f_{\Upsilon,\Lambda,\epsilon}$ satisfy
\bes{
\nm{\bm{c}^{\epsilon}}_2 = \| f_{\Upsilon,\Lambda,\epsilon} \|_{L^2(D,\nu)} \leq \frac{E_{\Lambda,\epsilon}(f)}{\epsilon}.
}
}

\prf{
Let $p = \cT_{\Lambda} \bm{c}$ for some $\bm{c} \in \bbC^N$.  Then, recalling the definitions of the constants $C'_{\Upsilon,\Lambda,\epsilon}$ and $C''_{\Upsilon,\Lambda,\epsilon}$, we have
\eas{
\nm{f - f_{\Upsilon,\Lambda,\epsilon}}_{L^2(\Omega,\mu)} & \leq \nm{f - p }_{L^2(\Omega,\mu)} + \nm{p_{\Upsilon,\Lambda,\epsilon}-f_{\Upsilon,\Lambda,\epsilon} }_{L^2(\Omega,\mu)} + \nm{ p - p_{\Upsilon,\Lambda,\epsilon}}_{L^2(\Omega,\mu)}
\\
& \leq \| f - p \|_{L^2(\Omega,\mu)} + C'_{\Upsilon,\Lambda,\epsilon} \| \cS_{\Upsilon}(f-p) \|_{2} + \epsilon C''_{\Upsilon,\Lambda,\epsilon} \| p \|_{L^2(D,\nu)}
\\
& \leq \left ( 1 + C'_{\Upsilon,\Lambda,\epsilon} \right ) \| f - p \|_{L^\infty(\Omega)} + \epsilon C''_{\Upsilon,\Lambda,\epsilon} \| p \|_{L^2(D,\nu)},
}
which gives the first result.  Note that in the third step we use \R{SUpsilon} to deduce that $\| \cS_{\Upsilon}(f-p) \|_{2} \leq \| f - p \|_{L^\infty(\Omega)}$ and the fact that $\mu$ is a probability measure, which implies that $\| f - p \|_{L^2(\Omega,\mu)} \leq \| f - p \|_{L^\infty(\Omega)}$.  For the second result, we first use Parseval's identity to give $\nm{\bm{c}^{\epsilon}}_2 = \| f_{\Upsilon,\Lambda,\epsilon} \|_{L^2(D,\nu)}$ and then write
\be{
\label{coeffsplit}
\| f_{\Upsilon,\Lambda,\epsilon} \|_{L^2(D,\nu)} \leq \| f_{\Upsilon,\Lambda,\epsilon}-p_{\Upsilon,\Lambda,\epsilon} \|_{L^2(D,\nu)} + \| p_{\Upsilon,\Lambda,\epsilon} \|_{L^2(D,\nu)}.
}
Consider the first term.  By \R{fapproxoperators} we have
\be{
\label{orange}
\begin{split}
\| f_{\Upsilon,\Lambda,\epsilon}-p_{\Upsilon,\Lambda,\epsilon} \|_{L^2(D,\nu)}  &= \nm{ \cT_{\Lambda} (\bm{A}_{\Upsilon,\Lambda,\epsilon})^{\dag} \cS_{\Upsilon}(f-p)  }_{L^2(D,\nu)} 
\\
&= \nm{ (\bm{A}_{\Upsilon,\Lambda,\epsilon})^{\dag} \cS_{\Upsilon}(f-p)  }_2 \leq \frac{1}{\epsilon} \| \cS_{\Upsilon}(f-p) \|_{2}
 \leq \frac{1}{\epsilon} \| f - p \|_{L^\infty(\Omega)}.
 \end{split}
}
Here in the second step we use Parseval's identity, in the third step we use standard properties of the SVD and in the fourth step we use \R{SUpsilon}.  Now consider the second term of \R{coeffsplit}.  Observe that $\bm{A}_{\Upsilon,\Lambda} = \cS_{\Upsilon} \cT_{\Lambda}$.
Hence, if $p = \cT_{\Lambda} \bm{c}$ then, using standard properties of the SVD once more, we get
\eas{
\| p_{\Upsilon,\Lambda,\epsilon} \|_{L^2(D,\nu)} &= \nm{\cT_{\Lambda} (\bm{A}_{\Upsilon,\Lambda,\epsilon})^{\dag} \cS_{\Upsilon} \cT_{\Lambda} \bm{c} }_{L^2(D,\nu)}
 = \nm{(\bm{A}_{\Upsilon,\Lambda,\epsilon})^{\dag} \bm{A}_{\Upsilon,\Lambda} \bm{c} }_2
 \leq \nm{\bm{c}}_{2}
 = \| p \|_{L^2(D,\nu)}.
}
Combining this with \R{orange} and substituting both into \R{coeffsplit} now gives the second result.
}

A few remarks are in order.  First, to guarantee accuracy and good (absolute) conditioning of the approximation we need to ensure that $C_{\Upsilon,\Lambda,\epsilon} \lesssim 1$.  This constant depends on the polynomial space, the data and the threshold $\epsilon$, but is independent of the function $f$.  In \S \ref{s:sampcomp} we derive bounds for this constant.  Second, once $C_{\Upsilon,\Lambda,\epsilon}$ is bounded, the approximation error is determined via the term $E_{\Lambda,\epsilon}(f) $, which depends on $f$ and the polynomial space but is independent of the data.  We estimate this term for functions in certain Sobolev spaces in \S \ref{s:approx_err}.  

Third, we notice that the coefficients of the ensuing approximation are bounded by the the approximation error divided by $\epsilon$.  Thus, although the coefficients may initially be $\ord{1/\epsilon}$, they are $\ord{1}$ in the limit as the dimension $N$  of the approximation space $P_{\Lambda}$ tends to infinity.  Note that bounded coefficients are particularly important for practical computations, since these are the values that will be stored.  Indeed, if the coefficients could grow arbitrarily large in relation to the function $f$ then the pointwise evaluation operator $\bm{c}^{\epsilon} \mapsto f_{\Upsilon,\Lambda,\epsilon}(\bm{x})$ would be ill-conditioned.

Fourth and finally, we note that $C_{\Upsilon,\Lambda,\epsilon} \leq 1/(\sqrt{v_{\Omega}} \epsilon)$ for any $\Upsilon$, $\Lambda$ and $\epsilon > 0$ \cite[Prop.\ 4.6]{BADHFramesPart2}.  In other words, the ill-conditioning of the reconstruction operator scales at worst like $1/\epsilon$.

\section{Approximation error for Legendre polynomial frames}\label{s:approx_err}

We now consider the approximation error $E_{\Lambda,\epsilon}(f)$, defined by \eqref{ELambdaDef}.  In doing so, we treat the following two scenarios separately:
\bull{
\item[(i)] $f$ defined and smooth over $D$,
\item[(ii)] $f$ undefined or nonsmooth over $D$.
}
We first require several notions of smoothness.  Let
\be{
\label{SobSpace}
H^{m}(\Omega,\mu) = \left \{ f \in L^2(\Omega,\mu) : \cD^{\bm{j}} f \in L^2(\Omega,\mu) : | \bm{j} |_1 \leq m \right \},
}
be the classical Sobolev spaces of index $m \geq 0$ on $\Omega$, with norm
\bes{
\nm{f}_{H^m(\Omega,\mu)} = \sqrt{\sum_{|\bm{j}|_1 \leq m} \nm{\cD^{\bm{j}} f}^2_{L^2(\Omega,\mu)} }.
}
Here $\cD^{\bm{j}} = \frac{\partial^{| \bm{j} |_1}}{\partial^{j_1}_{y_1} \cdots \partial^{j_d}_{y_d}}$ is the partial derivative operator of order $\bm{j}$.  These spaces are suitable for approximations using the tensor product or total degree spaces in low dimensions.  For moderate to high dimensions, we instead consider Sobolev spaces of dominating mixed smoothness:
\be{
\label{MixSobSpace}
H^{m}_{\mix}(D,\nu) = \left \{ f \in L^2(D,\nu) : \cD^{\bm{j}} f \in L^2(D,\nu) : | \bm{j} |_{\infty} \leq m \right \},
}
with norm
\bes{
\nm{f}_{H^m_{\mix}(D,\nu)} = \sqrt{\sum_{|\bm{j}|_{\infty} \leq m} \nm{\cD^{\bm{j}} f}^2_{L^2(D,\nu)} }.
}

\subsection{Results for the classical Sobolev spaces $H^m$}

We first consider the tensor product and total degree index sets:

\thm{
\label{t:approxerr1}
Let $P_{\Lambda}$ be constructed from the tensor Legendre polynomial basis on $L^2(D,\nu)$, where $D = (-1,1)^d$, $\nu$ is the uniform measure on $D$, and $\Lambda = \Lambda_{n}$ is either the tensor product \R{LambdaTP} or total degree \R{LambdaTD} index set of degree $n$.  If $\Omega \subseteq D$ and $f \in H^{m}(D,\nu)$ for some $m > d/2$, then
\bes{
E_{\Lambda,\epsilon}(f) \leq c_{m,d} \| f \|_{H^{m}(D,\nu)} n^{\theta(d)-m}+ \epsilon \| f \|_{L^2(D,\nu)},
}
where 
\be{
\label{thetadef}
\theta(d) = \left \{ \begin{array}{ll} \frac{d(2d+1)}{2d+2} & \mbox{odd $d$} \\ \frac{d(2d+3)}{2d+4} & \mbox{even $d$} \end{array} \right . .
}
Conversely, if $\Omega \subseteq D$ is Lipschitz and $f \in H^{m}(\Omega,\mu)$, where $\mu$ is the uniform measure on $\Omega$ and $m > d/2$, then
\bes{
E_{\Lambda,\epsilon}(f) \leq c_{m,d,\Omega} \left ( n^{\theta(d)-m} + \epsilon \right ) \| f \|_{H^{m}(\Omega,\mu)}.
}
}
\prf{
Let $\Lambda = \Lambda^{\mathrm{TP}}_{n}$.  In the first case, since $f$ is defined over the whole of $D$, we may let $p = f_{\Lambda}$ be its orthogonal projection onto $\spn \{ \psi_{\bm{n}} : \bm{n} \in \Lambda \} \subset L^2(D,\nu)$.  Then
\be{
\label{EOPbound}
E_{\Lambda,\epsilon}(f) \leq \| f - f_{\Lambda} \|_{L^\infty(D)} + \epsilon \| f_{\Lambda} \|_{L^2(D,\nu)} \leq \| f - f_{\Lambda} \|_{L^\infty(D)} + \epsilon \| f \|_{L^2(D,\nu)}.
}
It remains to estimate the first term.  For this, we first use the Gagliardo--Nirenberg inequality (see, for example, \cite{HackbuschTensor}) to give
\bes{
\| f - f_{\Lambda} \|_{L^\infty(D)} \leq c_{k,d} \nm{ f - f_{\Lambda} }^{\frac{d}{2k}}_{H^k(D,\nu)} \nm{ f - f_{\Lambda} }^{1-\frac{d}{2k}}_{L^2(D,\nu)},\qquad d < 2k \leq 2 m.
}
We now use the estimate
\bes{
\nm{ f - f_{\Lambda} }_{H^l(D,\nu)} \leq c_{l,m,d} n^{\sigma(l) - m} \| f \|_{H^{m}(D,\nu)},
}
where $\sigma(l) = 0 $ for $l=0$ and $\sigma(l) = 2l - 1/2$ for $l > 0$ (see, for example, \cite[(5.8.11)]{SMSD}).  Hence
\bes{
\| f - f_{\Lambda} \|_{L^\infty(D)} \leq c_{k,m,d} n^{\frac{d(2k-1/2-m)}{2k} - m(1-\frac{d}{2k})} \| f \|_{H^m(D,\nu)} = c_{k,m,d} n^{d(1-\frac{1}{4k}) - m}  \| f \|_{H^m(D,\nu)}.
}
Setting $k = \frac{d+1}{2}$ (odd $d$) or $k = \frac{d+2}{2}$ (even $d$) and substituting into \R{EOPbound} yields the first result for $\Lambda = \Lambda^{\mathrm{TP}}_{n}$.  For the total degree index set $\Lambda = \Lambda^{\mathrm{TD}}_{n}$ we first recall that $\Lambda^{\mathrm{TP}}_{n/d} \subseteq \Lambda^{\mathrm{TP}}_{n}$.  We therefore let $p = f_{\Lambda^{\mathrm{TP}}_{n/d}}$ so that
\bes{
E_{\Lambda^{\mathrm{TD}}_{n},\epsilon}(f) \leq \| f - f_{\Lambda^{\mathrm{TP}}_{n/d}} \|_{L^\infty(D)} + \epsilon \| f \|_{L^2(D,\nu)}.
}
The result for this index set now follows from the previous bound for $\Lambda^{\mathrm{TD}}_{n}$.

Now consider the case where $\Omega$ is Lipschitz and $f \in H^m(\Omega,\mu)$.  We follow the argument of \cite[Prop.\ 5.8]{BADHframespart}.  We first note that there is an extension $g$ of $f$ to $H^{m}(D,\nu)$ satisfying
\bes{
\| g \|_{H^m(D,\nu)} \leq c_{m,d,\Omega} \| f \|_{H^m(\Omega,\mu)}.
}
Now let $p = g_{\Lambda}$ be the orthogonal projection of $g$ onto $\spn \{ \psi_{\bm{n}} : \bm{n} \in \Lambda \}$.  Then
\bes{
\| p \|_{L^2(D,\nu)} \leq \| g \|_{L^2(D,\nu)} \leq c_{m,d,\Omega} \| f \|_{H^{m}(\Omega,\mu)},
}
and, by the previously-derived result,
\bes{
\| f - p \|_{L^\infty(\Omega)} \leq \| g - g_{\Lambda} \|_{L^\infty(D)} \leq c_{m,d,\Omega} n^{\theta(d)-m} \| g \|_{H^m(D,\nu)} \leq c_{m,d,\Omega} n^{\theta(d)-m} \| f \|_{H^m(\Omega,\mu)}.
}
This gives the second result.
}

Unsurprisingly, in scenario (i) one obtains a slightly better error bound, where the constant in the $\epsilon$ term involves the smaller $L^2$-norm as opposed to the $H^m$-norm.  For completeness, we now also consider the hyperbolic cross index set:

\thm{
\label{t:approxerrHC1}
Let $P_{\Lambda}$ be constructed from the tensor Legendre polynomial basis on $L^2(D,\nu)$, where $D = (-1,1)^d$, $\nu$ is the uniform measure on $D$, and $\Lambda = \Lambda_{n}$ is the hyperbolic cross index set \R{LambdaHC} of degree $n$.  If $\Omega \subseteq D$ and $f \in H^{m}(D,\nu)$ for some $m > d/2$, then
\bes{
E_{\Lambda,\epsilon}(f) \leq c_{m,d} \| f \|_{H^{md}(D,\nu)} n^{\frac{\theta(d)-m}{d}}+ \epsilon \| f \|_{L^2(D,\nu)},
}
where $\theta(d)$ is as in \R{thetadef}.  Conversely, if $\Omega \subseteq D$ is Lipschitz and $f \in H^{m}(\Omega,\mu)$, where $\mu$ is the uniform measure on $\Omega$ and $m > d/2$, then
\bes{
E_{\Lambda,\epsilon}(f) \leq c_{m,d,\Omega} \left ( n^{\frac{\theta(d)-m}{d}} + \epsilon \right ) \| f \|_{H^{m}(\Omega,\mu)}.
}
}
\prf{
Let $n^* = \lfloor (n+1)^{1/d} - 1 \rfloor$ and observe that $\Lambda^{\mathrm{TP}}_{n^*} \subseteq \Lambda^{\mathrm{HC}}_{n}$.
We now use the arguments from the proof of the previous theorem.
}

As is to be expected, Theorems \ref{t:approxerr1} and \ref{t:approxerrHC1}, which assume only classical Sobolev regularity, all exhibit the curse of dimensionality.  This can be seen by noting that
\bes{
n^{\theta(d)-m} \asymp N^{\frac{d-m}{d}},\qquad n \rightarrow \infty,
}
for fixed $d$ whenever $\Lambda$ is the total degree or tensor product index set, since in both cases $N = | \Lambda | \asymp n^d$.  Conversely, for the hyperbolic cross index set one has
\bes{
n^{\frac{d-m}{d}} \asymp N^{\frac{d-m}{d}} \left ( \log(N) \right )^{\frac{(m-d)(d-1)}{d}},
}
since in this case $N \asymp n (\log(n))^{d-1}$.

\subsection{Results for the mixed Sobolev spaces $H^{m}_{\mix}$}

Seeking to mitigate the curse of dimensionality when using the hyperbolic cross index set, we now consider the mixed Sobolev spaces $H^{m}_{\mix}(D,\nu)$:

\thm{
\label{t:approxerr2}
Let $P_{\Lambda}$ be constructed from the tensor Legendre polynomial basis on $L^2(D,\nu)$, where $D = (-1,1)^d$ and $\nu$ is the uniform measure on $D$.  If $f \in H^{m}_{\mix}(D,\nu)$ for some $m \geq 1$ then
\bes{
E_{\Lambda,\epsilon}(f) \leq \left \{ \begin{array}{lc} c_{m,d}  \| f \|_{H^{m}_{\mix}(D,\nu)} n^{1-m} + \epsilon \| f \|_{L^2(D,\nu)} & \mbox{$\Lambda = \Lambda^{\mathrm{TP}}_{n}$ or $\Lambda = \Lambda^{\mathrm{TD}}_{n}$} 
\\
c_{m,d}  \| f \|_{H^{m}_{\mix}(D,\nu)} n^{1-m} (\log n)^{\frac{d-1}{2}} + \epsilon \| f \|_{L^2(D,\nu)} & \mbox{$\Lambda = \Lambda^{\mathrm{HC}}_{n}$}
\end{array}
\right . . 
}
Furthermore, if $\Omega$ is compactly contained in $D$, then 
\bes{
E_{\Lambda,\epsilon}(f) \leq \left \{ \begin{array}{lc} c_{m,d,\Omega}  \| f \|_{H^{m}_{\mix}(D,\nu)} n^{1/2-m} + \epsilon \| f \|_{L^2(D,\nu)} & \mbox{$\Lambda = \Lambda^{\mathrm{TP}}_{n}$ or $\Lambda = \Lambda^{\mathrm{TD}}_{n}$} 
\\
c_{m,d,\Omega}   \| f \|_{H^{m}_{\mix}(D,\nu)} n^{1/2-m} (\log n)^{\frac{d-1}{2}} + \epsilon \| f \|_{L^2(D,\nu)} & \mbox{$\Lambda = \Lambda^{\mathrm{HC}}_{n}$}
\end{array}
\right . . 
}
}
\prf{
Since $f \in L^2(D,\nu)$ we may let $p = f_{\Lambda}$ be its orthogonal projection onto $\spn \{ \psi_{\bm{n}} : \bm{n} \in \Lambda \}$.  Then, using \R{1DLegDef}, \R{dDLegDef} and \R{Htilde_char}, we obtain
\ea{
\| f - f_{\Lambda} \|_{L^\infty(D)} &\leq \sum_{\bm{n} \notin \Lambda} \prod^{d}_{k=1} \sqrt{2n_k+1} \left | \ip{f}{\psi_{\bm{n}}}_{L^2(D,\nu)} \right | \nn
\\
& \leq \left ( \sum_{\bm{n} \notin \Lambda} \chi^{\mix}_{\bm{n},m} \left | \ip{f}{\psi_{\bm{n}}}_{L^2(D,\nu)} \right |^2 \right )^{1/2} \left ( \sum_{\bm{n} \notin \Lambda} \frac{\prod^{d}_{k=1} (2n_k+1)}{\chi^{\mix}_{\bm{n},m} } \right )^{1/2}, \nn
\\
& \leq \| f \|_{\tilde{H}^m_{\mix}(D,\nu)}  \left ( \sum_{\bm{n} \notin \Lambda} \frac{\prod^{d}_{k=1} (2n_k+1)}{\chi^{\mix}_{\bm{n},m} } \right )^{1/2}, \nn
}
where $\tilde{H}^m_{\mix}(D,\nu)$ and $\chi^{\mix}_{\bm{n},m}$ are as in \R{tildeHmix} and \R{chidef} respectively.  Observe that
\bes{
\chi^{\mix}_{\bm{n},m} \geq c_{m,d}\left ( \prod^{d}_{k=1} (n_k+1) \right )^{2m},
}
for some constant $c_{m,d}$, and therefore
\be{
\label{coeffstep}
\| f - f_{\Lambda} \|_{L^\infty(D)} \leq c_{m,d} \| f \|_{H^{m}_{\mix}(D,\nu)} \left ( \sum_{\bm{n} \notin \Lambda} \prod^{d}_{k=1} (n_k+1)^{1-2m} \right )^{1/2},
}
where here we also note that $\| f \|_{\tilde{H}^{m}_{\mix}(D,\nu)} \leq c_{m,d} \| f \|_{H^m_{\mix}(D,\nu)}$.  We now specify the index set.  First suppose that $\Lambda = \Lambda^{\mathrm{TP}}_n$.  Let $[d]$ denote the set of ordered tuples with entries in $\{1,\ldots,d\}$.  Then
\eas{
\sum_{\bm{n} \notin \Lambda} \prod^{d}_{k=1} (n_k+1)^{1-2m}  & = \sum_{\sigma \in [d] } \sum^{n}_{\substack{n_k =0 \\ k \notin \sigma}} \sum_{\substack{n_k > n \\ k \in \sigma}}  \prod^{d}_{k=1} (n_k+1)^{1-2m}
\\
&  = \sum_{\sigma \in [d]} \left (1 + \sum^{n}_{l=1} l^{1-2m} \right )^{d - |\sigma |} \left ( \sum_{l \geq n} l^{1-2m} \right )^{|\sigma|} \leq c_{m,d} n^{2-2m}.
}
Substituting into \R{coeffstep} now gives the result for $\Lambda = \Lambda^{\mathrm{TP}}_n$.  Moreover, the result for the total degree index set now follows as well, after noting that $\Lambda^{\mathrm{TD}}_{n} \supseteq \Lambda^{\mathrm{TP}}_{n/d}$.  Finally, for the hyperbolic cross index $\Lambda = \Lambda^{\mathrm{HC}}_{n}$ set we use, for example, \cite[Lem.\ 2.30]{BAthesis} to get
\bes{
\sum_{\bm{n} \notin \Lambda} \prod^{d}_{k=1} (n_k+1)^{1-2m} \leq c_{m,d} n^{2-2m} (\log(n))^{d-1},
}
as required.  

It remains to consider the case where $\Omega$ is compactly contained in $D$.  We first recall that univariate Legendre polynomials are uniformly bounded in compact subintervals of $(-1,1)$:
\bes{
| \psi_{n}(y)| \leq c_{r},\qquad -1+r \leq y \leq 1-r,\quad \forall n \in \bbN_0,\ 0 < r < 1,
}
for some $c_r>0$.  Hence $\| \phi_{\bm{n}} \|_{L^\infty(\Omega)} \leq c_{\Omega}$, $\forall \bm{n} \in \bbN^d_0$.
Letting $p = f_{\Lambda}$ and arguing as before, we get
\eas{
\| f - f_{\Lambda} \|_{L^\infty(\Omega)} &\leq c_{m,d,\Omega} \| f \|_{H^m_{\mix}(D,\nu)} \left ( \sum_{\bm{n} \notin \Lambda} \frac{1}{\chi^{\mix}_{\bm{n},m}} \right )^{1/2} 
 \\
& \leq c_{m,d,\Omega} \| f \|_{H^m_{\mix}(D,\nu)}\left ( \sum_{\bm{n} \notin \Lambda} \prod^{d}_{k=1} (n_k+1)^{-2m}\right )^{1/2}.
}
We now proceed in the same way, replacing the exponent $1-2m$ by $-2m$ throughout.
}

\section{Sample complexity}\label{s:sampcomp}
In this section we consider the efficiency of the approximation.  In view of Theorem \ref{t:acc_stab_cond} this requires estimating the constant $C_{\Upsilon,\Lambda,\epsilon}$ defined in \R{Cconst}.   Our main results are twofold.  First, in \S \ref{ss:Nikolskii} we show that when the sample points are drawn independently according to a suitable measure on $\Omega$ then the sample complexity can always be related to the constant of a certain Nikolskii inequality for the polynomial space $P_{\Lambda}$.  Second, in \S \ref{ss:lambdarectangle} we show that for domains satisfying a suitable property this constant is at most log-quadratic in the dimension $N$ of the polynomial space $P_{\Lambda}$.

\subsection{The constant $C_{\Upsilon,\Lambda}$}
It is difficult to analyze $C_{\Upsilon,\Lambda,\epsilon}$ directly, since it is defined in terms of the singular values and singular vectors of the matrix $\bm{A}$.  In order to provide concrete bounds, we now consider
\bes{
C_{\Upsilon,\Lambda}= \sup \left \{ \nm{p}_{L^2(\Omega,\mu)} : p \in P_{\Lambda}, \frac1M \sum_{\bm{y} \in \Upsilon} | p(\bm{y}) |^2 = 1 \right \}.
}
{Note that $C_{\Upsilon,\Lambda}$ depends only on the samples $\Upsilon$ and the space $P_{\Lambda}$.  Unlike $C_{\Upsilon,\Lambda,\epsilon}$, it is independent of functions $\phi_{\bm{n}}$, $\bm{n} \in \Lambda$, used to span this space and consequently the domain $D$ as well.  We also have the following:}

\lem{
\label{l:CUB}
Let $C_{\Upsilon,\Lambda,\epsilon}$ be as in \R{Cconst}.  Then $C_{\Upsilon,\Lambda,\epsilon} \leq C_{\Upsilon,\Lambda}$ and moreover $C_{\Upsilon,\Lambda} = C_{\Upsilon,\Lambda,\epsilon}$ whenever the minimum singular value of $\bm{A} = \bm{A}_{\Upsilon,\Lambda}$ satisfies $\sigma_{\min}(\bm{A}) > \epsilon$.
}
\prf{
Recall that $C_{\Upsilon,\Lambda,\epsilon} $ is the maximum of $C'_{\Upsilon,\Lambda,\epsilon}$ and $C''_{\Upsilon,\Lambda,\epsilon}$.  Consider $C'_{\Upsilon,\Lambda,\epsilon}$.  Let $\bm{b} \in \bbC^M$, $\nm{\bm{b}}_2 = 1$, and notice that $\cT_{\Lambda} (\bm{A}_{\Upsilon,\Lambda,\epsilon})^{\dag} \bm{b} \in P_{\Lambda}$.  Hence, by the definition of $C_{\Upsilon,\Lambda}$ and \R{SUpsilon}, we have
\bes{
\nm{\cT_{\Lambda} (\bm{A}_{\Upsilon,\Lambda,\epsilon})^{\dag} \bm{b}}_{L^2(\Omega,\mu)} \leq C_{\Upsilon,\Lambda} \nm{\cS_{\Upsilon} \cT_{\Lambda} (\bm{A}_{\Upsilon,\Lambda,\epsilon})^{\dag} \bm{b}}_{2}.
}
Note that $\bm{A} = \bm{A}_{\Upsilon,\Lambda} = \cS_{\Upsilon} \cT_{\Lambda}$.  By standard properties of the SVD,
\bes{
\nm{\cS_{\Upsilon} \cT_{\Lambda} (\bm{A}_{\Upsilon,\Lambda,\epsilon})^{\dag} \bm{b}}_{2} \leq \nm{\bm{b}}_{2} = 1,
}
and therefore $\nm{\cT_{\Lambda} (\bm{A}_{\Upsilon,\Lambda,\epsilon})^{\dag} \bm{b}}_{L^2(\Omega,\mu)} \leq C_{\Upsilon,\Lambda}$.  Since $\bm{b}$ was arbitrary we get $C'_{\Upsilon,\Lambda,\epsilon} \leq C_{\Upsilon,\Lambda}$.

On the other hand, suppose that $\sigma_{\min}(\bm{A}) > \epsilon$.  Let $p = \cT_{\Lambda} \bm{c} \in P_{\Lambda}$ with $\nm{\cS_{\Upsilon} p}_{2} = 1$.  Let $\bm{b} = \cS_{\Upsilon} p$ and write $p = \cT_{\Lambda} \bm{c}$.  Then
\bes{
C'_{\Upsilon,\Lambda,\epsilon} \geq \nm{\cT_{\Lambda} (\bm{A}_{\Upsilon,\Lambda,\epsilon})^{\dag} \cS_{\Upsilon} \cT_{\Lambda} \bm{c} }_{L^2(\Omega,\mu)} = \nm{\cT_{\Lambda} (\bm{A}_{\Upsilon,\Lambda,\epsilon})^{\dag} \bm{A}_{\Upsilon,\Lambda} \bm{c} }_{L^2(\Omega,\mu)} = \nm{\cT_{\Lambda} \bm{c}}_{L^2(\Omega,\mu)} = \| p \|_{L^2(\Omega,\mu)},
}
where in the third step we use the fact that $\bm{A}_{\Upsilon,\Lambda}$ is full rank.  Hence, since $p$ was arbitrary, we get $C'_{\Upsilon,\Lambda,\epsilon} \geq C_{\Upsilon,\Lambda} $, and therefore $C'_{\Upsilon,\Lambda,\epsilon}  = C_{\Upsilon,\Lambda}$ in this case. 

Finally, consider $C''_{\Upsilon,\Lambda,\epsilon}$.  Let $\bm{d} \in \bbC^N$, $\nm{\bm{d}}_{2} = 1$.  Then
\eas{
\nm{\cT_{\Lambda} \bm{d} - \cT_{\Lambda}(\bm{A}_{\Upsilon,\Lambda,\epsilon})^{\dag} \cT_{\Lambda} \bm{d} }_{L^2(\Omega,\mu)} & \leq C_{\Upsilon,\Lambda} \nm{\cS_{\Upsilon} \left ( \cT_{\Lambda} \bm{d} - \cT_{\Lambda}(\bm{A}_{\Upsilon,\Lambda,\epsilon})^{\dag} \bm{A}_{\Upsilon,\Lambda} \bm{d} \right ) }_{2}
\\
& = C_{\Upsilon,\Lambda}\nm{\left ( \bm{A}_{\Upsilon,\Lambda} - \bm{A}_{\Upsilon,\Lambda} (\bm{A}_{\Upsilon,\Lambda,\epsilon})^{\dag} \bm{A}_{\Upsilon,\Lambda} \right ) \bm{d} }_{2}
\\
& = C_{\Upsilon,\Lambda}\nm{\bm{U} \left ( \bm{\Sigma} - \bm{\Sigma}_{\epsilon} \right ) \bm{V}^* \bm{d} }_{2}
 \leq C_{\Upsilon,\Lambda}\epsilon \nm{\bm{d}}_{2}.
}
Since $\bm{d}$ was arbitrary, we get $C''_{\Upsilon,\Lambda,\epsilon} \leq C_{\Upsilon,\Lambda}$ as required.  On the other hand, if $\sigma_{\min}(\bm{A}) > \epsilon$ then $\bm{\Sigma} - \bm{\Sigma}_{\epsilon} = 0$.  Hence $C''_{\Upsilon,\Lambda,\epsilon} = 0$.
}

\subsection{Random sampling for compact domains and Nikolskii inequalities}\label{ss:Nikolskii}
We now show that $C_{\Upsilon,\Lambda}$ can be bounded using the constant of a suitable Nikolskii inequality for the space $P_{\Lambda} \subset L^2(\Omega,\mu)$.  To this end, let $N(P_{\Lambda},\Omega,\mu)$ be the smallest positive number in the $(L^2(\Omega,\mu),L^\infty(\Omega))$-Nikolskii inequality
\be{
\label{NikolskiiConst}
\| p \|_{L^\infty(\Omega)} \leq \rN(P_{\Lambda},\Omega,\mu)\| p \|_{L^2(\Omega,\mu)},\qquad \forall p \in P_{\Lambda}.
}
Then we have the following result:

\thm{
\label{t:SampNikolskii}
Let $0 < \delta,\gamma < 1$ and suppose that $\bm{y}_1,\ldots,\bm{y}_M$ are independent and randomly drawn according to the probability measure $\mu$ defined by \R{mu_def}.  If
\bes{
M \geq \left(\rN(P_{\Lambda},\Omega,\mu)\right )^2  \left ( (1-\delta) \log(1-\delta)+\delta \right )^{-1} \log(N/\gamma),
}
where $N = | \Lambda |$ and $\rN(P_{\Lambda},\Omega,\mu)$ is the constant of the Nikolskii inequality \R{NikolskiiConst}, then with probability at least $1-\gamma$ the quantity $C_{\Upsilon,\Lambda}$ satisfies
\bes{
C_{\Upsilon,\Lambda}\leq \frac{1}{\sqrt{1-\delta}}.
}
}
\prf{
Our proof is based on {essentially the same arguments as those used in previous works (see, for example, \cite{DavenportEtAlLeastSquares})}.  First, let $\{ \Phi_{\bm{n}} \}_{\bm{n} \in \Lambda}$ be an orthonormal basis for $P_{\Lambda}$ in $L^2(\Omega,\mu)$.  Let $p \in P_{\Lambda}$ be arbitrary and write $p = \sum_{\bm{n} \in \Lambda} c_{\bm{n}} \Phi_{\bm{n}}$,
so that
\bes{
\| p \|^2_{L^2(\Omega,\mu)} = \int_{\Omega} | p(\bm{y}) |^2 \D \mu(\bm{y}) = \nm{\bm{c}}^2_{2},
}
where $\bm{c} = (c_{\bm{n}} )_{\bm{n} \in \Lambda}$, and
$
\frac1M \sum_{\bm{y} \in \Upsilon} | p(\bm{y}) |^2 = \bm{c}^* \bm{B} \bm{c},
$
where $\bm{B} \in \bbC^{N \times N}$ is the self-adjoint matrix with
\bes{
(\bm{B})_{\bm{m},\bm{n}} = \frac1M \sum_{\bm{y} \in \Upsilon} \overline{\Phi_{\bm{m}}(\bm{y})} \Phi_{\bm{n}}(\bm{y}),\qquad \bm{m},\bm{n} \in \Lambda.
}
It follows that
$
C_{\Upsilon,\Lambda} = 1/\sqrt{\lambda_{\min}(\bm{B})}, 
$
where $\lambda_{\min}(\bm{B})$ is the minimal eigenvalue of $\bm{B}$.  We estimate this quantity by writing it in the usual way as the sum of random matrices:
\bes{
\bm{B} = \sum^{M}_{m=1} \bm{X}_{m},\qquad \bm{X}_{m} = \left \{ \frac1M \overline{\Phi_{\bm{m}}(\bm{y}_m)} \Phi_{\bm{n}}(\bm{y}_m) \right \}_{\bm{m},\bm{n} \in \Lambda}.
}
By construction, these matrices are independent, nonnegative definite and satisfy $\bbE(\bm{X}_m) = \frac1M \bm{I}$, where $\bm{I}$ is the identity matrix.  Moreover, for any $\bm{c} \in \bbC^N$ we have
\bes{
\bm{c}^* \bm{X}_{m} \bm{c} = \frac{1}{M} \left | \sum_{\bm{n} \in \Lambda} c_{\bm{n}} \Phi_{\bm{n}}(\bm{y}_m) \right |^2 \leq \frac{(\rN(P_{\Lambda},\Omega,\mu))^2}{M} \nm{ \sum_{\bm{n} \in \Lambda} c_{\bm{n}} \Phi_{\bm{n}} }^2_{L^2(\Omega,\mu)} = \frac{(\rN(P_{\Lambda},\Omega,\mu))^2 }{M} \nm{\bm{c}}^2_2.
}
The Matrix Chernoff bound (see, for example, \cite[Thm.\ 1.1]{TroppUserFriendly}) now gives 
\bes{
\bbP \left ( \lambda_{\min}(\bm{X}) \leq (1-\delta) \right ) \leq N \exp \left ( - \frac{(1-\delta) \log(1-\delta)+\delta}{M^{-1} (\rN(P_{\Lambda},\Omega,\mu))^2 } \right ).
}
Setting the right-hand side equal to $\gamma$ and rearranging yields the result.
}

This leads to the following result on accuracy of the truncated SVD least-squares approximation:

\cor{
\label{c:acc_stab_prob}
Let $0 < \delta,\gamma < 1$ and suppose that $\bm{y}_1,\ldots,\bm{y}_M$ are independent and randomly drawn according to the probability measure $\mu$ defined by \R{mu_def}.  Let
\bes{
M \geq \left(\rN(P_{\Lambda},\Omega,\mu)\right )^2  \left ( (1-\delta) \log(1-\delta)+\delta \right )^{-1} \log(N/\gamma),
}
where $N = | \Lambda |$ and $\rN(P_{\Lambda},\Omega,\mu)$ is the constant of the Nikolskii inequality \R{NikolskiiConst}.  Then with probability at least $1-\gamma$ the truncated SVD least-squares approximation $f_{\Upsilon,\Lambda,\epsilon}$ of $f \in L^\infty(\Omega)$ satisfies
\bes{
\| f - f_{\Upsilon,\Lambda,\epsilon} \|_{L^2(\Omega,\mu)} \leq \left ( 1 + \frac{1}{\sqrt{1-\delta}} \right ) E_{\Lambda,\epsilon}(f),
}
where $E_{\Lambda,\epsilon}(f)$ is as in \R{ELambdaDef}.
}

\subsection{The $\lambda$-rectangle property and log-quadratic sample complexity}\label{ss:lambdarectangle}

We now consider $N(P_{\Lambda},\Omega,\mu)$.  Estimating this constant for general irregular domains in arbitrarily-many dimensions is an open problem.  We shall not attempt to resolve it in full generality here (see \S \ref{s:conclusion} for some further discussion).  Instead, we show that this constant is at most log-quadratic for a large class of irregular domains whenever $\mu$ is the uniform measure.

\begin{figure}
\begin{center}
\begin{tabular}{ccccc}
\begin{tikzpicture}[>=stealth]
\filldraw[fill=green!20,draw=green!50!black, thick] 
(0,0) -- (3,0) -- (3,1.5) -- (1.5,1.5) -- (1.5,3) -- (0,3) -- (0,0);
\end{tikzpicture}
&\hspace{3pc}&
\begin{tikzpicture}[>=stealth]
\filldraw[fill=green!20,draw=green!50!black, thick] 
(0,0) -- (3,0) -- (3,1.5) -- (1.5,3) -- (0,3) -- (0,0);
\end{tikzpicture}
&\hspace{3pc}&
\begin{tikzpicture}[>=stealth]
\filldraw[fill=green!20,draw=green!50!black, thick] 
(0,0) -- (3,0) -- (3,3) -- (0,3) -- (0,0) ;
\filldraw[fill=white,draw=green!50!black, thick] 
(1.5,1.5) circle (0.75);
\end{tikzpicture}
\\
(a) && (b) && (c)
\end{tabular}
\end{center}
\vspace*{-5mm}
\caption{Examples of domains that have the $\lambda$-rectangle property: (a) is an $L$-shaped domain, (b) is a domain with a linear constraint, and (c) is a domain with an exclusion.}
\label{f:lambdarectangle2}
\end{figure}
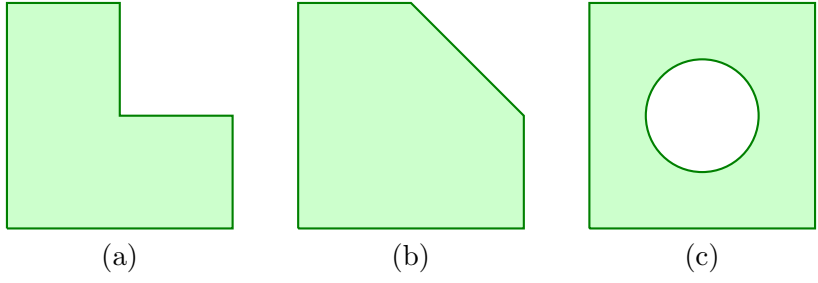

The types of domain we now consider are those satisfying the so-called following property:

\defn{
[$\lambda$-rectangle property]
A compact domain $\Omega$ has the $\lambda$-rectangle property for some $0<\lambda<1 $ if it can be written as a (possibly overlapping and uncountable) union
\bes{
\Omega = \bigcup_{R \in \cR} R,
}
of hyperrectangles $R$ satisfying
\bes{
 \inf_{R \in \cR} \mathrm{Vol}(R) = \lambda \mathrm{Vol}(\Omega).
}
}
There are many domains of interest that have this property.  We now list several examples:

\begin{itemize}
\item \textit{$L$-shaped domains.} These are unions of two rectangles, so clearly have this property.
\item \textit{Domains with linear constraints.}  The domain
\bes{
\Omega = \{ - 1 \leq y_1 , y_2 \leq 1, y_1 + y_2 \leq 1 \},
}
along with its various higher-dimensional generalizations, can be expressed as
\bes{
\Omega = \bigcup_{x \in [0,1]} R_x,\qquad R_x = [-1,x] \otimes [-1,1-x].
}
Hence it has the $\lambda$-rectangle property with $\lambda = 4/7$.  Note that such domains can occur in problems such as surrogate forwards model construction in parameter studies; for instance, whenever two parameters $y_1$ and $y_2$, rather than being independent, satisfy a (possibly \textit{a priori} unknown) linear relation.
\item \textit{Domains with exclusions.} The domain
\bes{
\Omega = \{ -1 \leq y_1,y_2 \leq 1, y^2_1 + y^2_2 \geq 1/2 \},
}
along with various generalizations, also satisfies the $\lambda$-rectangle property. Note that such domains correspond to practical scenarios where, due to certain physical constraints, $f(\bm{y})$ can only be evaluated for $\bm{y}$ not too close to zero.
\end{itemize}
See Fig.\  \ref{f:lambdarectangle2} for illustrations.  On the other hand, there are a number of notable domains that do not have this property.  These include the unit Euclidean ball $\{ \bm{y} \in \bbR^d : \nm{\bm{y}}_2 \leq 1 \}$ and the simplex $\{ \bm{y} \in \bbR^d : 0 \leq y_1,\ldots,y_{d-1} \leq 1, 0 \leq y_d \leq 1 - (y_1+\ldots+y_{d-1}) \}$.  See \S \ref{s:conclusion} for additional details.

\thm{
\label{t:Nikolskii_lambda_rect}
Suppose that $\Omega \subseteq (-1,1)^d$ has the $\lambda$-rectangle property and let $P_{\Lambda}$ be constructed from the tensor Legendre polynomial basis on $(-1,1)^d$, where $\Lambda$ is any lower set (see Definition \ref{d:lowerset}) of cardinality $| \Lambda | = N$.  Let $\mu$ be the uniform probability measure on $\Omega$ and $\rN(P_{\Lambda},\Omega,\mu)$ be the constant in the Nikolskii inequality \R{NikolskiiConst}.  Then
\bes{
\left(\rN(P_{\Lambda},\Omega,\mu) \right)^2 \leq \frac{N^2}{\lambda} .
}
}
\prf{
We first claim that $P_{\Lambda} = \spn \{ \phi_{\bm{n}} : \bm{n} \in \Lambda \}$ coincides with the space
\bes{
\bbP_{\Lambda} = \spn \left \{ \bm{y} \in \Omega \mapsto \bm{y}^{\bm{n}} : \bm{n} = (n_1,\ldots,n_d) \in \Lambda \right \}.
}
Here we use the notation $\bm{y}^{\bm{n}} = y^{n_1}_1 \cdots y^{n_d}_{d}$.
Since $\phi_{\bm{n}}(\bm{y})$ is a tensor Legendre polynomial we have
\bes{
\phi_{\bm{n}}(\bm{y}) = \prod^{d}_{k=1} \psi^{(k)}_{n_k}(y_k) = \prod^{d}_{k=1} \left ( \sum^{n_k}_{m_k=0} a_{m_k,n_k} y^{n_k}_k \right ) = \sum^{n_1}_{m_1=0} \cdots \sum^{n_d}_{m_d=0} a_{\bm{m},\bm{n}} \bm{y}^{\bm{m}},
}
where $a_{m_k,n_k}$ are the coefficients of $\psi^{(k)}_{n_k}$ in the monomial basis and $a_{\bm{m},\bm{n}} = \prod^{d}_{k=1} a_{m_k,n_k}$.  Since $m_k \leq n_k$ for all $k$, it follows from the lower set assumption that $\bm{m} \in \Lambda$ and therefore $\phi_{\bm{n}} \in \bbP_{\Lambda}$.  Hence $P_{\Lambda} \subseteq \bbP_{\Lambda}$.  In a similar manner, one also finds that $\bm{y} \mapsto \bm{y}^{\bm{n}}$ is in $P_{\Lambda}$, and therefore $\bbP_{\Lambda} \subseteq P_{\Lambda}$, as required.

Now let $p \in \bbP_{\Lambda}$ and $\bm{y} \in \Omega$ with $\bm{y} \in R$ for some $R \in \cR$.  Define the uniform measure on $R$ as
\bes{
\D \tilde{\mu}(\bm{y}) = \frac{1}{\mathrm{Vol}(R)} \D \bm{y},
}
and note that $| p(\bm{y}) | \leq \rN(\bbP_{\Lambda},R,\tilde{\mu}) \nm{p}_{L^2(R,\tilde{\mu})}$,
where $\rN(\bbP_{\Lambda},R,\tilde{\mu})$ is the Nikolskii constant for the space $\bbP_{\Lambda}$ with respect to $L^2(R,\tilde{\mu})$.  It  is known that $(\rN(\bbP_{\Lambda},R,\tilde{\mu}))^2 \leq N^2$ \cite[Thm.\ 6]{MiglioratiJAT}.  Also
\bes{
\nm{p}^2_{L^2(R,\tilde{\mu})} = \frac{1}{\mathrm{Vol}(R)} \int_{R} |p(\bm{y})|^2 \D \bm{y} \leq \frac{\mathrm{Vol}(\Omega)}{\mathrm{Vol}(R)} \int_{\Omega} |p(\bm{y})|^2 \D \mu(\bm{y}) \leq \frac{1}{\lambda} \nm{p}^2_{L^2(\Omega,\mu)}.
}
Hence
$
| p(\bm{y}) |^2 \leq \frac{N^2}{\lambda } \| p \|^2_{L^2(\Omega,\mu)}.
$
Since $\bm{y} \in \Omega$ and $p \in \bbP_{\Lambda}$ were arbitrary, we now get the result.
}

Combining this with Corollary \ref{c:acc_stab_prob} now gives the following:

\cor{
\label{c:sample_complexity}
Suppose that $\Omega \subseteq (-1,1)^d$ has the $\lambda$-rectangle property and let $P_{\Lambda}$ be constructed from the tensor Legendre polynomial basis on $(-1,1)^d$, where $\Lambda$ is any lower set of cardinality $| \Lambda | = N$.  Let $0 < \delta,\gamma < 1$ and suppose that $\bm{y}_1,\ldots,\bm{y}_M$ are independent and randomly drawn according to the uniform probability measure on $\Omega$.  If
\bes{
M \geq N^2 \lambda^{-1} \left ( (1-\delta) \log(1-\delta)+\delta \right )^{-1} \log(N/\gamma),
}
where $N = | \Lambda |$, then with probability at least $1-\gamma$ the truncated SVD least-squares approximation $f_{\Upsilon,\Lambda,\epsilon}$ of $f \in L^{\infty}(\Omega,\mu)$ satisfies
\bes{
\| f - f_{\Upsilon,\Lambda,\epsilon} \|_{L^2(\Omega,\mu)} \leq \left ( 1 + \frac{1}{\sqrt{1-\delta}} \right ) E_{\Lambda,\epsilon}(f),
}
where $E_{\Lambda,\epsilon}(f)$ is as in \R{ELambdaDef}.
}

\section{Truncated estimators and $L^2$-error bounds}\label{s:L2bounds}

The error bounds proved in Corollary \ref{c:acc_stab_prob} and elsewhere have the limitation of relating (in probability) the $L^2$-norm of the error to an approximation error $E_{\Lambda,\epsilon}(f)$ measured in the $L^\infty$-norm.  In this penultimate section we show that it is possible to bound the $L^2$-norm of a related estimator in expectation in terms of the $L^2$-norm approximation error
\be{
\label{tildeELambdaDef}
\tilde{E}_{\Lambda,\epsilon}(f) = \inf \left \{ \nm{f - p }_{L^2(\Omega,\mu)} + \epsilon \nm{p}_{L^2(D,\nu)} : \ p \in P_{\Lambda} \right \}.
}
We follow the approach of \cite{DavenportEtAlLeastSquares}.  First, suppose that $f \in L^{\infty}(\Omega,\mu)$ and let $L \geq 0$ be such that $\| f \|_{L^\infty(\Omega) } \leq L$.
Now define the truncation operator
\bes{
T_{L}(g)(\bm{y}) = \sgn(g(\bm{y})) \min \{ | g(\bm{y}) | , L \},
}
where $\sgn(z)$ denotes the complex sign of $z \in \bbC$.  If $f_{\Upsilon,\Lambda,\epsilon}$ is the truncated SVD least-squares approximation we now consider the new approximation
\be{
\label{ftruncest}
f_{\Upsilon,\Lambda,\epsilon,L} = T_{L} \left ( f_{\Upsilon,\Lambda,\epsilon} \right ).
}
Our main result is now the following:

\thm{
\label{t:truncated_est}
Let $0 < \delta,\gamma < 1$ and $f \in L^\infty(\Omega)$ with $\| f \|_{L^\infty(\Omega) } \leq L$ for some $L \geq 0$.  Let $\bm{y}_1,\ldots,\bm{y}_M$ be independent and randomly drawn according to $\mu$ and $f_{\Upsilon,\Lambda,\epsilon,L}$ be as in \R{ftruncest}.  If $\tilde{E}_{\Lambda,\epsilon}(f)$ is as in \R{tildeELambdaDef} and
\be{
\label{meascondtrunc}
M \geq \left(\rN(P_{\Lambda},\Omega,\mu)\right )^2  \left ( (1-\delta) \log(1-\delta)+\delta \right )^{-1} \log(N/\gamma),
}
where $N = | \Lambda |$ and $\rN(P_{\Lambda},\Omega,\mu)$ is the constant of the Nikolskii inequality \R{NikolskiiConst}, then
\bes{
\bbE \left ( \| f - f_{\Upsilon,\Lambda,\epsilon,L} \|^2_{L^2(\Omega,\mu)} \right ) \leq 3  \frac{2-\delta}{1-\delta} \left(\tilde{E}_{\Lambda,\epsilon}(f)\right)^2 + 4 L^2 \gamma.
}
}
\prf{
The proof is based on \cite[Thm.\ 2]{DavenportEtAlLeastSquares}.  Let $E$ be the event $C_{\Upsilon,\Lambda,\epsilon} \leq \frac{1}{\sqrt{1-\delta}}$, where $C_{\Upsilon,\Lambda,\epsilon}$ is as in \R{Cconst}.  Lemma \ref{l:CUB}, Theorem \ref{t:SampNikolskii} and the measurement condition \R{meascondtrunc} give that $\bbP(E^c) \leq \gamma$.
Now let $\D\mu$ be the uniform measure on $\Omega$ and $\D \mu_{M} = \D \mu \otimes \cdots \otimes \D \mu$ be the probability measure of the draw $\bm{y}_1,\ldots,\bm{y}_M$.  Then
\ea{
\bbE \left ( \| f - f_{\Upsilon,\Lambda,\epsilon,L} \|^2_{L^2(\Omega,\mu)} \right ) \D \mu_{M} &= \int_{E}  \| f - f_{\Upsilon,\Lambda,\epsilon,L} \|^2_{L^2(\Omega,\mu)} \D \mu_{M} + \int_{E^c}  \| f - f_{\Upsilon,\Lambda,\epsilon,L} \|^2_{L^2(\Omega,\mu)} \D \mu_{M} \nn
\\
& \leq \int_{E}  \| f - f_{\Upsilon,\Lambda,\epsilon,L} \|^2_{L^2(\Omega,\mu)} \D \mu_{M} + 4 L^2 \gamma. \label{truncstep1}
}
It remains to bound the first term.  Assume the event $E$ occurs and let $p \in P_{\Lambda}$ be such that $\| f - p \|_{L^2(\Omega,\mu)} + \epsilon \| p \|_{L^2(D,\nu)} = \tilde{E}_{\Lambda,\epsilon}(f)$ (it is straightforward to show that such a minimizer exists, since $P_{\Lambda}$ is finite dimensional).  Then, arguing as in the proof of Theorem \ref{t:acc_stab_cond} and using the fact that $C_{\Upsilon,\Lambda,\epsilon} \leq \frac{1}{\sqrt{1-\delta}}$, we have
\eas{
\| f - f_{\Upsilon,\Lambda,\epsilon,L} \|^2_{L^2(\Omega,\mu)} &\leq \left ( \| f - p \|_{L^2(\Omega,\mu)} + C'_{\Upsilon,\Lambda,\epsilon}\| \cS_{\Upsilon}(f-p) \|_{2} + \epsilon C''_{\Upsilon,\Lambda,\epsilon} \| p \|_{L^2(D,\nu)} \right )^2
\\
&\leq 3 \| f - p \|^2_{L^2(\Omega,\mu)} + \frac{3}{1-\delta} \| \cS_{\Upsilon}(f-p) \|^2_{2} + \frac{3\epsilon^2}{1-\delta} \| p \|^2_{L^2(D,\nu)}.
}
Hence
\bes{
 \int_{E}  \| f - f_{\Upsilon,\Lambda,\epsilon,L} \|^2_{L^2(\Omega,\mu)} \D \mu_{M} \leq 3 \| f - p \|^2_{L^2(\Omega,\mu)} +\frac{3}{1-\delta} \bbE \left ( \| \cS_{\Upsilon}(f-p) \|^2_{2} \right ) + \frac{3 \epsilon^2}{1-\delta}  \| p \|^2_{L^2(D,\nu)}.
}
Observe that $\bbE \left ( \| \cS_{\Upsilon}(f-p) \|^2_{2} \right ) = \bbE | f(\bm{y}) - p(\bm{y}) |^2 = \| f - p \|^2_{L^2(\Omega,\mu)}$.
Therefore we obtain
\eas{
 \int_{E}  \| f - f_{\Upsilon,\Lambda,\epsilon,L} \|^2_{L^2(\Omega,\mu)} \D \mu_{M} &\leq 3 \frac{2-\delta}{1-\delta} \left ( \| f - p \|^2_{L^2(\Omega,\mu)} + \epsilon^2 \| p \|^2_{L^2(D,\nu)} \right ) \leq 3 \frac{2-\delta}{1-\delta} \left ( \tilde{E}_{\Lambda,\epsilon}(f) \right )^2.
}
Substituting this into \R{truncstep1} now gives the result.
}

Much like in \S \ref{s:approx_err}, we can establish bounds for $\tilde{E}_{\Lambda,\epsilon}(f)$ under different regularity conditions:

\thm{
\label{t:approxerr1L2}
Let $P_{\Lambda}$ be constructed from the tensor Legendre polynomial basis on $L^2(D,\nu)$, where $D = (-1,1)^d$, $\nu$ is the uniform measure on $D$, and $\Lambda = \Lambda_{n}$ is either the tensor product \R{LambdaTP} or total degree \R{LambdaTD} index set of degree $n$.  If $\Omega \subseteq D$ and $f \in H^{m}(D,\nu)$ for some $m \geq 1$, then
\bes{
\tilde{E}_{\Lambda,\epsilon}(f) \leq c_{m,d} \| f \|_{H^{m}(D,\nu)} n^{-m}+ \epsilon \| f \|_{L^2(D,\nu)}.
}
Conversely, if $\Omega \subseteq (-1,1)^d$ is Lipschitz and $f \in H^{m}(\Omega,\mu)$ where $\mu$ is the uniform measure on $\Omega$ and $m \geq 1$, then
\bes{
\tilde{E}_{\Lambda,\epsilon}(f) \leq c_{m,d,\Omega} \left ( n^{-m} + \epsilon \right ) \| f \|_{H^{m}(\Omega,\mu)}.
}
Finally, if $\Lambda = \Lambda^{\mathrm{HC}}_n$ is the hyperbolic cross index set \R{LambdaHC} then the same results hold with $n^{-m}$ replaced by $n^{-m/d}$.
}
\prf{
As in the proof of Theorem \ref{t:approxerr1}, if $f \in H^{m}(D,\nu)$ we let $f_{\Lambda}$ be the orthogonal projection of $f$ onto $\spn \{ \psi_{\bm{n}} : \bm{n} \in \Lambda \} \subset L^2(D,\nu)$.  Then by Parseval's identity and \R{Htilde_char},
\eas{
\| f - f_{\Lambda} \|^2_{L^2(D,\nu)} = \sum_{\bm{n} \notin \Lambda} \left | \ip{f}{\psi_{\bm{n}}}_{L^2(D,\nu)} \right |^2  &\leq \frac{1}{\min_{\bm{n} \notin \Lambda} \{ \chi_{\bm{n},m} \} }  \sum_{\bm{n} \notin \Lambda} \chi_{\bm{n},m} \left | \ip{f}{\psi_{\bm{n}}}_{L^2(D,\nu)} \right |^2 
\\
& \leq  \frac{1}{\min_{\bm{n} \notin \Lambda} \{ \chi_{\bm{n},m} \} } \| f \|^2_{\tilde{H}^m(D,\nu)}.
}
It remains to bound $\min_{\bm{n} \notin \Lambda} \{ \chi_{\bm{n},m} \} $ for the three index sets.  Using \R{chidef}, we first observe that
\bes{
\chi_{\bm{n},m} = \sum_{|\bm{j}|_1 \leq m} \prod^{d}_{k=1} (n_k(n_k+1))^{j_k} \geq | \bm{n} |^{2m}_{\infty} > n^{2m},\qquad \bm{n} \notin \Lambda^{\mathrm{TP}}_{n}.
}
Similarly, for the total degree index set
\bes{
\chi_{\bm{n},m} \geq c_{m,d} | \bm{n} |^{2m}_{1} > c_{m,d} n^{2m},\qquad \bm{n} \notin \Lambda^{\mathrm{TD}}_{n},
}
and for the hyperbolic cross
\bes{
\chi_{\bm{n},m} \geq c_{m,d} | \bm{n} |^{2m/d}_{\mathrm{hc}} > c_{m,d} n^{2m/d},\qquad \bm{n} \notin \Lambda^{\mathrm{HC}}_{n}.
}
This gives the first result.  For the second result, we argue as in the proof of Theorem \ref{t:approxerr1} to construct an extension $g \in H^{m}(D,\nu)$ of $f$, and then use the previously-derived bounds.
}

\thm{
\label{t:approxerr2L2}
Let $P_{\Lambda}$ be constructed from the tensor Legendre polynomial basis on $L^2(D,\nu)$, where $D = (-1,1)^d$ and $\nu$ is the uniform measure on $D$.  If $f \in H^{m}_{\mix}(D,\nu)$ for some $m \geq 1$ then
\bes{
\tilde{E}_{\Lambda,\epsilon}(f) \leq c_{m,d}  \| f \|_{H^{m}_{\mix}(D,\nu)} n^{-m} + \epsilon \| f \|_{L^2(D,\nu)},
}
when $\Lambda = \Lambda^{\mathrm{TP}}_{n}$, $\Lambda = \Lambda^{\mathrm{TD}}_{n}$ or $\Lambda = \Lambda^{\mathrm{HC}}_{n}$.
}
\prf{
Consider the setup of the previous proof.  We have
\bes{
\| f - f_{\Lambda} \|^2_{L^2(D,\nu)} \leq  \frac{1}{\min_{\bm{n} \notin \Lambda} \{ \chi^{\mix}_{\bm{n},m} \} } \| f \|^2_{\tilde{H}^m_{\mix}(D,\nu)},
}
where $\chi^{\mix}_{\bm{n},m}$ is as in \R{chidef}.  We now observe that
\bes{
\chi^{\mix}_{\bm{n},m} = \sum_{|\bm{j}|_{\infty} \leq m} \prod^{d}_{k=1} (n_k(n_k+1))^{j_k} \geq c_{m,d} | \bm{n} |^{2m}_{\mathrm{hc}} >  c_{m,d} n^{2m},\qquad \bm{n} \notin \Lambda,
}
where $\Lambda$ is any of the three index sets consider.  The result now follows immediately.
}

\section{Numerical results {and discussion}}\label{s:numerics}
We conclude this paper with several numerical experiments illustrating the theoretical results.  Unless otherwise stated we use Legendre polynomials on $D = (-1,1)^d$, hyperbolic cross index sets, samples drawn independently from the uniform measure on $\Omega$ and a threshold parameter $\epsilon = 10^{-8}$.

\subsection{Function regularity}
\begin{figure}[t]
\begin{center}
\begin{tabular}{cc}
 \includegraphics[width=8cm]{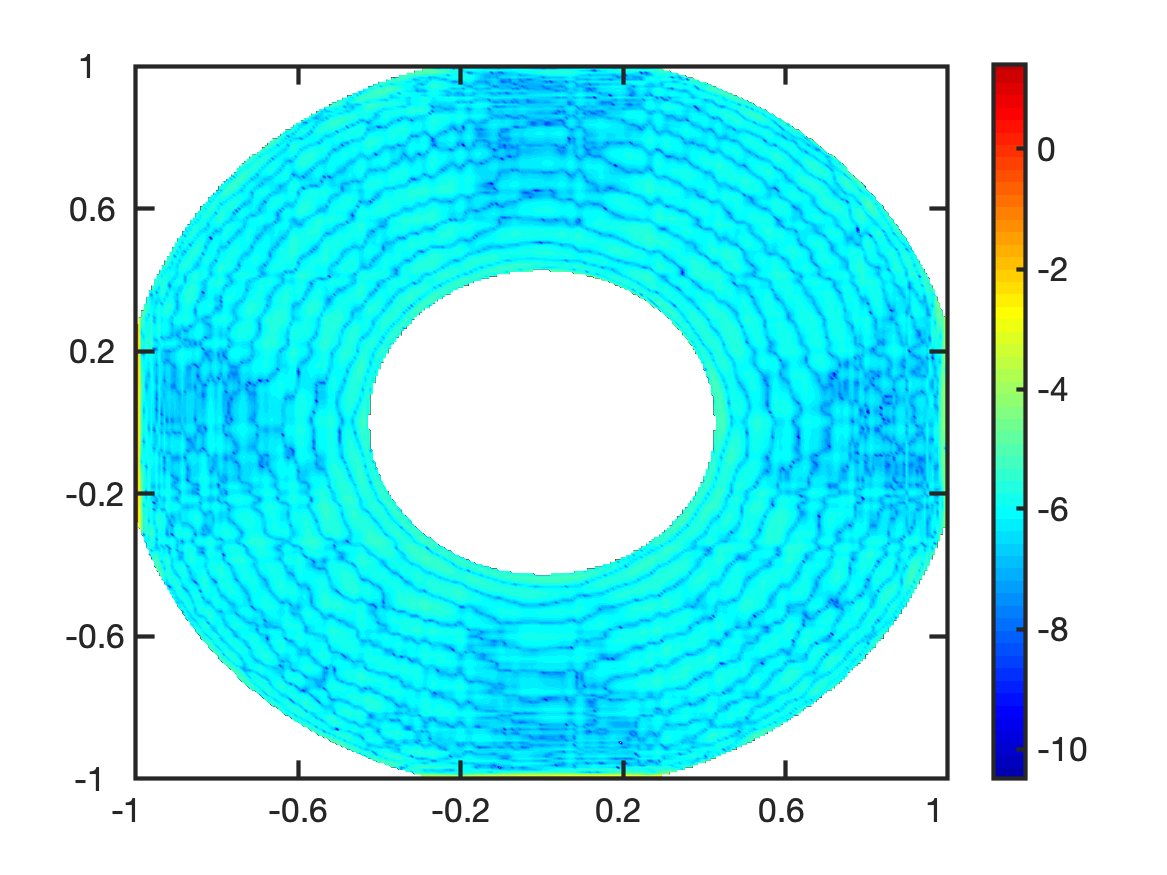} &  \includegraphics[width=8cm]{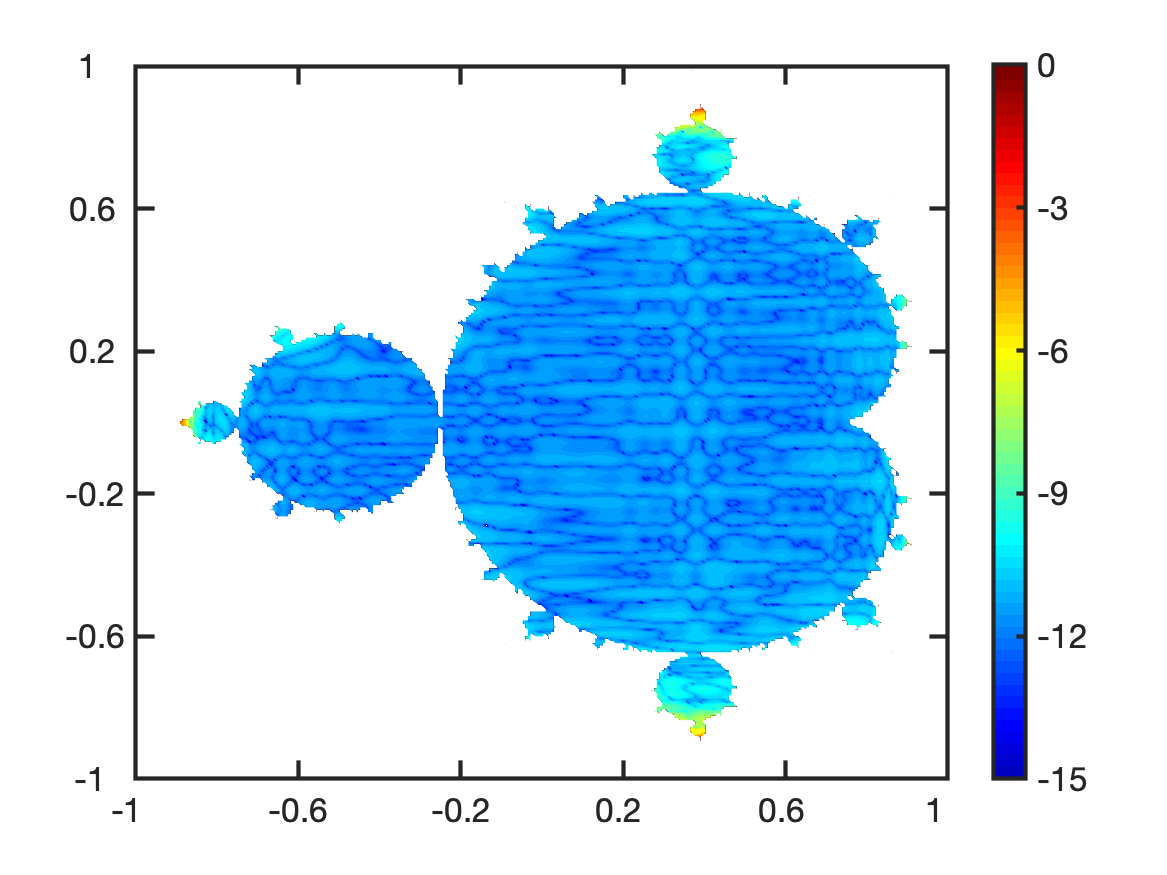}
 \\
 $f(y_1,y_2) = \log(8 (y^2_1+y^2_2)) - 2(y^2_1+y^2_2)$ & $f(y_1,y_2) = \cos(2 y_1) \sin(y_2)$
 \\
 $\Omega = \{ -1 \leq y_1 , y_2 \leq 1 : f(y_1,y_2) \geq 0 \}$ & $ \Omega  = \mbox{Mandelbrot set}$
 \\
 $\Lambda = \Lambda^{\mathrm{HC}}_{200}$, $N = 1102$, $M =  5510$ &  $\Lambda = \Lambda^{\mathrm{HC}}_{100}$, $N = 484$, $M = 2420$
 \end{tabular}
\end{center}
\vspace*{-5mm}
\caption{Pointwise error for polynomial frame approximations over two bivariate domains.}
\label{f:RegularityTest}
\end{figure}

We first consider the approximation of several bivariate functions.  The left panel of Fig.\ \ref{f:RegularityTest} shows the approximation of a smooth function on the domain $\Omega = \{\bm{y} : f(\bm{y}) \geq 0 \}$.  The function is singular on $D\backslash \Omega$.  Yet, as predicted by the results of \S \ref{s:approx_err}, this does not hamper its approximation on $\Omega$.  The right panel shows the approximation of a function defined on the Mandelbrot set.  This domain is not Lipschitz, but since the function has a smooth extension to the whole of $D$, an accurate approximation is obtained.  This also agrees with the results of \S \ref{s:approx_err}.  Note that in neither case does the domain need to be known in advance in order to compute the approximation.  It is defined implicitly by the data.

\subsection{Sample complexity}

\begin{figure}[t]
\begin{center}
\begin{tabular}{ccc}
 \includegraphics[width=5.0cm]{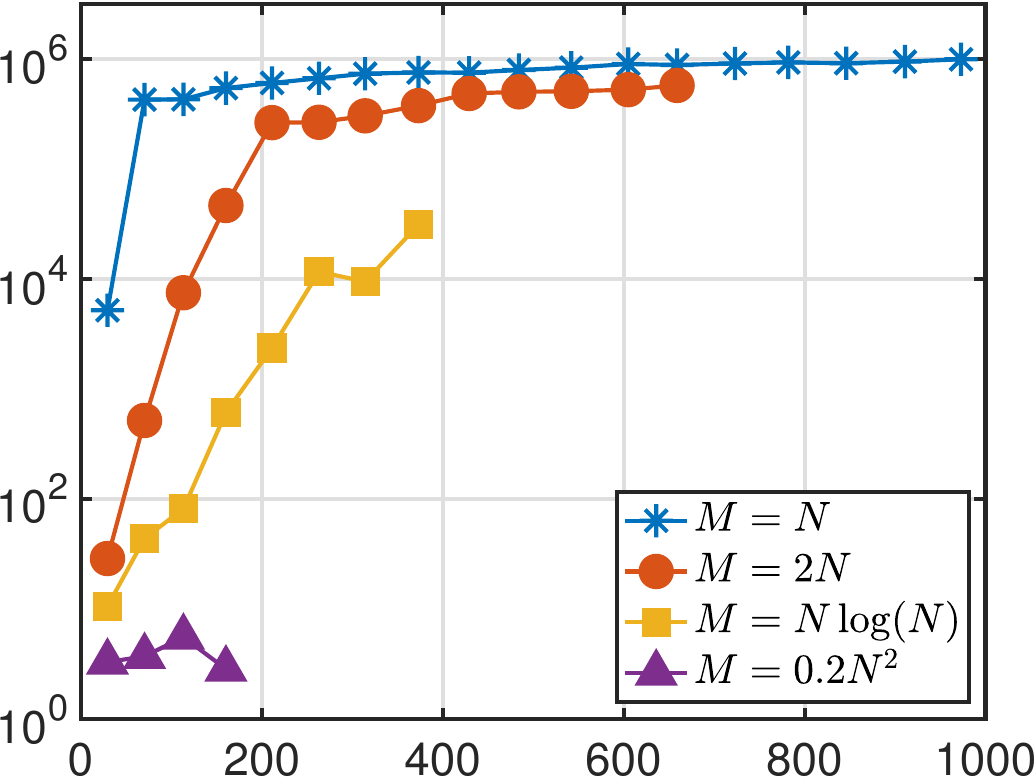} & 
  \includegraphics[width=5.0cm]{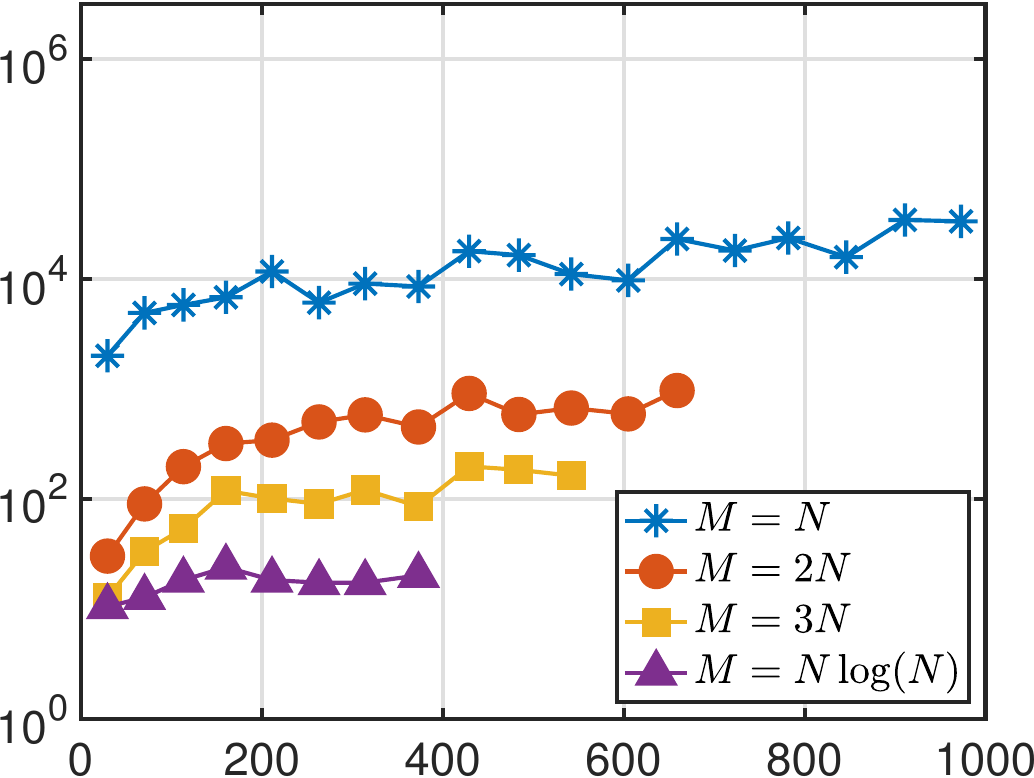} &
  \includegraphics[width=5.0cm]{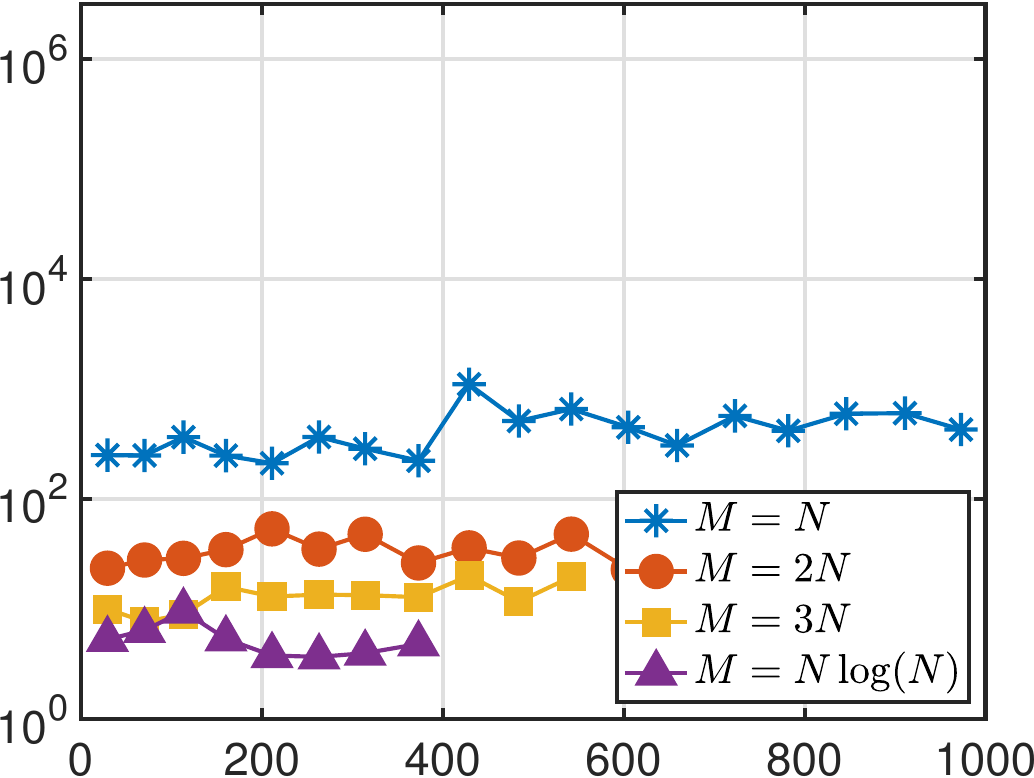} 
  \\  
 \includegraphics[width=5.0cm]{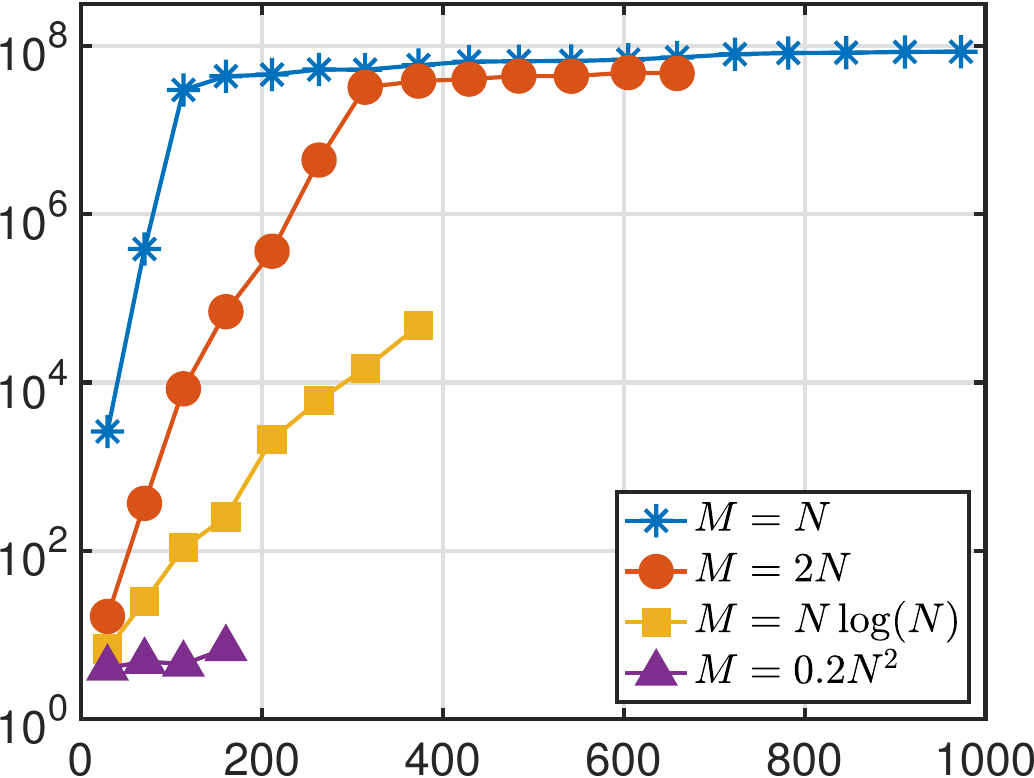} &
 \includegraphics[width=5.0cm]{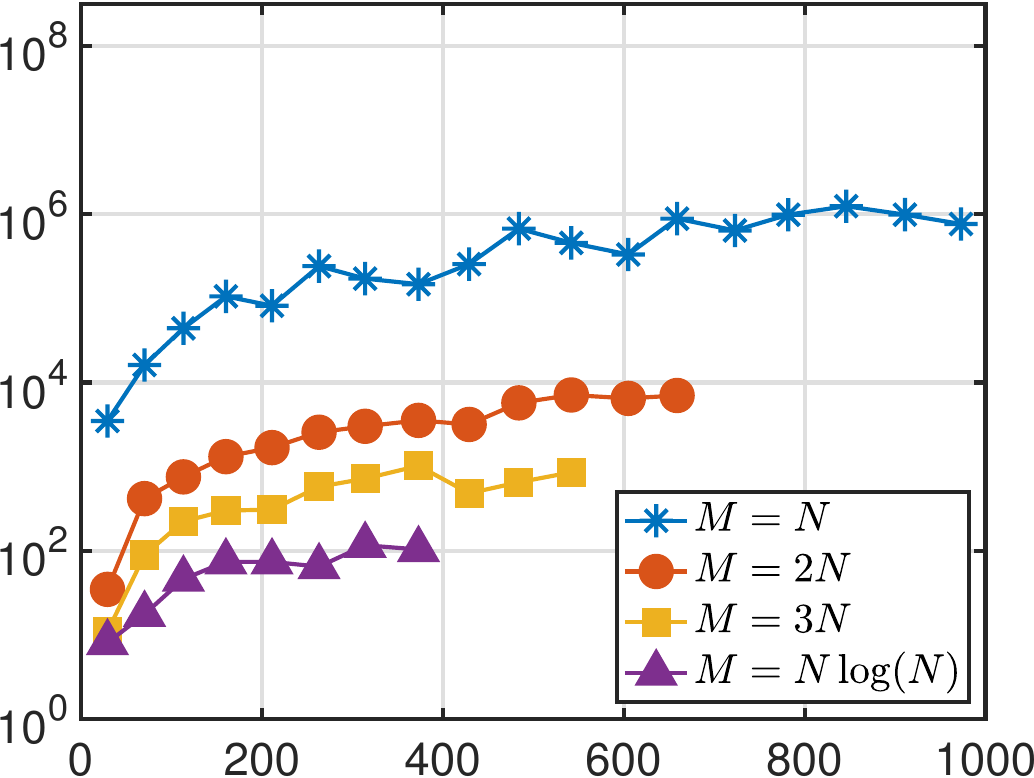} &   
   \includegraphics[width=5.0cm]{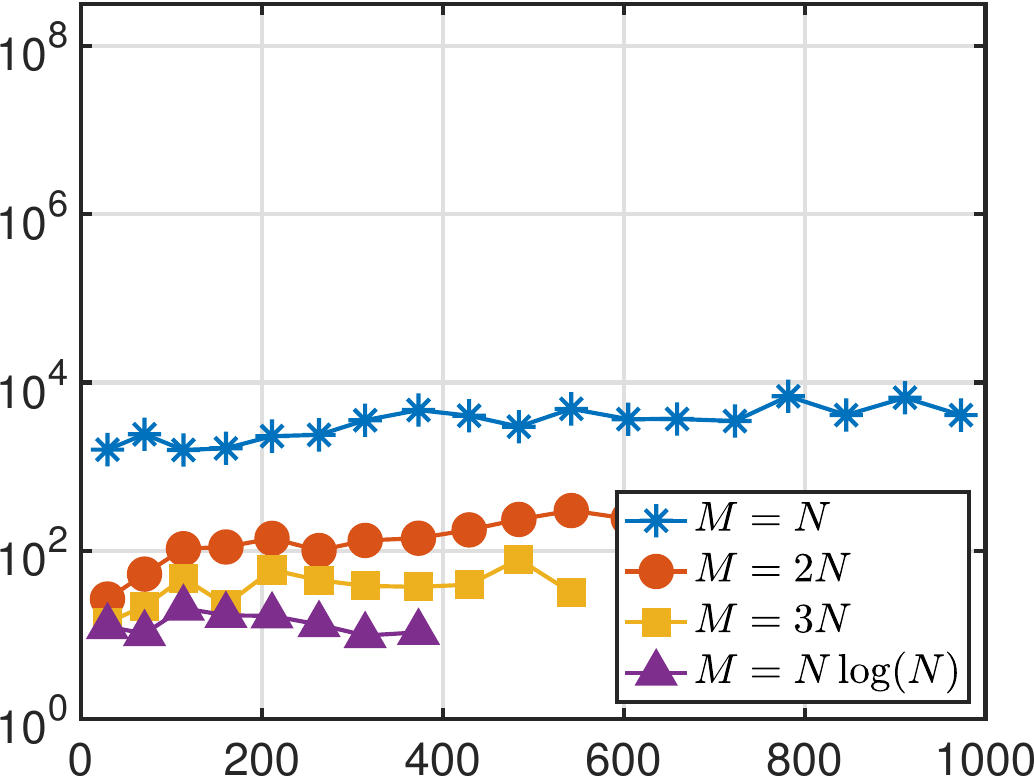}
   \\
       \includegraphics[width=5.0cm]{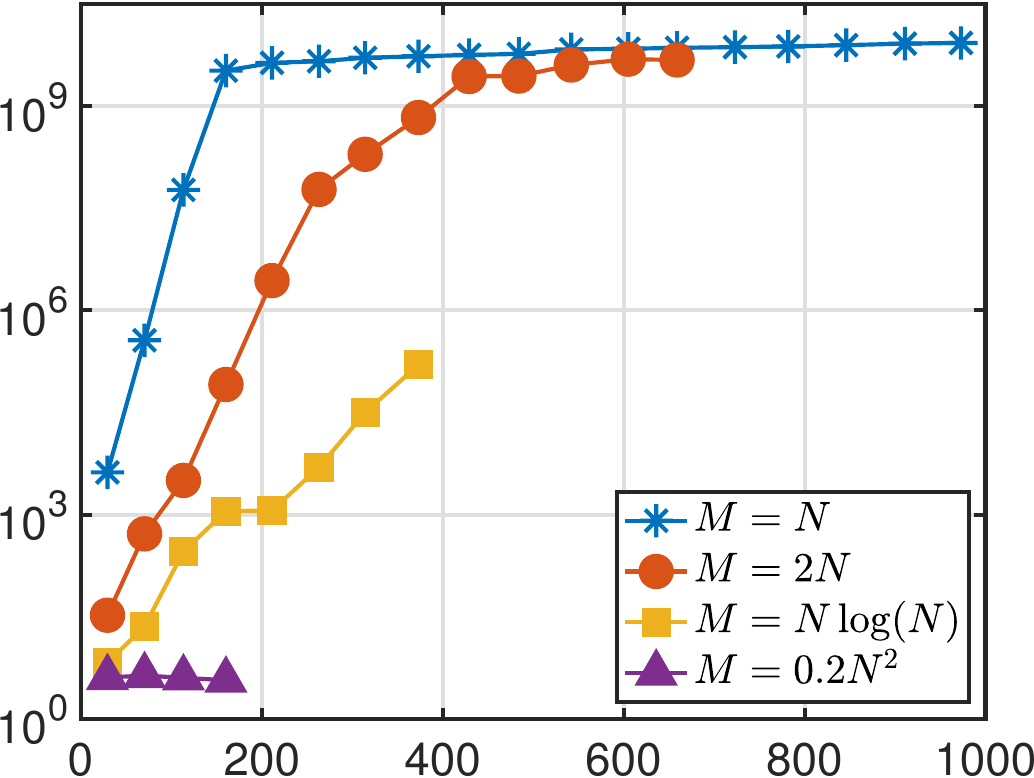} &
        \includegraphics[width=5.0cm]{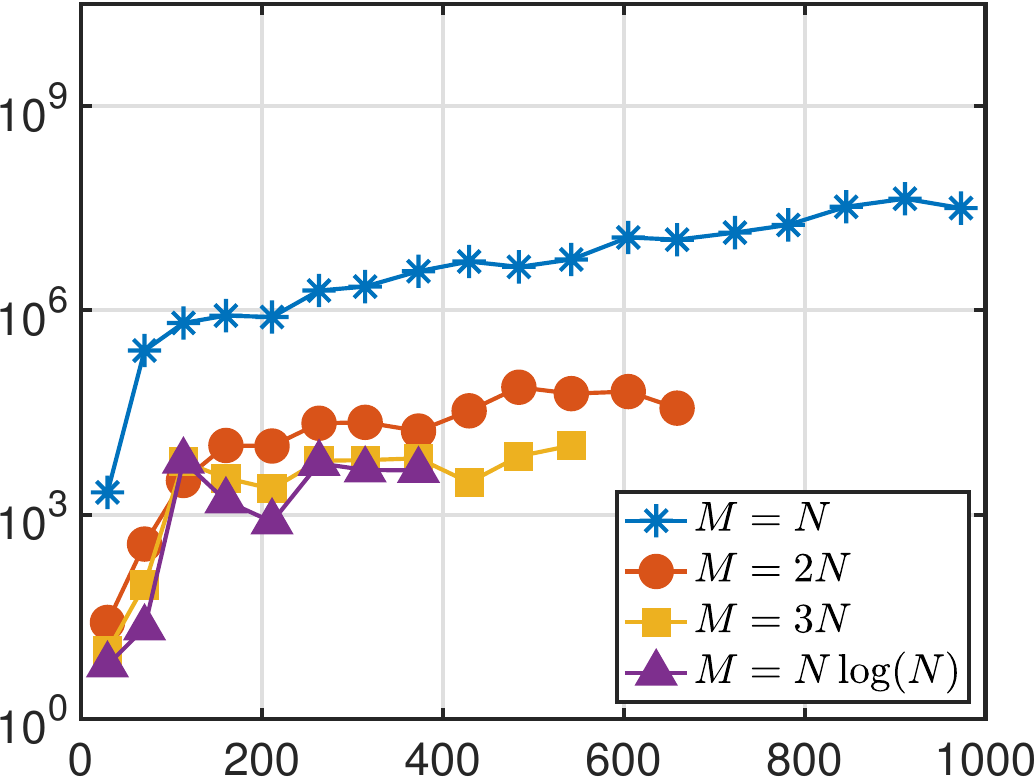} &
 \includegraphics[width=5.0cm]{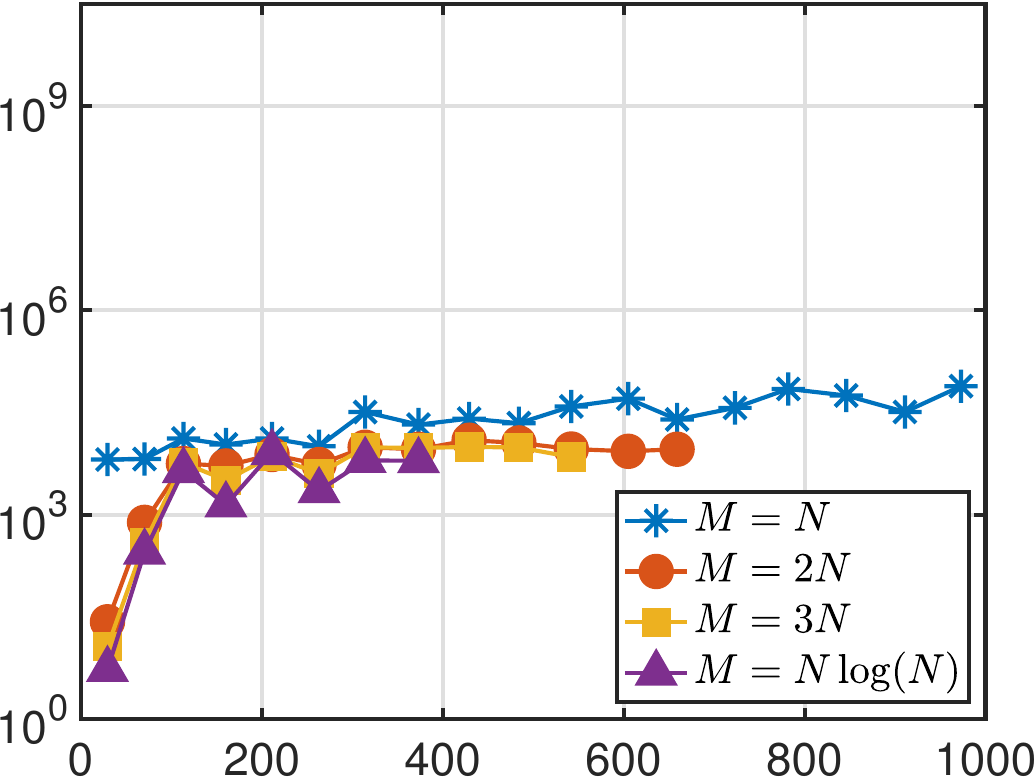} 
   \\
 (a): $r = 1$ & (b): $r = 3/4$ & (c): $r = 1/2$
  \end{tabular}
\end{center}
\vspace*{-5mm}
\caption{
{The constant $C_{\Upsilon,\Lambda,\epsilon}$ against $N$ for $\epsilon = 10^{-6}$, (top row) $\epsilon = 10^{-8}$ (middle row) and $\epsilon = 10^{-10}$ (bottom row).  The domain $\Omega$ is a circle of radius $r$ in $d=2$ dimensions.  Computations were averaged over $20$ trials with the median value taken.  The computation of $C_{\Upsilon,\Lambda,\epsilon}$ was done as in Remark \ref{r:constcompute}, using a precomputed grid of $K = 10000$ Monte--Carlo points in $\Omega$.
}
}
\label{f:RegularityTest2}
\end{figure}

{In Fig.\  \ref{f:RegularityTest2}} we examine the sample complexity of polynomial frame approximations for a two-dimensional circular domain.  This requires computing the constant $C_{\Upsilon,\Lambda,\epsilon}$, which is discussed in the remark below.  Fig.\ \ref{f:RegularityTest2}(a) suggests that quadratic oversampling is sufficient in this case, even though the domain is not of $\lambda$-rectangle type.  Moreover, linear or log-linear oversampling results in exponential increase of $C_{\Upsilon,\Lambda,\epsilon}$, up to roughly $1/\epsilon$ (recall that $C_{\Upsilon,\Lambda,\epsilon} \lesssim 1/\epsilon$; see \S \ref{ss:mainresaccstab}).  On the other hand, Figs.\ \ref{f:RegularityTest2}(b),(c) suggest that log-linear oversampling is sufficient whenever domain $\Omega$ does not touch the bounding cube $D$.  {Furthermore, the constant $C_{\Upsilon,\Lambda,\epsilon}$ gets smaller (for the same level of oversampling) as the radius $r$ of the domain decreases, or in other words, as the distance between the boundary of $\Omega$ and the boundary of $D$ grows.  These interesting observations, which are} at odds with the log-quadratic rates predicted in \S \ref{s:sampcomp}, {have been thoroughly documented in} the {one-dimensional case in the related setting when trigonometric polynomials are used instead of algebraic polynomials \cite{FEStability,FEParameterSelection}.}

While we currently have no proof, it is possible to give an intuitive explanation for this phenomenon.  The sample complexity relate{s} to the maximal growth of a polynomial (in an $L^2$-sense) on $\Omega$ when it is bounded at $M$ points in $\Omega$.  A polynomial that grows large in this sense must also be large on $D\backslash\Omega$, and therefore have large coefficients when represented in the Legendre basis.  Yet, when regularizing via the truncated SVD (which prohibits large coefficients), such polynomials are excluded from the resulting approximation space.  This also explains why {the constant $C_{\Upsilon,\Lambda,\epsilon}$ decreases as $r$ decreases}: for $r=1$ the boundaries of $\Omega$ and $D$ intersect, but {as $r$ decreases these boundaries are increasingly separated.} Formalizing this intuition into a proof is an open problem.

\rem{
\label{r:constcompute}
As shown in \cite{BADHFramesPart2}, the constants $C'_{\Upsilon,\Lambda,\epsilon}$ and $C''_{\Upsilon,\Lambda,\epsilon}$ can be expressed as
\be{
\label{C12matrix}
C'_{\Upsilon,\Lambda,\epsilon} = \sqrt{\lambda_{\max}\left ((\bm{B}')^* \bm{G} \bm{B}' \right ) },\qquad C''_{\Upsilon,\Lambda,\epsilon} = \epsilon^{-1} \sqrt{\lambda_{\max}\left ((\bm{B}'')^* \bm{G} \bm{B}''\right ) },
}
where $\bm{G} = \bm{G}_{\Lambda}$ is the Gram matrix of the truncated frame \R{Grammatrix}, $\bm{B}' = ( \bm{A}_{\Upsilon,\Lambda,\epsilon}  )^{\dag} = \bm{V} \left ( \bm{\Sigma}_{\epsilon} \right )^{\dag} \bm{U}^*$ and $\bm{B}'' = \bm{V} \bm{I}^{\perp}_{\epsilon}\bm{V}^{*}$.  Here $\bm{U} \bm{\Sigma} \bm{V}^*$ is the SVD of $A$, and $\bm{I}^{\perp}_{\epsilon}$ is the diagonal matrix with $\bm{n}^{\rth}$ entry $1$ if $\sigma_{\bm{n}} \leq \epsilon$ and zero otherwise.  Computing the Gram matrix $\bm{G}$ over an irregular domain is difficult, but it can be done approximately via Monte--Carlo integration.  Specifically, $\bm{G} \approx \bm{H}^* \bm{H}$ where 
\bes{
\bm{H} = \bm{H}_{K,\Lambda} = \left ( \frac{1}{\sqrt{K}} \phi_{\bm{n}}(\bm{z}_k) \right )_{k=1,\ldots,K,\bm{n} \in \Lambda} \in \bbC^{K \times N},
}
and $\bm{z}_1,\ldots,\bm{z}_K$ are drawn independently from $\mu$.  Replacing $\bm{G}$ by $\bm{H}^* \bm{H}$ in \R{C12matrix} and using standard properties of singular values leads to the simpler approximate expressions
\bes{
C'_{\Upsilon,\Lambda,\epsilon} \approx \nm{\bm{H} \bm{V} \left ( \bm{\Sigma}_{\epsilon} \right )^{\dag} }_{2},\qquad C''_{\Upsilon,\Lambda,\epsilon} \approx  \epsilon^{-1} \nm{\bm{H} \bm{V} \bm{I}^{\perp}_{\epsilon}}_{2},
}
where $\nm{\cdot}_{2}$ denotes the matrix $2$-norm.
}

\subsection{Higher dimensions}
\begin{figure}[t]
\begin{center}
\begin{tabular}{ccc}
  \includegraphics[width=5.0cm]{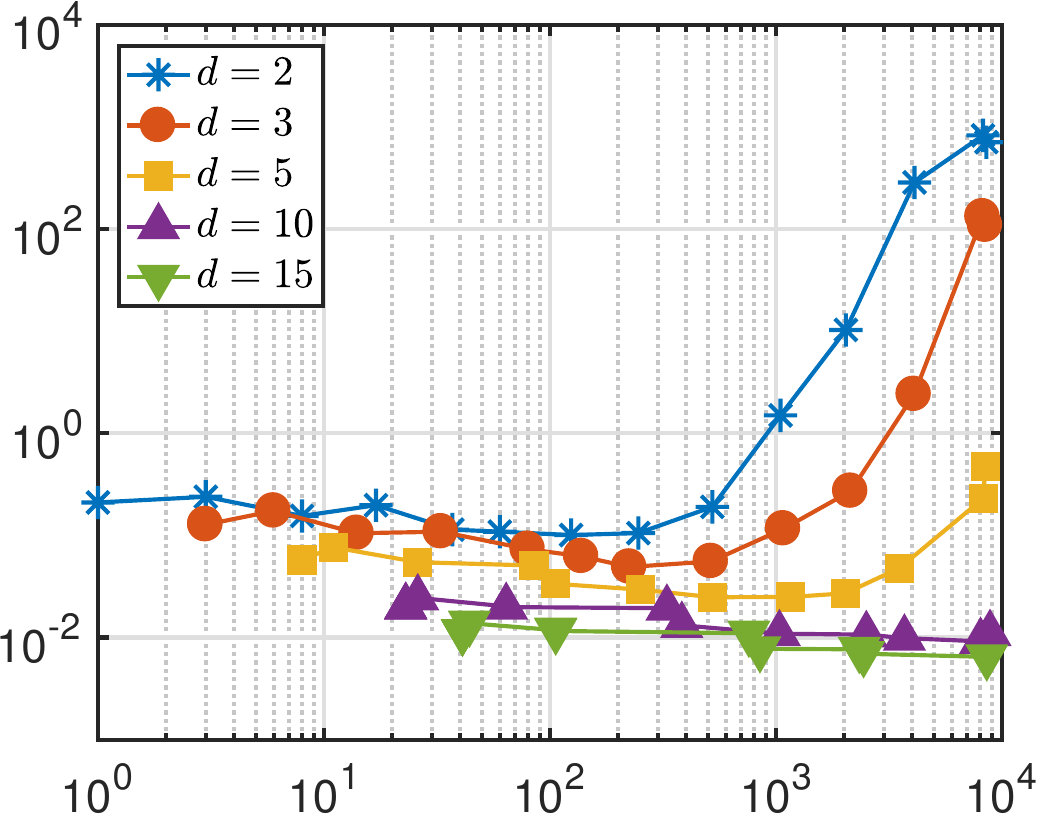} &
    \includegraphics[width=5.0cm]{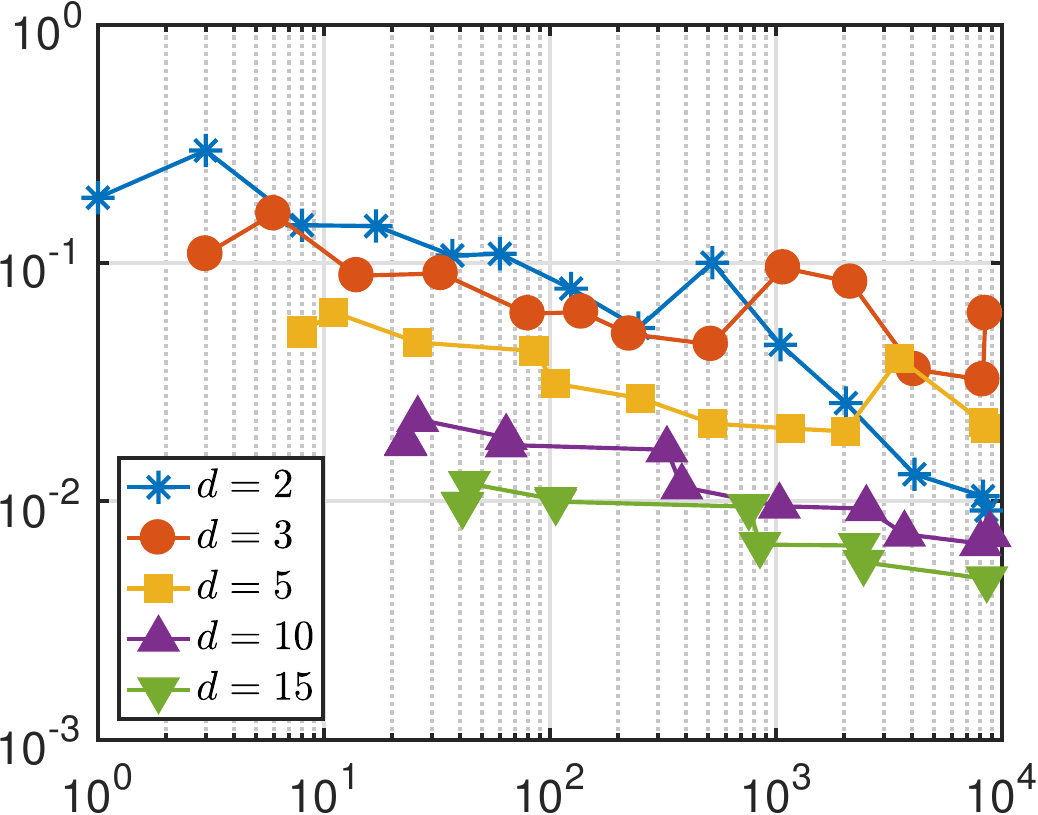} 
&
  \includegraphics[width=5.0cm]{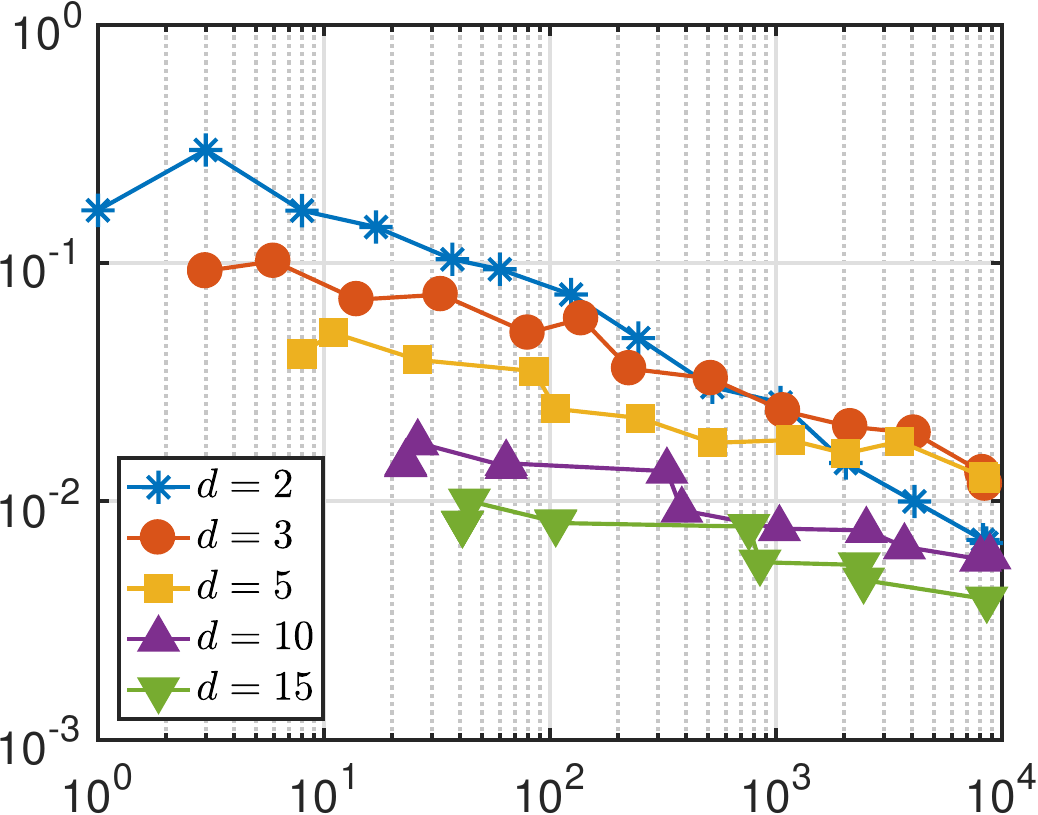}
   \\
 (a): $r = 1$ & (b): $r = 3/4$ & (c): $r = 1/2$
  \end{tabular}
\end{center}
\vspace*{-5mm}
\caption{
The median error over $20$ trials versus $M$ for approximating the function $f(\bm{y}) = 1/\sum^{d}_{i=1} \sqrt{|y_i|} $ on the annular domain $\Omega = \left \{ \bm{y} : r/4 \leq \nm{\bm{y}}_2 \leq r \right \}$.  For each $M$, the value of $N$ is chosen as the largest such that $N \log(N) \leq M$.
}
\label{f:DimensionTest}
\end{figure}

In Fig.\ \ref{f:DimensionTest} we consider the approximation error in various different dimensions.  This figure shows the approximation error versus $M$ for an annular region of several different radii.  In view of the previous discussion, log-linear oversampling was used throughout.  It is noticeable that when $r=1$, meaning that $\Omega$ touches the boundary of $D$, the approximation is ill-conditioned, and the error duly increases for large enough $M$.  As is to be expected, this increase is most severe in lower dimensions (since the cardinality of the polynomial space is largest in this setting).  Conversely, as soon as $\Omega$ is compactly contained in $D$, the approximation error decreases as $M$ increases.  Note that the function being approximated is smooth in $\Omega$ but singular at $\bm{y} = \bm{0} \in D \backslash \Omega$.  As predicted by the results of \S \ref{s:approx_err}, the approximation error decreases rapidly despite this singularity.

{
\subsection{Choice of $\epsilon$}\label{ss:epschoice}
}

{In this section, we discuss the influence of the regularization parameter on the approximation.  First, we note that the approximation is fairly robust to the choice of the parameter $\epsilon$.  In the noiseless setting, $\epsilon$ can be considered a \textit{target accuracy} for the method: namely, for sufficiently large $M$ and $N$, the approximation error will be on the order of $\epsilon$ (provided, of course, $\epsilon$ is larger than machine epsilon, since floating point error will always limit the best achievable accuracy in practice).  Indeed, under the mild conditions that the subspaces $\Lambda = \Lambda_N$ satisfy $\Lambda_1 \subset \Lambda_2 \subset \cdots $ and $\cup_{N} \Lambda_N = \bbN^d_0$ (which certainly holds for all choices considered in this paper), one has
\bes{
\limsup_{N \rightarrow \infty} \tilde{E}_{\Lambda_N,\epsilon}(f) \leq \epsilon \nm{f}_{L^2(\Omega,\mu)}.
}
This follows by choosing $p = g_{\Lambda}$ in \R{tildeELambdaDef}, where $g \in L^2(D,\nu)$ is the extension of $f$ by zero to $D$ and $g_{\Lambda}$ is its orthogonal projection (a similar conclusion holds for $E_{\Lambda,\epsilon}(f)$ under slightly stronger regularity assumptions, since this quantity involves an $L^{\infty}$-norm as opposed to the $L^2$-norm).}

{
This robustness is in stark contrast to the setting of ill-posed problems, where a careful choice of regularization parameter is usually crucial (see, for example, \cite{HansenIllPosed,NeumaierIllCond}). In such problems, one is typically interested in a specific solution $\bm{c}$ of the linear system, and the regularization parameter needs to be carefully chosen to strike a balance between the residual of the linear system and some desired property of $\bm{c}$ (e.g.\ smoothness).
On the other hand, our concern lies not with the vector $\bm{c}$, but rather with how well the function $f_{\Upsilon,\Lambda,\epsilon}$ approximates $f$, without preference for one set of coefficients over another, and this implies that success is measured largely by the size of the residual only. Furthermore, success is guaranteed for any $f$ by increasing $N$ due to the completeness of the polynomial frame.
}

{
The situation is slightly different if the function samples $f(\bm{y})$, $\bm{y} \in \Upsilon$, are corrupted by noise. In the setting of ll-posed problems, an optimal choice of the regularization parameter often involves the corner of the L-curve \cite{PCHansenLCurve}. 
Yet, the method of this paper remains robust in the noisy setting: the presence of noise merely implies that the limiting accuracy is determined by the maximum of $\epsilon$ and the noise level.
}

{
The above discussion assumes sufficient oversampling so that the constant $C_{\Upsilon,\Lambda,\epsilon}$ in Theorem \ref{t:acc_stab_cond} satisfies $C_{\Upsilon,\Lambda,\epsilon} \lesssim 1$.  The parameter $\epsilon$ also affects this constant.  Generally, $C_{\Upsilon,\Lambda,\epsilon}$ increases as $\epsilon$ decreases, reflecting the fact that as $\epsilon$ decreases more singular values are retained and the regularized approximation space becomes larger.  Hence, smaller $\epsilon$ generally means worse conditioning and accuracy of the approximation for fixed $M$.  Or equivalently, a larger $M$ is required to maintain the same level of conditioning and accuracy.  Note that this is not reflected in the sample complexity analysis conducted in \S \ref{s:sampcomp}, wherein the dependence on $\epsilon$ was ignored (recall Lemma \ref{l:CUB}).
Nevertheless, this has a practical impact.  If one requires only low accuracy (or if accuracy is limited by noise in data), then it is disadvantageous to take $\epsilon$ any smaller than needed.
}

{
These assertions are confirmed in Fig.\ \ref{f:RegularityTest2}.  For all choices of $M$ and $N$, a larger $\epsilon$ implies a smaller constant $C_{\Upsilon,\Lambda,\epsilon}$.  
These phenomena have also been investigated for the closely related \textit{Fourier extension} approximation in the one-dimensional setting, see \cite{FEParameterSelection}.  See also \S \ref{s:conclusion} for some further comments and formulation of open problems.
}

\subsection{Other bases}
\begin{figure}[t]
\begin{center}
\begin{tabular}{ccc}
  \includegraphics[width=5.0cm]{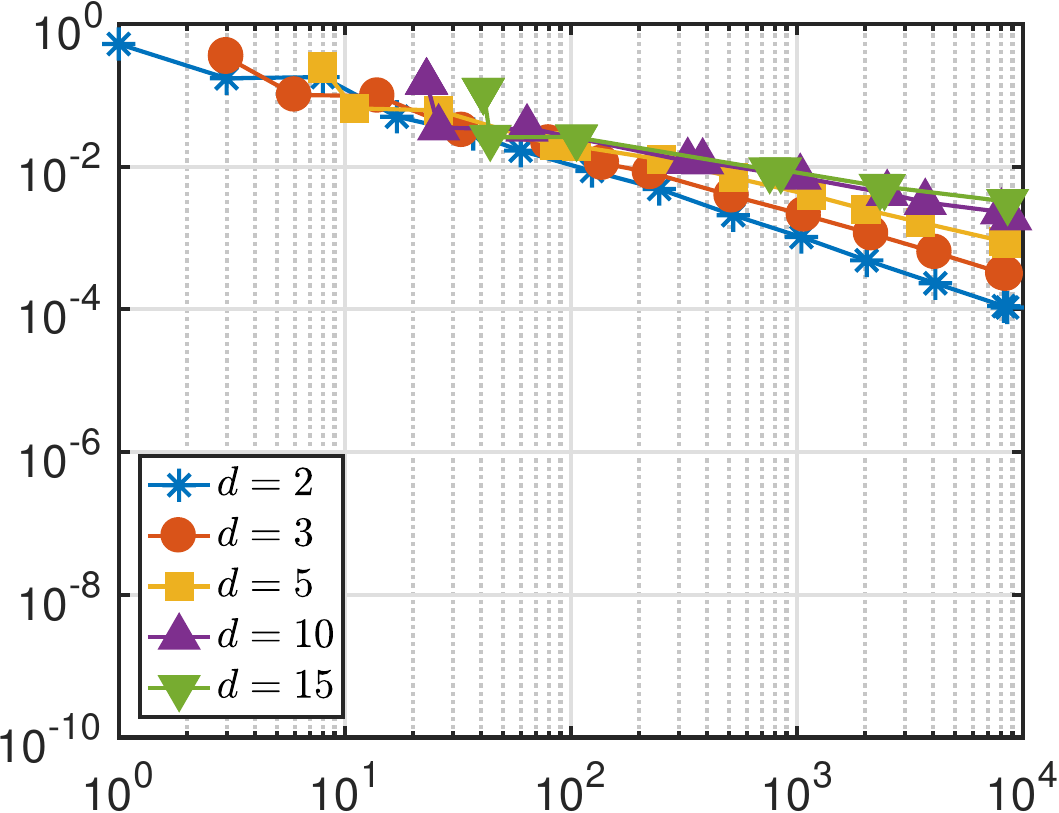} &
    \includegraphics[width=5.0cm]{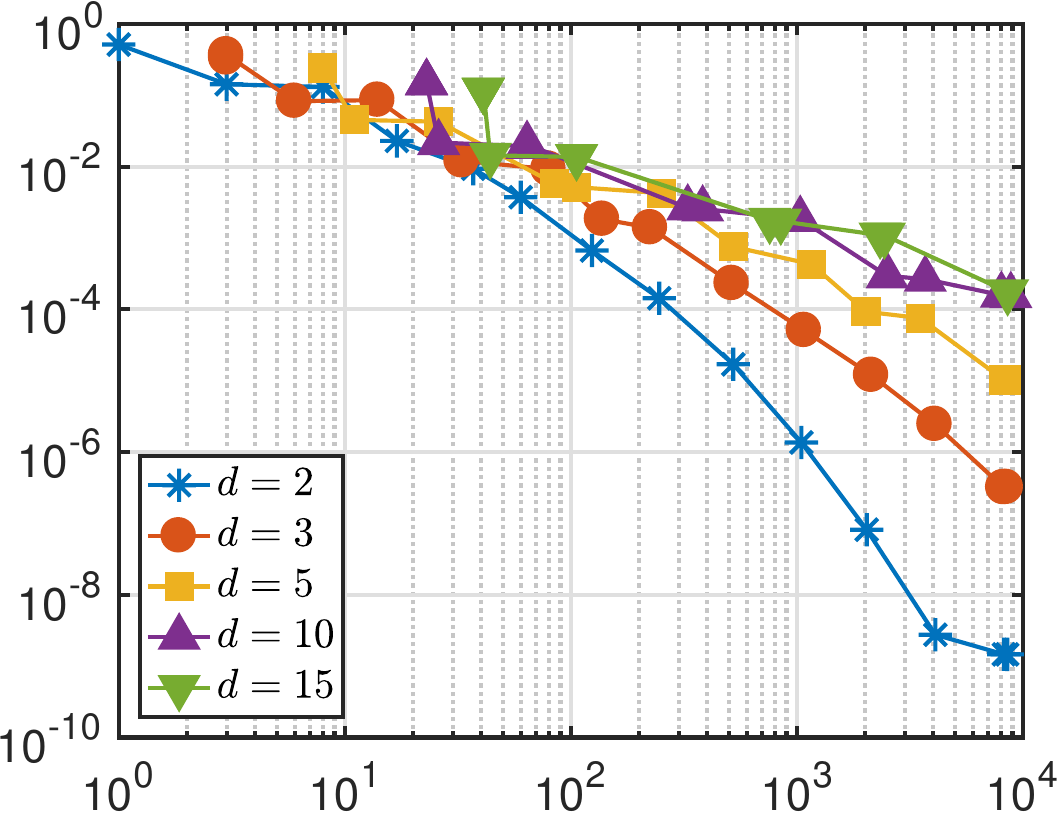}
  &
    \includegraphics[width=5.0cm]{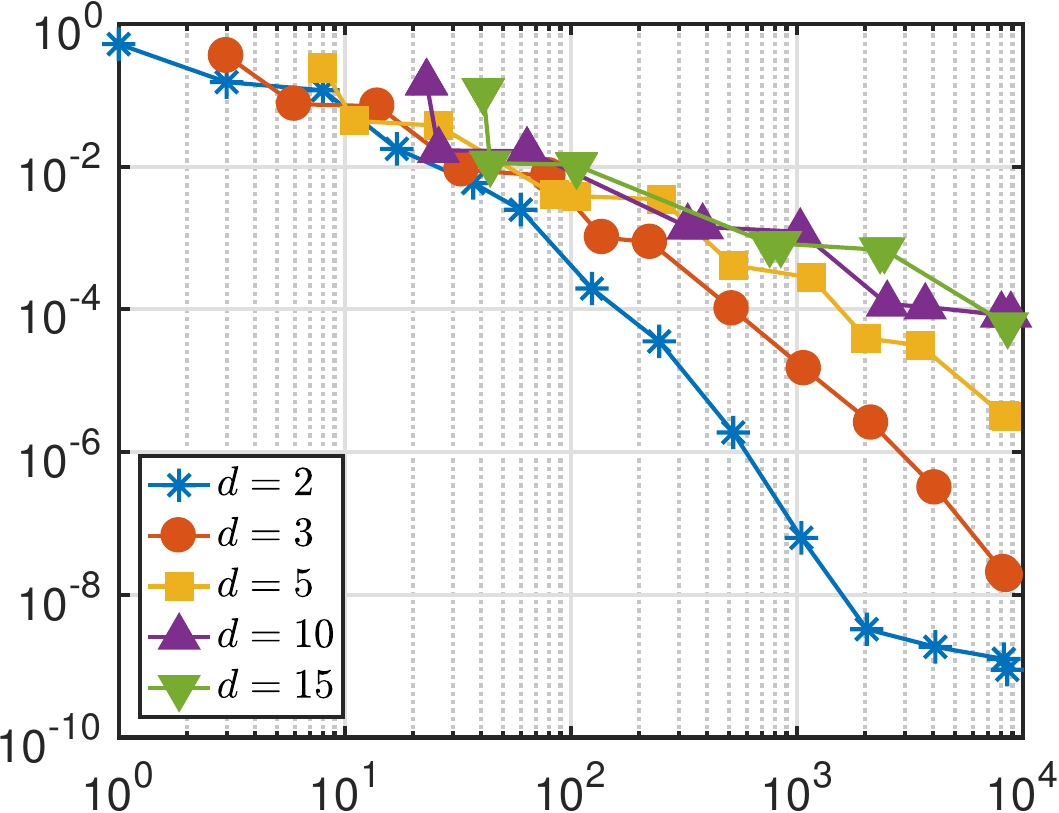} 
   \\
 (a): $T=1$ & (b): $T = 3/2$ & (c): $T = 2$
  \end{tabular}
\end{center}
\vspace*{-5mm}
\caption{
The median error over $20$ trials versus $M$ for approximating the function $f(\bm{y}) = \exp \left( - \sum^{d}_{i=1} y_i / d \right )$ on the corner domain $\Omega = \left \{ \bm{y} \in (-1,1)^d : y_1 + \ldots + y_d \leq 1 \right \}$ using the tensor cosine basis on $[-T,T]^d$.  For each $M$, the value of $N$ is chosen as the largest such that $N \log(N) \leq M$.
}
\label{f:DimensionTest1}
\end{figure}

\begin{figure}[t]
\begin{center}
\begin{tabular}{ccc}
  \includegraphics[width=5.0cm]{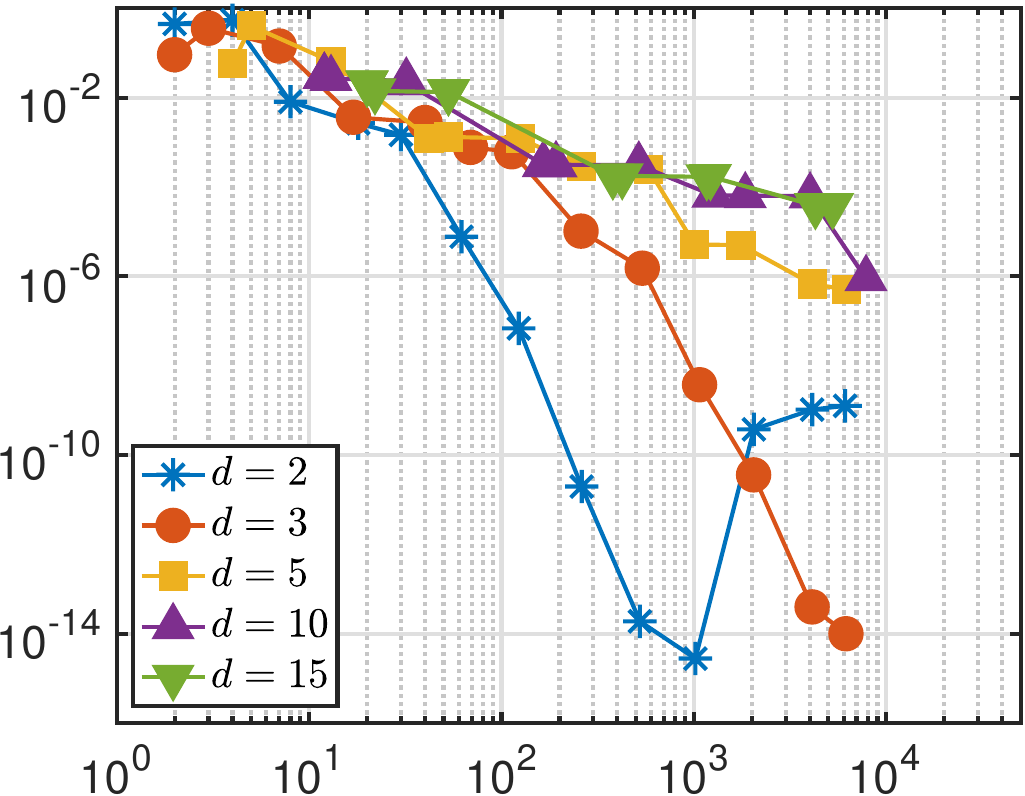} &
    \includegraphics[width=5.0cm]{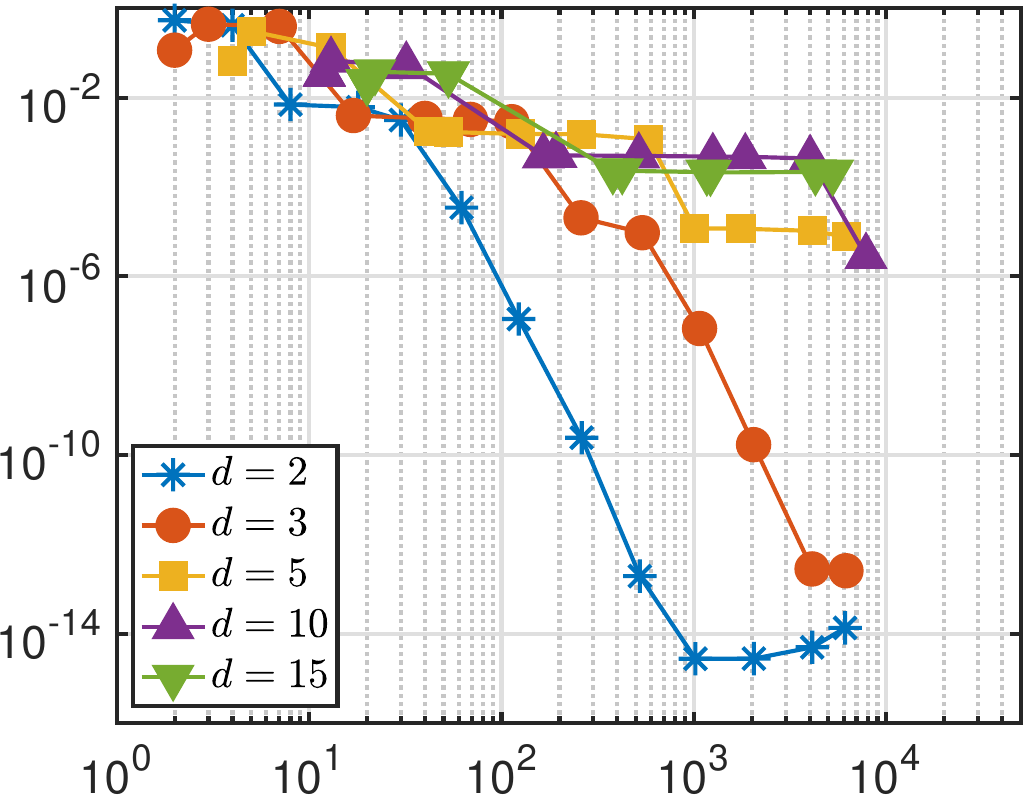} &
   \includegraphics[width=5.0cm]{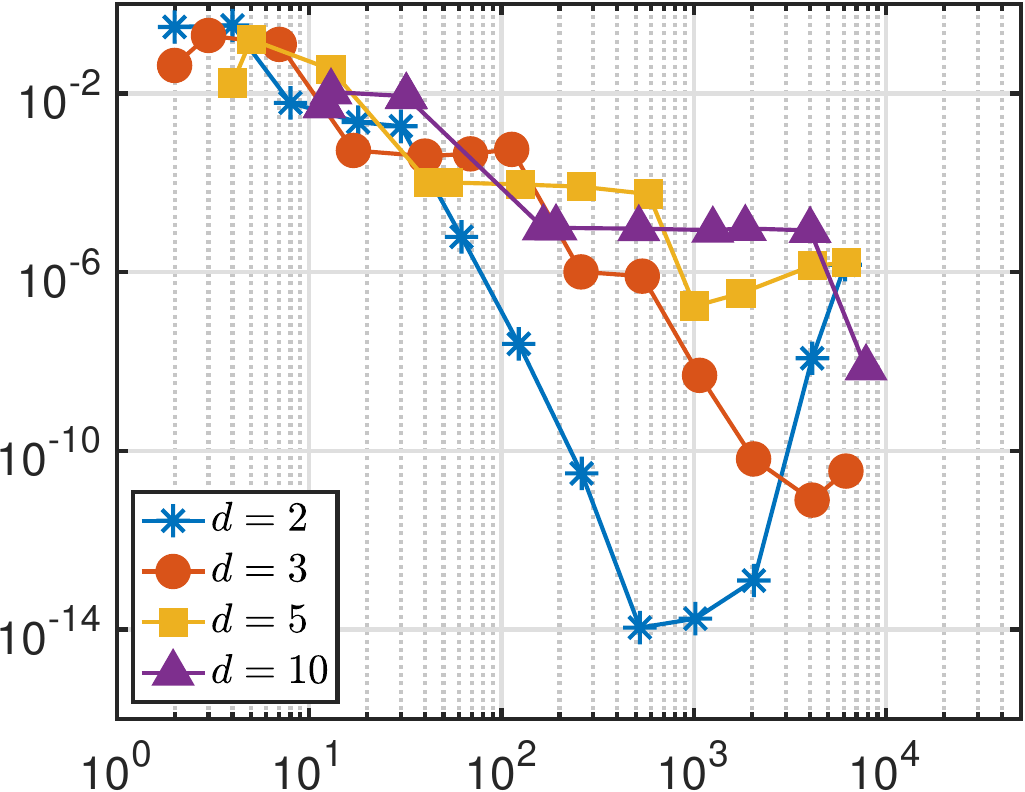}
  \end{tabular}
\end{center}
\vspace*{-5mm}
\caption{
The median error over $20$ trials versus $M$ for approximating the {function $f(\bm{y}) = \cos (\sum^{d}_{i=1} y_i / d )$ on the domains $\Omega = \{ \bm{y} \in (-1,1)^d : y_1+\ldots+y_d \leq 1 \}$ (left), $\Omega = \{ \bm{y} \in (-1,1)^d : \nm{\bm{y}}_2 \geq 1/2 \}$ (middle) and $\Omega = \{ \bm{y} \in (-1,1)^d : \nm{\bm{y}}_2 \leq 1 \}$ (right) using Chebyshev polynomials on $D = (-1,1)^d$ and samples drawn randomly according to the Chebyshev measure restricted to $\Omega$.  For each $M$, the value of $N$ is chosen as the largest such that $\tfrac12 N \log(N) \leq M$.}
}
\label{f:DimensionTest2}
\end{figure}

Finally, in Figs.\ \ref{f:DimensionTest1} \& \ref{f:DimensionTest2} we use different orthogonal bases on the extended domain.  First, in Fig.\ \ref{f:DimensionTest1} we consider the tensor cosine basis defined on $D = (-T,T)^d$, where $T \geq 1$ is a parameter.   The basis elements in this case are tensor-products on the univariate functions
\bes{
\phi_n(y) = \cos(n \pi ( y+T)/(2T)),\qquad n = 0,1,2,\ldots.
}
When $T = 1$, the domain $\Omega$ is not compactly contained in $D$ and the approximation error decreases slowly, at rate of $N^{-1}$.  This stems from the fact that cosine expansions, much like Fourier expansions, only converge rapidly for smooth functions that satisfy additional boundary conditions \cite{BAthesis}.  When there is no gap between $\Omega$ and the boundary of $D$, there are no smooth extensions of $f$ satisfying these boundary conditions.  Conversely, once $T > 1$ and $\Omega$ is compactly-contained in $D$, such extensions exist, and we witness correspondingly faster convergence.  Error estimates similar to those proved in \S \ref{s:approx_err} can also be established for these approximations.  See \cite{BA2} for further details.

{
In Fig.\ \ref{f:DimensionTest2} we consider Chebyshev polynomials on $D = (-1,1)^d$ and random sampling according to the tensor Chebyshev density restricted to $\Omega$ with log-linear oversampling. Notice that the $d = 2$ approximation exhibits instability.  We conjecture that this is related to distribution of the samples points.
Points drawn on a cube according to the Chebyshev density cluster quadratically near the boundary of the cube, a property which generally permits a lower sample complexity.  However, points drawn according to the same density when restricted to a subdomain $\Omega$ do not necessarily cluster in this way over the boundary of $\Omega$.  Unless $\Omega$ is compactly contained in $D$, it appears the severity of the instability is related to the amount of boundary $\Omega$ and $D$ share. 
}

\section{Conclusions and challenges}\label{s:conclusion}

In this work, we have introduced and analyzed a method, known as polynomial frame approximation, for approximating multivariate functions on irregular domains.  Amongst the various results proved, we have shown that for functions of mixed Sobolev regularity the regularized least-squares polynomial frame approximation is well-conditioned and converges algebraically-fast down to a given threshold parameter $\epsilon$.  Moreover, for a large class of domains, the sample complexity is provably quadratic in the dimension of the approximation space, up to a log factor.

This paper marks only a first foray into the broader topic of multivariate polynomial approximation on irregular domains.  Consequently, there are a number of interesting challenges for future research.  We conclude by highlighting three directions for further work:

\pbk
\textit{1.\ Sample complexity estimates.} When sampling from the uniform measure, we have shown log-quadratic sample complexity for $\lambda$-rectangle domains, {with the factor $\lambda^{-1}$ appearing in the sample complexity bound.  It is unknown whether or not this factor is sharp.  Moreover,} as mentioned, many domains do not have this property.  We conjecture that the same sample complexity holds for a much more general class of domains which includes spheres and simplices (two notable domains which do not have the $\lambda$-rectangle property) and which is potentially also invariant under rotations (rotations generally destroy the $\lambda$-rectangle property).  This remains an open problem.  Moreover, as discussed in \S \ref{s:numerics}, log-linear sample complexity appears to be sufficient whenever $\Omega$ is compactly contained in $D$.  While there is intuition behind this observation, we currently have no proof.

\pbk
\textit{2.\ Choice of $\epsilon$.} As discussed in \S \ref{ss:epschoice}, the method is robust to the choice of $\epsilon$.  Yet this parameter does affect the sample complexity.  Understanding the intricate relationship between the sample complexity, the domains $\Omega$ and $D$, the subspace $P_{\Lambda}$ and the parameter $\epsilon$ is very much an open problem.
As also noted in \S \ref{ss:epschoice}, in this paper we consider a fixed $\epsilon$ chosen according to some desired target accuracy (in our experiments we have simply taken $\epsilon = 10^{-8}$).  The possibility of \textit{adaptive} strategies, choosing $\epsilon$ depending on $\Lambda$ and $f$, is a topic for future work.

\pbk
\textit{3.\ Optimal sampling.}  {Recent work} has identified densities for random sampling which achieve near-optimal log-linear sample complexities for least-squares approximations \cite{MiglioratiCohenOptimal}.  {While these densities can be defined over irregular domains, it becomes challenging to sample efficiently from them in the case where the domain is not of tensor-product type.  One solution to this problem is to employ discrete measures, supported over a fine grid that suitably fills $\Omega$.  This strategy, which uses ideas of \cite{MiglioratiCohenOptimal}, has been recently developed in \cite{BAJMCoptimal,MiglioratiIrregular}.  Yet this procedure requires the domain $\Omega$ to be known in advance, and requires a fine grid to first be generated.  This may not be possible in all applications, especially in higher dimensions.  For instance, the case $\Omega = \{ \bm{y} : f(\bm{y}) \geq 0 \}$, which arises in practical surrogate model construction problems (see \S \ref{s:introduction}), presents clear difficulties.  Developing efficient sampling procedures for such problems remains a topic for future investigation.}

\pbk
\textit{4.\ Compressed sensing-based polynomial approximations.} Polynomial-based compressed sensing approaches have recently proved effective for high-dimensional approximation in regular domains (see \cite{AdcockCSFunInterp,BASBCWMatheon,ChkifaDownwardsCS,NarayanZhouCCP,PengHamptonDoostantweighted,YanGuoXui_l1UQ} and references therein).  A problem for future work is to extend these approaches to irregular domains.  Note that since polynomial frames are redundant, the usual compressed sensing theory for orthogonal bases does not apply.

\section*{Acknowledgements}
A preliminary version of this work was presented at the Oberwolfach conference on ``Multiscale and High-Dimensional Problems''.  The authors would like to thank the organizers and participants for the useful discussions and feedback received during the conference.  They would also like to thank Claudio Canuto, Juan M.\ Cardenas, Albert Cohen, Paul Constantine, Ingrid Daubechies, Bert Debusschere, Alireza Doostan, Wolfgang Hackbusch, Sebastian Moraga, Vladimir Temlyakov and Tino Ullrich.  The first author is supported by NSERC grant 611675, as well as an Alfred P.\ Sloan Research Fellowship. The second author is supported by FWO-Flanders projects G.0641.11 and G.A004.14, as well as by KU Leuven project C14/15/055.

\appendix

\section{Background on Legendre polynomials}\label{s:Legbackground}
This section contains some ancillary results on Legendre polynomials used earlier in the paper.  Let $\{ \psi_{n} \}^{\infty}_{n=0}$ be the orthonormal Legendre polynomial basis with respect to the uniform measure on $(-1,1)$.  This is defined by
\be{
\label{1DLegDef}
\psi_{n}(y) = \sqrt{2n+1} P_{n}(y),
}
where $P_n$ is the classical Legendre polynomial with normalization $P_n(1) = 1$.

\subsection{One dimensional Legendre--Sobolev spaces}
Recall that $\psi_{n}(y)$ are the eigenfunctions of the Sturm--Liouville operator $\cL$, defined by
\bes{
\cL f(y) = \left ( (1-y^2) f'(y) \right )'.
}
Specifically, $\cL\psi_{n}(y) = n(n+1) \psi_{n}(y)$.
The operator $\cL$ is compact, self-adjoint and nonnegative definite.  Note that
\bes{
\ip{\cL f}{g}_{L^2(D,\nu)} = \ip{f'}{g'}_{L^2(D,\rho)} = \ip{f}{\cL g}_{L^2(D,\nu)},
}
where $D = (-1,1)$, $\nu$ is the uniform measure on $(-1,1)$ and $\D \rho(y) = \frac{1-y^2}{2} \D y$.
The operator $\cL$ has a well-defined square root, which we write as $\cS = \cL^{1/2}$.  Note that
\bes{
\nm{\cS f }^2_{L^2(D,\nu)} = \| f' \|_{L^2(D,\rho)} = \ip{\cL f}{f}_{L^2(D,\nu)}.
}
With this in hand, for $j \in \bbN$ let $\cS^{j} = \cS \circ \cS \circ  \cdots \circ  \cS $ be the $j$-fold composition of $\cS$ and define the Legendre--Sobolev space
\bes{
\tilde{H}^{m}(D,\nu) = \left \{ f \in L^2(D,\nu) : \cS^{j} f \in L^2(D,\nu),\ j=0,\ldots,m \right \},
} 
with inner product and norm
\bes{
\ip{f}{g}_{\tilde{H}^m(D,\nu)} = \sum^{m}_{j=0} \ip{\cS^{j} f}{\cS^{j} g}_{L^2(D,\nu)},
\qquad
\nm{f}_{\tilde{H}^m(D,\nu)} = \sqrt{\sum^{m}_{j=0} \nm{\cS^{j} f }^2_{L^2(D,\nu)} }.
}
The set $\{ \psi_{n} \}_{n \in \bbN_0}$ is an orthogonal basis for $\tilde{H}^m(D,\nu)$, and one has the expression
\bes{
\nm{f}_{\tilde{H}^m(D,\nu)} = \sqrt{ \sum^{\infty}_{n=0} \chi_{n,m} \left | \ip{f}{\psi_n}_{L^2(D,\nu)} \right |^2 },\qquad \chi_{n,m} = \sum^{m}_{j=0} \left ( n(n+1) \right )^j.
}
Here we use the convention $0^0 = 1$.

\subsection{Multidimensional Legendre--Sobolev spaces}
Let $D = (-1,1)^d$ be the unit cube and define the tensor Legendre polynomial basis $\{ \psi_{\bm{n}} \}_{\bm{n} \in \bbN^d_0}$ as
\be{
\label{dDLegDef}
\psi_{\bm{n}}(\bm{y}) = \prod^{d}_{k=1} \psi_{n_k}(y_k),\qquad \bm{n} = (n_1,\ldots,n_d) \in \bbN^d_0,\ \bm{y} = (y_1,\ldots,y_d) \in D.
}
For $k=1,\ldots,d$, let $\cL_{k}$ be the compact, self-adjoint nonnegative definite operator
\bes{
\cL_k f(\bm{y}) = \frac{\partial}{\partial y_k} \left ( (1 - y^2_k) \frac{\partial f}{\partial y_k} \right ),
}
with corresponding square-root $\cS_k = \cL^{1/2}_k$ and powers $\cS^{j}_{k} = \cS_k \circ \cdots \circ \cS_k$.  Now let $\bm{j} = (j_1,\ldots,j_d) \in \bbN^d_0$ be a multi-index.  We define the operator
\bes{
\cS^{\bm{j}} = \cS^{j_1}_{1} \circ \cdots \circ \cS^{j_d}_{d}.
}
With this in hand, we now define the $d$-dimensional Legendre--Sobolev spaces
\bes{
\tilde{H}^{m}(D,\nu) = \left \{ f \in L^2(D,\nu) : \cS^{\bm{j}} f \in L^2(D,\nu),\ | \bm{j} |_1 \leq m \right \},
}
with inner product and norm
\bes{
\ip{f}{g}_{\tilde{H}^m(D,\nu)} = \sum_{|\bm{j} |_1 \leq m} \ip{\cS^{\bm{j}} f}{\cS^{\bm{j}} g}_{L^2(D,\nu)},
\qquad
\nm{f}_{\tilde{H}^m(D,\nu)} = \sqrt{\sum_{|\bm{j} |_1 \leq m} \nm{\cS^{\bm{j}} f }^2_{L^2(D,\nu)}}.
}
We also define the mixed $d$-dimensional Legendre--Sobolev spaces as
\be{
\label{tildeHmix}
\tilde{H}^{m}_{\mix}(D,\nu) = \left \{ f \in L^2(D,\nu) : \cS^{\bm{j}} f \in L^2(D,\nu),\ | \bm{j} |_{\infty} \leq m \right \},
}
with inner product and norm
\bes{
\ip{f}{g}_{\tilde{H}^m_{\mix}(D,\nu)} = \sum_{|\bm{j} |_{\infty} \leq m} \ip{\cS^{\bm{j}} f}{\cS^{\bm{j}} g}_{L^2(D,\nu)},
\qquad
\nm{f}_{\tilde{H}^m_{\mix}(D,\nu)} = \sqrt{\sum_{|\bm{j} |_{\infty} \leq m} \nm{\cS^{\bm{j}} f }^2_{L^2(D,\nu)}}.
}
Both these norms can be characterized in terms of Legendre polynomial coefficients.  Specifically,
\be{
\label{Htilde_char}
\begin{split}
\nm{f}_{\tilde{H}^m(D,\nu)} = \sqrt{ \sum_{\bm{n} \in \bbN^d_0} \chi_{\bm{n},m} \left | \ip{f}{\psi_{\bm{n}}}_{L^2(D,\nu)} \right |^2 },
\\
\nm{f}_{\tilde{H}^m_{\mix}(D,\nu)} = \sqrt{ \sum_{\bm{n} \in \bbN^d_0} \chi^{\mix}_{\bm{n},m} \left | \ip{f}{\psi_{\bm{n}}}_{L^2(D,\nu)} \right |^2 },
\end{split}
}
where
\be{
\label{chidef}
\chi_{\bm{n},m} = \sum_{|\bm{j}|_1 \leq m} \prod^{d}_{k=1} (n_k(n_k+1))^{j_k},\qquad \chi^{\mix}_{\bm{n},m} = \sum_{|\bm{j}|_\infty \leq m} \prod^{d}_{k=1} (n_k(n_k+1))^{j_k}.
}
Finally, we note that one has the continuous embeddings $H^{m}(D,\nu) \hookrightarrow \tilde{H}^m(D,\nu)$ and $H^{m}_{\mix}(D,\nu) \hookrightarrow \tilde{H}^m_{\mix}(D,\nu)$.

\bibliographystyle{abbrv}
\small
\bibliography{PEHighDimRefs}

\begin{thebibliography}{10}

\bibitem{BAthesis}
B.~Adcock.
\newblock {\em Modified Fourier expansions: theory, construction and
  applications}.
\newblock PhD thesis, University of Cambridge, 2010.

\bibitem{BA2}
B.~Adcock.
\newblock Multivariate modified {F}ourier series and application to boundary
  value problems.
\newblock {\em Numer. Math.}, 115(4):511--552, 2010.

\bibitem{AdcockCSFunInterp}
B.~Adcock.
\newblock Infinite-dimensional compressed sensing and function interpolation.
\newblock {\em Found. Comput. Math.}, 18(3):661--701, 2018.

\bibitem{BASBCWMatheon}
B.~Adcock, S.~Brugiapaglia, and C.~G. Webster.
\newblock Compressed sensing approaches for polynomial approximation of
  high-dimensional functions.
\newblock In {\em Compressed Sensing and Its Applications}. Birkh\"auser, 2017.

\bibitem{BAJMCoptimal}
B.~Adcock and J.~M. Cardenas.
\newblock Near-optimal sampling strategies for multivariate function
  approximation on general domains.
\newblock {\em SIAM J. Math. Data Sci. (to appear)}, 2020.

\bibitem{BADHFramesPart2}
B.~Adcock and D.~Huybrechs.
\newblock Frames and numerical approximation {II}: generalized sampling.
\newblock {\em arXiv:1802.01950}, 2018.

\bibitem{BADHframespart}
B.~Adcock and D.~Huybrechs.
\newblock Frames and numerical approximation.
\newblock {\em SIAM Rev.}, 61(3):443--473, 2019.

\bibitem{FEStability}
B.~Adcock, D.~Huybrechs, and J.~Mart{\'\i}n-Vaquero.
\newblock On the numerical stability of {F}ourier extensions.
\newblock {\em Found. Comput. Math.}, 14(4):635--687, 2014.

\bibitem{AdcockNecSamp}
B.~Adcock, R.~Platte, and A.~Shadrin.
\newblock Optimal sampling rates for approximating analytic functions from
  pointwise samples.
\newblock {\em IMA J. Num. Anal.}, 39(3):1360--1390, 2019.

\bibitem{FEParameterSelection}
B.~Adcock and J.~Ruan.
\newblock Parameter selection and numerical approximation properties of
  {F}ourier extensions from fixed data.
\newblock {\em J. Comput. Phys.}, 273:453--471, 2014.

\bibitem{albin2011}
N.~Albin and O.~P. Bruno.
\newblock A spectral {FC} solver for the compressible {N}avier--{S}tokes
  equations in general domains {I}: {E}xplicit time-stepping.
\newblock {\em J. Comput. Phys.}, 230(16):6248--6270, 2011.

\bibitem{boyd2005fourier}
J.~Boyd.
\newblock Fourier embedded domain methods: extending a function defined on an
  irregular region to a rectangle so that the extension is spatially periodic
  and ${C}^{\infty}$.
\newblock {\em Appl. Math. Comput.}, 161(2):591--597, 2005.

\bibitem{BoydFourCont}
J.~P. Boyd.
\newblock A comparison of numerical algorithms for {F}ourier {E}xtension of the
  first, second, and third kinds.
\newblock {\em J. Comput. Phys.}, 178:118--160, 2002.

\bibitem{bruno2010high}
O.~Bruno and M.~Lyon.
\newblock High-order unconditionally stable {FC}-{AD} solvers for general
  smooth domains {I}. {B}asic elements.
\newblock {\em J. Comput. Phys.}, 229(6):2009--2033, 2010.

\bibitem{brunoFEP}
O.~P. Bruno, Y.~Han, and M.~M. Pohlman.
\newblock Accurate, high-order representation of complex three-dimensional
  surfaces via {F}ourier continuation analysis.
\newblock {\em J. Comput. Phys.}, 227(2):1094--1125, 2007.

\bibitem{bunggrieb}
H.-J. Bungartz and M.~Griebel.
\newblock Sparse grids.
\newblock {\em Acta Numer.}, 13:147--269, 2004.

\bibitem{SMSD}
C.~Canuto, M.~Y. Hussaini, A.~Quarteroni, and T.~A. Zang.
\newblock {\em Spectral methods: Fundamentals in {S}ingle {D}omains}.
\newblock Springer, 2006.

\bibitem{ChkifaEtAl}
A.~Chkifa, A.~Cohen, G.~Migliorati, F.~Nobile, and R.~Tempone.
\newblock Discrete least squares polynomial approximation with random
  evaluations - application to parametric and stochastic elliptic {pde}s.
\newblock {\em ESAIM Math. Model. Numer. Anal.}, 49(3):815--837, 2015.

\bibitem{ChkifaEtAlBreaking}
A.~Chkifa, A.~Cohen, and C.~Schwab.
\newblock Breaking the curse of dimensionality in sparse polynomial
  approximation of parametric {PDEs}.
\newblock {\em J. Math. Pures Appl.}, 103:400--428, 2015.

\bibitem{ChkifaDownwardsCS}
A.~Chkifa, N.~Dexter, H.~Tran, and C.~G. Webster.
\newblock Polynomial approximation via compressed sensing of high-dimensional
  functions on lower sets.
\newblock {\em Math. Comp.}, 87:1415--1450, 2018.

\bibitem{christensen2003introduction}
O.~Christensen.
\newblock {\em An Introduction to Frames and {R}iesz Bases}.
\newblock Birkhauser, 2003.

\bibitem{DavenportEtAlLeastSquares}
A.~Cohen, M.~A. Davenport, and D.~Leviatan.
\newblock On the stability and accuracy of least squares approximations.
\newblock {\em Found. Comput. Math.}, 13:819--834, 2013.

\bibitem{CohenDeVoreApproxPDEs}
A.~Cohen and R.~A. DeVore.
\newblock Approximation of high-dimensional parametric {PDE}s.
\newblock {\em Acta Numer.}, 24:1--159, 2015.

\bibitem{CohenDeVoreSchwabFoCM}
A.~Cohen, R.~A. DeVore, and C.~Schwab.
\newblock {Convergence rates of best $N$-term Galerkin approximations for a
  class of elliptic sPDEs}.
\newblock {\em Found. Comput. Math.}, 10:615--646, 2010.

\bibitem{MiglioratiCohenOptimal}
A.~Cohen and G.~Migliorati.
\newblock Optimal weighted least-squares methods.
\newblock {\em SMAI Journal of Computational Mathematics}, 3:181--203, 2017.

\bibitem{ConstantineBook}
P.~G. Constantine.
\newblock {\em Active Subspaces: Emerging Ideas for Dimension Reduction in
  Parameter Studies}.
\newblock SIAM, 2015.

\bibitem{DoostanOwhadiSparse}
A.~Doostan and H.~Owhadi.
\newblock A non-adapted sparse approximation of {PDEs} with stochastic inputs.
\newblock {\em J. Comput. Phys.}, 230(8):3015--3034, 2011.

\bibitem{DynFloaterDownward}
N.~Dyn and M.~Floater.
\newblock Multivariate polynomial interpolation on lower sets.
\newblock {\em J. Approx. Theory}, 177:34--42, 2013.

\bibitem{GanzburgRemez}
M.~I. Ganzburg.
\newblock Polynomial inequalities on measurable sets and their applications.
\newblock {\em Constr. Approx.}, 17:275--306, 2001.

\bibitem{HackbuschTensor}
W.~Hackbusch.
\newblock {$L^\infty$} estimation of tensor truncations.
\newblock {\em Numer. Math.}, 125:419--440, 2013.

\bibitem{PCHansenLCurve}
P.~C. Hansen.
\newblock Analysis of discrete ill-posed problems by means of the l-curve.
\newblock {\em SIAM Rev.}, 34(4):561--580, 1992.

\bibitem{HansenIllPosed}
P.~C. Hansen.
\newblock {\em Rank-Deficient and Discrete Ill-Posed Problems: Numerical
  Aspects of Linear Inversion}.
\newblock SIAM, 1998.

\bibitem{lyon2010high}
M.~Lyon and O.~Bruno.
\newblock High-order unconditionally stable {FC}-{AD} solvers for general
  smooth domains {II}. {E}lliptic, parabolic and hyperbolic {PDE}s; theoretical
  considerations.
\newblock {\em J. Comput. Phys.}, 229(9):3358--3381, 2010.

\bibitem{MiglioratiThesis}
G.~Migliorati.
\newblock {\em Polynomial approximation by means of the random discrete {$L^2$}
  projection and application to inverse problems for {PDEs} with stochastic
  data}.
\newblock PhD thesis, Politecnico di Milano, 2013.

\bibitem{MiglioratiJAT}
G.~Migliorati.
\newblock {Multivariate Markov-type and Nikolskii-type inequalities for
  polynomials associated with downward closed multi-index sets}.
\newblock {\em J. Approx. Theory}, 189:137--159, 2015.

\bibitem{MiglioratiIrregular}
G.~Migliorati.
\newblock Multivariate approximation of functions on irregular domains by
  weighted least-squares methods.
\newblock {\em arXiv:1907.12304}, 2019.

\bibitem{MiglioratiEtAlFoCM}
G.~Migliorati, F.~Nobile, E.~von Schwerin, and R.~Tempone.
\newblock Analysis of the discrete {$L^2$} projection on polynomial spaces with
  random evaluations.
\newblock {\em Found. Comput. Math.}, 14:419--456, 2014.

\bibitem{NarayanJakemanZhouChristoffelLS}
A.~Narayan, J.~D. Jakeman, and T.~Zhou.
\newblock A {C}hristoffel function weighted least squares algorithm for
  collocation approximations.
\newblock {\em Math. Comp.}, 86:1913--1947, 2017.

\bibitem{NarayanZhouCCP}
A.~Narayan and T.~Zhou.
\newblock Stochastic collocation on unstructured multivariate meshes.
\newblock {\em Commun. Comput. Phys.}, 18(1):1--36, 2015.

\bibitem{NeumaierIllCond}
A.~Neumaier.
\newblock Solving ill-conditioned and singular linear systems: a tutorial on
  regularization.
\newblock {\em SIAM Rev.}, 40(3):636--666, 1998.

\bibitem{pasquettiFourEmbed}
R.~Pasquetti and M.~Elghaoui.
\newblock A spectral embedding method applied to the advection--diffusion
  equation.
\newblock {\em J. Comput. Phys.}, 125:464--476, 1996.

\bibitem{PengHamptonDoostantweighted}
J.~Peng, J.~Hampton, and A.~Doostan.
\newblock A weighted {$\ell_1$}-minimization approach for sparse polynomial
  chaos expansions.
\newblock {\em J. Comput. Phys.}, 267:92--111, 2014.

\bibitem{TrefPlatteIllCond}
R.~Platte, L.~N. Trefethen, and A.~Kuijlaars.
\newblock Impossibility of fast stable approximation of analytic functions from
  equispaced samples.
\newblock {\em SIAM Rev.}, 53(2):308--318, 2011.

\bibitem{SargsyanEtAlUQDisc}
K.~Sargsyan, C.~Safta, B.~J. Debusschere, and H.~Najm.
\newblock Uncertainty quantification given discontinuous model response and a
  limited number of model runs.
\newblock {\em {SIAM} J. Sci. Comput.}, 34(1):B44--B64, 2012.

\bibitem{SargsyanEtAlDimUQ}
K.~Sargsyan, C.~Safta, H.~Najm, B.~J. Debusschere, D.~Ricciuto, and
  P.~Thornton.
\newblock Dimensionality reduction for complex models via {B}ayesian
  compressive sensing.
\newblock {\em Int. J. Uncertain. Quantif.}, 4(1):63--93, 2014.

\bibitem{StinsonZonotop}
K.~Stinson, D.~F. Gleich, and P.~G. Constantine.
\newblock A randomized algorithm for enumerating zonotope vertices.
\newblock {\em arXiv:1602.06620}, 2016.

\bibitem{TemlyakovRemezHC}
V.~Temlyakov and S.~Tikhonov.
\newblock Remez-type inequalities for the hyperbolic cross polynomials.
\newblock {\em Constr. Approx.}, 46(3):593--615, 2017.

\bibitem{TroppUserFriendly}
J.~A. Tropp.
\newblock User friendly tail bounds for sums of random matrices.
\newblock {\em Found. Comput. Math.}, 12:389--434, 2012.

\bibitem{WitteveenIaccarinoSimplex}
J.~A.~S. Witteveen and G.~Iaccarino.
\newblock Simplex stochastic collocation with random sampling and extrapolation
  for nonhypercube probability spaces.
\newblock {\em {SIAM} J. Sci. Comput.}, 34(2):A814--A838, 2012.

\bibitem{XiuKarniadakisPC}
D.~Xiu and G.~E. Karniadakis.
\newblock {The Wiener--Askey polynomial chaos for stochastic differential
  equations}.
\newblock {\em {SIAM} J. Sci. Comput.}, 24(2):619--644, 2002.

\bibitem{YanGuoXui_l1UQ}
L.~Yan, L.~Guo, and D.~Xiu.
\newblock Stochastic collocation algorithms using {$\ell_1$}-minimization.
\newblock {\em Int. J. Uncertain. Quantif.}, 2(3):279--293, 2012.

\bibitem{KarniadakisUQCS}
X.~Yang and G.~E. Karniadakis.
\newblock Reweighted {$\ell_1$} minimization method for stochastic elliptic
  differential equations.
\newblock {\em J. Comput. Phys.}, 248:87--108, 2013.

\bibitem{ZhouNarayanXiu}
T.~Zhou, A.~Narayan, and D.~Xiu.
\newblock Weighted discrete least-squares polynomial approximation using
  randomized quadratures.
\newblock {\em J. Comput. Phys.}, 1:787--800, 2015.

\end{thebibliography}

\end{document}